%% file: OptINNs.tex
\documentclass[]{article}
\title{Optimality-Informed Neural Networks for Solving Parametric Optimization Problems}
\newcommand{\institute}{
              Saarland University, \\
              Saarbrücken,\\
              Germany
}

\author{
	Matthias K. Hoffmann\\ \institute \\ 
              \texttt{matthias.hoffmann@uni-saarland.de}
	\And 
	Amine Othmane \\ \institute
	\And 
	Kathrin Fla\ss{}kamp \\ \institute
}

\input{general_packages}
\usepackage{standalone}
\input{macros}
\input{listoftodos}

\usepackage[outline]{contour} 
\contourlength{1.4pt}

\graphicspath{{figs}}

\input{inc/acronyms}
\newcolumntype{Y}{>{\centering\arraybackslash}X}
\begin{document}

\maketitle

\begin{abstract}\noindent
Many engineering tasks require solving families of nonlinear constrained optimization problems, parametrized in setting-specific variables.
This is computationally demanding, particularly if solutions have to be computed across strongly varying parameter values, e.g., in real-time control or for model-based design.
Thus, we propose to learn the mapping from parameters to the primal optimal solutions and to their corresponding duals using neural networks, giving a dense estimation in contrast to gridded approaches.
Our approach, Optimality-Informed Neural Networks (OptINNs), combines (i) a KKT-residual loss that penalizes violations of the first-order optimality conditions under standard constraint qualification assumptions, and (ii) problem-specific output activations that enforce simple inequality constraints (e.g., box-type/positivity) by construction.
This design reduces data requirements, allows the prediction of dual variables, and improves feasibility and closeness to optimality compared to penalty-only training.
Taking quadratic penalties as a baseline, as previously proposed for this problem class in the literature, our method simplifies hyperparameter tuning and attains tighter adherence to optimality conditions.
We evaluate OptINNs on different nonlinear optimization problems ranging from low to high dimensions and compare them with equally sized neural networks trained using the quadratic-penalty method.
On small problems, OptINNs match a quadratic-penalty baseline in primal accuracy while additionally predicting dual variables with low error.
On larger problems, OptINNs achieve lower constraint violations, lower primal error, and lower errors in the predicted costs.
These results suggest that embedding feasibility and optimality into the network architecture and loss can make learning-based surrogates more accurate, feasible, and data-efficient for parametric optimization.
\end{abstract}



\section*{Acknowledgments}
M.~K.~Hoffmann acknowledges funding by the German Research Foundation (DFG), Project No.~501928699 and 564445282.
A.~Othmane is grateful for the support from the German Federal Ministry of Research, Technology and Space, Project No.~05M22TSA.

\input{inc/introduction_new}
\input{inc/p_optim_problems}

\input{inc/optinns}
\input{inc/comparison_pmnn.tex}
\input{inc/numerical_results}
\input{inc/summary}

\newpage
\printbibliography[]
\newpage
\section*{Appendix}
\appendix
\input{inc/appendix}
\end{document}

%% file: general_packages.tex
\usepackage{arxiv}

\usepackage{amsmath}
\usepackage{amsthm}
\usepackage{amsfonts}
\usepackage{amssymb}
\usepackage[nolist]{acronym}
\usepackage[utf8]{inputenc}
\usepackage{microtype}
\usepackage[T1]{fontenc}

\usepackage{booktabs}

\usepackage{listings}
\usepackage{matlab-prettifier}
\usepackage{multicol}
\usepackage[labelfont=bf, labelsep=space]{caption}
\usepackage{subcaption}
\usepackage{enumerate}
\usepackage{graphicx}  
\usepackage{tabularx}
\usepackage{xcolor}
\usepackage{float}
\usepackage{csquotes}

\makeatletter
\let\blx@noerroretextools\relax
\makeatother
\usepackage[style=authoryear, 
			natbib=true,
			maxbibnames=3, 
			minbibnames=1, 
			backend=bibtex, 
			url=false, 
			doi=true, 
			isbn=false,
			eprint=false,
			date=year,
			dashed=false,
			]{biblatex}
\renewcommand*{\cite}{\parencite}
\DeclareDelimFormat{nameyeardelim}{\addspace}

\AtEveryBibitem{%
  \clearfield{note}%
}
\usepackage{siunitx}
\usepackage{xspace}

\usepackage[ruled]{algorithm2e}

\usepackage{mathtools}

\usepackage{multirow}
\usepackage[hidelinks]{hyperref} 
\makeatletter
\def\cl@chapter{}
\makeatother
\usepackage[capitalise,nameinlink]{cleveref}
\crefname{section}{section}{sections}
\Crefname{section}{Section}{Sections}

\usepackage{eurosym} 

\usepackage{nicefrac}

\usepackage{tikz}
\usepackage{pgfplots}
\usepackage{url}

\usepackage{autonum}

\usetikzlibrary{arrows,shapes,positioning,calc,intersections,fit,decorations.markings, arrows.meta, backgrounds, spy, arrows.meta, bending, backgrounds, patterns}
\usepackage{adjustbox}
\usepgfplotslibrary{groupplots, patchplots, statistics, fillbetween}
\usepackage{pgfplotstable}

\pgfplotsset{compat=newest}

\usepackage{ifplatform}
\usepackage{ifthen}

\tikzset{
	partial ellipse/.style args={#1:#2:#3}{
		insert path={+ (#1:#3) arc (#1:#2:#3)}
	}
}

\usepackage{listofitems}

\iffalse 
\usetikzlibrary{external}
\else
\PackageWarning{}{Shell escape is not enabled. Externalization won't work.}
\providecommand{\tikzsetnextfilename}[1]{}
\fi

\definecolor{colorinnertube}{RGB}{0,72,119}
\definecolor{coloroutertube}{RGB}{200,34,84}
\definecolor{udsblue}{RGB}{0,72,119}
\definecolor{udsred}{RGB}{200,34,84}
\definecolor{udsgreen}{RGB}{10,102,10}
\definecolor{udsyellow}{RGB}{215,223,35}
\definecolor{udsorange}{HTML}{ffae00}

\newtheorem{assumption}{Assumption}



%% file: macros.tex
\DeclareMathOperator*{\relu}{ReLU}

\DeclareMathOperator*{\rk}{RK4}

\newcommand{\st}{\ensuremath{\text{s.\ t.}}}
\definecolor{matlabyellow}{RGB}{237,177,32}

\newcommand{\gracemath}[1]{\ensuremath{#1}\xspace}

\newcommand{\R}{\ensuremath{\mathbb{R}}\xspace}
\newcommand{\N}{\ensuremath{\mathcal{N}}\xspace}

\newcommand{\nameStationarity}{Stat}
\newcommand{\namePrimalG}{FeasG}
\newcommand{\namePrimalH}{FeasH}
\newcommand{\nameComplSlack}{CSl}

\newcommand{\calL}{\mathcal{L}}
\newcommand{\calD}{\mathcal{D}}
\newcommand{\namedL}[1]{\calL^\mathrm{#1}}
\newcommand{\namedD}[1]{\calD^\mathrm{#1}}
\newcommand{\namedo}[1]{\omega^\mathrm{#1}}
\newcommand{\Ls}{\gracemath{\namedL{\nameStationarity}}}
\newcommand{\Lpg}{\gracemath{\namedL{\namePrimalG}}}
\newcommand{\Lph}{\gracemath{\namedL{\namePrimalH}}}
\newcommand{\Lcs}{\gracemath{\namedL{\nameComplSlack}}}
\newcommand{\Ds}{\gracemath{\namedD{\nameStationarity}}}
\newcommand{\Dpg}{\gracemath{\namedD{\namePrimalG}}}
\newcommand{\Dph}{\gracemath{\namedD{\namePrimalH}}}
\newcommand{\Dcs}{\gracemath{\namedD{\nameComplSlack}}}

\newcommand{\Lkkt}{\gracemath{\namedL{KKT}}}
\newcommand{\Dmse}{\gracemath{\namedD{MSE}}}
\newcommand{\Lmse}{\gracemath{\namedL{MSE}}}
\newcommand{\os}{\gracemath{\namedo{\nameStationarity}}}
\newcommand{\opg}{\gracemath{\namedo{\namePrimalG}}}
\newcommand{\oph}{\gracemath{\namedo{\namePrimalH}}}
\newcommand{\ocs}{\gracemath{\namedo{\nameComplSlack}}}

\newcommand{\namedP}[1]{P^\mathrm{#1}}
\newcommand{\Ps}{\gracemath{\namedP{\nameStationarity}}}
\newcommand{\Ppg}{\gracemath{\namedP{\namePrimalG}}}
\newcommand{\Pph}{\gracemath{\namedP{\namePrimalH}}}
\newcommand{\Pcs}{\gracemath{\namedP{\nameComplSlack}}}

\newcommand{\nn}{\gracemath{\mathcal{N}}}
\newcommand{\nnt}[2][{}]{\gracemath{\nn_{#1}(#2,\tv)}}

\newcommand{\bound}[1]{\gracemath{\overline{\underline{#1}}}}

\newcommand{\makev}[1]{\boldsymbol{#1}}
\newcommand{\xv}{{\makev{x}}}
\newcommand{\yv}{{\makev{y}}}
\newcommand{\tv}{{\makev{\theta}}}
\newcommand{\lv}{{\makev{\lambda}}}
\newcommand{\mv}{{\makev{\mu}}}
\newcommand{\pv}{{\makev{p}}}
\newcommand{\gv}{{\makev{g}}}
\newcommand{\hv}{{\makev{h}}}

\newtheorem{theorem}{Theorem}[section]
\newtheorem{definition}{Definition}[section]
\newtheorem{lemma}{Lemma}[section]
\newtheorem{corollary}{Corollary}[section]

\newtheorem{remark}{Remark}[section]
\newtheorem{example}{Example}[section]

\newcommand{\Rset}[1]{\R^{#1}}
\newcommand{\setx}{\Rset{n_x}}
\newcommand{\setp}{\Rset{n_p}}
\newcommand{\setg}{\Rset{n_g}}
\newcommand{\seth}{\Rset{n_h}}


%% file: listoftodos.tex
\usepackage[colorinlistoftodos,textsize=tiny]{todonotes} 
\usepackage{ifthen} 

\newboolean{correctionsApproved} 
\setboolean{correctionsApproved}{true} 
\newboolean{remarksApproved} 
\setboolean{remarksApproved}{true} 
\newboolean{marknewparts}
\setboolean{marknewparts}{false}

\ifthenelse{\boolean{marknewparts}}{
\newcommand{\colorKFcorrection}{blue} 
\newcommand{\colorKFremark}{blue!80} 
\newcommand{\colorMKHcorrection}{udsred} 
\newcommand{\colorMKHremark}{udsred!80} 
\newcommand{\colorAOcorrection}{udsorange} 
\newcommand{\colorAOremark}{udsorange!80} 
}
{
\newcommand{\colorKFcorrection}{red} 
\newcommand{\colorKFremark}{red!80} 
\newcommand{\colorMKHcorrection}{red} 
\newcommand{\colorMKHremark}{red!80} 
\newcommand{\colorAOcorrection}{red} 
\newcommand{\colorAOremark}{red!80} 
}

\ifthenelse{\boolean{correctionsApproved} \and \boolean{remarksApproved}} 
{}
{} 

\makeatletter
\newcommand{\correctionstart}{%
	\ifhmode
		\ifdim\lastskip>\z@ \let\correction@next\ignorespaces
		\else \let\correction@next\relax
		\fi
	\else \let\correction@next\ignorespaces
	\fi
	\correction@next}
\makeatother

\newcommand{\KFcorrection}[2][]{%
	\ifthenelse{\boolean{correctionsApproved}}%
	{\correctionstart#2\ignorespaces\xspace}
	{%
		{\color{\colorKFcorrection}{#2}}
		\addcontentsline{tdo}{todo}{\color{\colorKFcorrection}{KF correction:} #1\xspace}
	}%
}

\newcommand{\KFremark}[2][]{%
	\ifthenelse{\boolean{remarksApproved}}%
	{~\ignorespaces}
	{%
		\footnote{\color{\colorKFremark}{#2}}
		\addcontentsline{tdo}{todo}{\color{\colorKFremark}{KF remark:} #1}
	}%
}

\newcommand{\MKHcorrection}[2][]{%
	\ifthenelse{\boolean{correctionsApproved}}%
	{\correctionstart#2\ignorespaces\xspace}
	{%
		{\color{\colorMKHcorrection}{#2}}
		\addcontentsline{tdo}{todo}{\color{red}{MKH correction:} #1\xspace}
	}%
}

\newcommand{\MKHremark}[2][]{%
	\ifthenelse{\boolean{remarksApproved}}%
	{~\ignorespaces}
	{%
		\footnote{\color{\colorMKHremark}{#2}}
		\addcontentsline{tdo}{todo}{\color{\colorMKHremark}{MKH remark:} #1}
	}%
}

\newcommand{\AOcorrection}[2][]{%
	\ifthenelse{\boolean{correctionsApproved}}%
	{\correctionstart#2\ignorespaces\xspace}
	{%
		{\color{\colorAOcorrection}{#2}}
		\addcontentsline{tdo}{todo}{\color{red}{AO correction:} #1\xspace}
	}%
}

\newcommand{\AOremark}[2][]{%
	\ifthenelse{\boolean{remarksApproved}}%
	{~\ignorespaces}
	{%
		\footnote{\color{\colorAOremark}{#2}}
		\addcontentsline{tdo}{todo}{\color{\colorAOremark}{AO remark:} #1}
	}%
}

\let\oldlistoftodos\listoftodos
\renewcommand{\listoftodos}{%
	\ifthenelse{\boolean{correctionsApproved} \and \boolean{remarksApproved}}%
	{}
	{\oldlistoftodos}
}

%% file: inc/acronyms.tex
\begin{acronym}
	\acro{nn}[NN]{neural network}
	\acro{mlp}[MLP]{multi-layer perceptron}
	\acro{licq}[LICQ]{linear independence constraint qualification}
	\acro{pinn}[PINN]{Physics-Informed Neural Network}
\end{acronym}

%% file: inc/introduction_new.tex

\section{Introduction}
Many optimization problems which arise in engineering applications are parametric, i.e., the cost criteria as well as the constraints depend on external parameters.
In multi-objective optimization, for instance, any scalarization technique introduces a family of optimization problems, which is parametrized via weights or alternative hyperparameters modeling prioritization among multiple objectives~\cite{Marler2004}.
Parametric optimization problems also show up in model predictive control, in which subsequent optimal control problems vary in their current initial state, which contributes to the set of constraints~\cite{Bemporad2002explicitmpc}.
In general, many types of external influences from varying environmental conditions can be modeled as a vector of parameters.
These enter the constraints and the objective of an optimization problem as presented by works from economics~\cite{Markowitz1952portfolio} to engineering~\cite{martins2021engineering}.
A historical review of parametric optimization can be found in~\cite{Oberdieck2016pprogramming}.
For families of nonlinear constrained optimization problems, one is interested in finding not just the solutions for a small number of parameter values, but the functional relationship between the problem instances and their respective optimal solutions.
Under regularity assumptions, classical results guarantee that this relationship is, at least locally, a smooth curve~\cite{fiacco1983sensitivity}.
Still, its computation raises numerical challenges, since a potentially high-dimensional, multivariate optimization problem cannot be solved to high accuracy on a dense parameter grid.
Instead, a suitable representation for the curve of optimal solutions has to be found.
To this end, we propose a neural network (NN) representation and ensure approximation quality by a combined architecture and loss function design.
This approach provides a parameter-to-optimal-solution mapping, which allows efficient evaluation as required in, e.g., optimization-based computer-aided design, decision-making in multi-objective optimization, or explicit model predictive control.
Specifically, it provides a mapping for optimal primal-dual pairs, which are necessary for, e.g., sensitivity analysis~\cite{rao1989sensitivity}.
In safety-critical applications, where approximation errors of the NN model are intolerable, the solution can still serve as an initial guess for warm-starting any gradient-based optimization method.

\subsection{Related Work}\label{sec:rel_works_nn}
Solution maps for parametric optimization problems have classically been computed by pathfollowing and surrogate-based techniques.
More recently, intensive research has been carried out on using neural-network-based methods instead.

\paragraph{Pathfollowing methods}
Pathfollowing methods seek to update the solution when the parameter vector varies by computing a new solution from a given solution tuple via continuation~\cite{guddat1990jumps}.
Continuation methods rely on the smoothness of the solution with respect to the parameters, as they use the implicit function theorem to predict the evolution of the optimal solution.
For application to parametric optimization, the conditions of optimality are formulated as an implicit relation, allowing one to track the corresponding solution curve to compute a new estimate~\cite{rao1989continuation}.
A common approach is a predictor-corrector algorithm, which follows the tangent of the current solution and then re-enters the solution curve via Newton's method.
Detailed insights into the realization of this tracking are provided in~\cite{allgower1990continuation}.
While the continuation path is uniquely defined in one-dimensional parameter spaces, the lack of a natural continuation direction complicates the tracking process in higher dimensions~\cite{rheinboldt1988manifold}.
Such multiparametric optimization problems arise in optimal control problems with parametric initial conditions or multi-objective optimization problems with three or more objectives.
These can be solved by pre-defining a path parametrized in arc-length to reduce the dimensionality to one, or by tracking the full manifold~\cite{rheinboldt1988manifold}.
These methods, however, can be computationally intensive, which is a limitation for real-time applications, and give no dense representation of the solution curve.

\paragraph{Surrogate-based methods}
Distinct from pathfollowing, surrogate-based methods aim to find surrogate models for the mapping from parameters to the optimal solutions.
In the context of optimal control with linear-quadratic cost and linear dynamics, the solution can be represented by a continuous, piecewise affine mapping as a function of the initial condition~\cite{Bemporad2002explicitmpc}.
Notably, only one \MKHcorrection{solution of the optimization problem per set of active inequality constraints} is sufficient to construct this mapping.
In general, however, the structure of the solution mapping is unknown, creating a need for more versatile approximation methods.
In~\cite{demarchi2023functionapprox}, this broader scope is addressed by inferring the solution mapping without assuming a specific functional form.
This motivates their use of surrogate models fitted via residual minimization of the KKT conditions.
Minimizing residuals is analogous to performing root-finding on the KKT conditions, as it is done \MKHcorrection{for modified KKT conditions} in Primal-Dual Interior Point methods~\cite{boyd2004convex}.
In \cite{demarchi2023functionapprox}, the authors propose fitting a linear combination of Radial Basis Functions (RBFs) without using data samples by evaluating the residuals on a fixed grid of parameter values.
These two approaches represent contrasting paradigms: the former leverages specific structural properties and \MKHcorrection{solutions of the optimization problem} for control tasks, while the latter offers a general, grid-based framework that operates purely on residuals without incorporating solution data.
However, the reliance on a predefined sample grid for the RBF limits the scalability of the latter approach.
Specifically, RBF methods are more susceptible to the curse of dimensionality than neural networks, as the number of basis functions required to cover the space grows exponentially with the dimension~\cite{Bengio2005curse_rbf, Barron1993curse}.

\paragraph{Neural-network-based methods for parametric optimization problems}
Neural networks have seen significant advancements over the past two decades, with growing impact in engineering applications. Their ability to model complex, nonlinear relationships such as dynamic system behavior~\cite{chen2018node} and to process large datasets makes them well-suited for tasks like predictive maintenance~\cite{rivas2020predictive}.
However, their requirement for large amounts of training data has long been a limiting factor, also when applying them to parametric optimization problems.

Early investigations primarily concentrated on developing recurrent neural networks or general differential equations designed to converge to local optima.
In 1985, Hopfield and Tank demonstrated that Hopfield networks could minimize quadratic energy functions to address optimization problems, effectively estimating low-cost solutions to the Traveling Salesman Problem~\cite{hopfield1985neural}.
Building on this, Dhingra and Rao introduced mechanisms to handle constraints via penalty terms, approximating solutions for constrained optimization problems in mechanical engineering~\cite{rao1992nnmechanical}.
Similar approaches were explored in~\cite{Lillo1993penaltynn}.
Zhang and Constantinides introduced Lagrange Programming Neural Networks by integrating the Karush-Kuhn-Tucker (KKT) conditions to address necessary conditions of optimality~\cite{zhang1992lnn}.
Their resulting ordinary differential equation applied gradient descent on the primal and gradient ascent on the dual variables simultaneously by applying numerical integration methods.
Wu and Tam specialized these methodologies for quadratic optimization problems~\cite{wu1999nnquadratic}.

More recently, \acp{mlp} have become the predominant architecture for neural-network-based optimization.
Nikbakht et al.\ proposed minimizing the Lagrangian directly by approximating the Lagrange multipliers through cross-validation~\cite{nikbakht2020unsupervised}.
\MKHcorrection{Liu et al.\ trained neural networks that map the problem parameters to the optimal solution in a single forward pass, either on solver data or on the cost function with the constraints as quadratic penalty terms~\cite{liu2024teaching}.}
Although quadratic-penalty methods are widely used, they are often considered suboptimal in modern optimization theory due to their tendency to produce solutions outside the feasible space~\cite{Lange2013,nocedal2006numericaloptimization}.
By relying on constraint-penalizing losses, these approaches do not fully exploit fundamental principles from optimization theory.
In particular, they cannot provide Lagrange multipliers, which are helpful, e.g., as initial guesses for classical solvers, or even required, such as for sensitivity analysis~\cite{rao1989sensitivity}.

\MKHcorrection{Instead of penalizing constraint violations, other approaches include correction or even optimization steps in the network itself.
Differentiable optimization layers include the solution of a convex optimization problem, whose parameters are given by the outputs of the preceding layers, as a layer of the network~\cite{amos2017optnet,agrawal2019cvxpylayers}.
The problem is solved numerically in the forward pass, and its solution is differentiated with respect to the problem's parameters in the backward pass.
For general, possibly nonconvex problems, \MKHcorrection{correction- and projection-based} methods correct the prediction of the network towards the feasible set.
The method DC3 predicts a subset of the decision variables, computes the remaining ones by solving the equality constraints, and reduces violations of the inequality constraints by gradient descent, both in training and at inference~\cite{donti2021dc3}.
In NLPOpt-Nets, layers that minimize a quadratic approximation of the cost function subject to the constraints approximately, by a sequence of quadratic programs with linearized constraints, are appended.
These layers additionally predict Lagrange multipliers~\cite{roy2026nlpoptnet}.
In Deep ALM, the network predicts only the Lagrange multipliers and is trained to maximize the dual function, while the decision variables are computed by minimizing the augmented Lagrangian subject to box constraints~\cite{kotary2024deepalm}.
All these methods require iterative computations in each evaluation of the network, also at inference.}

\paragraph{Theory-informed neural networks}
To overcome the issue of requiring a large amount of data in the training of NNs, theory-informed approaches have been developed.
A foundational example was proposed by Munos, Baird, and Moore, who used NNs to approximate solutions to the Hamilton-Jacobi-Bellman (HJB) equation~\cite{munos1999hjbresidual}.
For optimal control, this equation describes the necessary and sufficient conditions of optimality, and training is performed by minimizing its residuals.
Similar approaches for solving the HJB equation appear in other works~\cite{tao2005hjb, dario2019hjbinterplanetary, sirignano2018dgm, verma2024nn_pocp}.
Complementary to HJB-based methods, Effati and Pakdaman incorporate the Pontryagin Maximum Principle (PMP) to provide necessary conditions for controlling a dynamical system~\cite{effati2013ocp}.
\MKHcorrection{Seo demonstrated the use of \acp{pinn} for optimal control problems by embedding the conditions of the PMP into the loss as squared residuals~\cite{seo2024neuralcomputing}.}
Recent work by Nakamura-Zimmerer et al.\ exemplifies hybrid approaches, \MKHcorrection{training on PMP-based data of the value function and its gradient instead of on the residuals of the HJB equation}~\cite{nakamura-zimmerer_adaptive_2021}.
In this framework, we see an opportunity to leverage the Karush-Kuhn-Tucker (KKT) conditions, which provide necessary conditions for optimality in constrained optimization problems.
\MKHcorrection{For parametric convex optimization problems, for which, under Slater's condition, the KKT conditions are necessary and sufficient for global optimality, losses based on the KKT conditions have also been proposed in~\cite{femine2024KKTinformed,arving2024kktnets}.
In~\cite{femine2024KKTinformed}, the network is trained without training data on the deviations from the KKT conditions, which are minimized jointly by Jacobian descent, and a ReLU function ensures non-negative Lagrange multipliers.
In~\cite{arving2024kktnets}, a weighted sum of a KKT-loss and an MSE-loss with fixed weights is minimized for linear programs with two decision variables, and the predicted primal and dual variables are reported to be not yet sufficiently close to the optimal solutions.}

\subsection{Contributions of this paper}
Addressing the challenges of the aforementioned classical, surrogate-based and neural-network-based methods for parametric optimization problems, we propose \emph{Optimality-Informed Neural Networks} (OptINNs).
In this framework, optimality conditions are directly embedded into the training process and simple constraints integrated in the design of the network's architecture.
This work is situated within the growing field of theory-informed machine learning, building on approaches such as \acp{pinn} for solving ordinary or partial differential equations (ODEs/PDEs), see~\cite{raissi2019PINN}, or geometry-informed neural networks, see~\cite{berzins2024ginn}.
Similar to other surrogate-based approaches, neural networks allow for a dense representation of the solution space, while suffering less from the curse of dimensionality than kernel-based approaches~\cite{Bengio2005curse_rbf, Barron1993curse}.
In contrast to methods based on the penalty method, incorporating optimality in the training allows for the estimation of the corresponding Lagrange multipliers.
\MKHcorrection{In contrast to differentiable optimization layers, projection-based methods, and Deep ALM~\cite{amos2017optnet,donti2021dc3,roy2026nlpoptnet,kotary2024deepalm}, a trained OptINN does not embed an iterative algorithm, such that estimates of the optimal solution and the associated Lagrange multipliers are obtained in a single forward pass.}

For OptINNs, we propose a scheme that integrates data from optimal solutions with evaluations of optimality residuals, with respect to the KKT conditions, at sampled parameter instances. 
The scheme remains applicable even when no data from optimal solutions are available.
Thereby, the method balances between a pure data-based training and the approach suggested in \cite{demarchi2023functionapprox} for approximating the solution mapping of parametrized optimization problems.
By framing our approach within the theory-informed machine learning paradigm, we benefit from utilizing established training strategies~\cite{wang_experts_2023, krishnapriyan2021pinnfailure}.
We derive advantages of OptINNs in comparison to the established quadratic-penalty method: OptINNs do not show the known shortcoming of quadratic-penalty methods, which suffer from offsets in optimal solutions due to the designed constraint penalization.
Unlike quadratic-penalty-method-based neural networks, OptINNs predict both the primal and dual variables of the patched solution. 
Despite this added task, numerical experiments demonstrate that same-sized OptINNs show little to no deterioration in performance or accuracy in small optimization problems as well as better performance and accuracy in larger optimal control problems.

\subsection{Outline}
The rest of this paper is structured as follows.
In \Cref{sec:popt}, the foundations of parametric nonlinear optimization, including the conditions of optimality, the local existence of solution branches, and the construction of patched solutions, are laid.
In \Cref{sec:optinn}, we present our main contributions. 
Namely, we discuss the design of loss functions incorporating optimality conditions and data fitting. Further, we discuss the architecture of OptINNs and their training, in particular the balancing of different components in the loss function.
Then, in \Cref{sec:constraint-penalization}, we briefly introduce neural network training via constraint-penalization and discuss its disadvantages compared to the OptINN approach.
\Cref{sec:numerical_results} presents numerical experiments on four optimization problems, for which we also compare OptINNs with quadratic-penalty-method-based neural networks.
Our focus lies on low-data regimes, in particular.
Finally, \Cref{sec:sum} summarizes our findings and outlines directions for future research.

\subsection{Notation}
In this work, vectors, matrices and vector-valued functions are denoted with bold symbols.
We denote the $i$-th element of a vector $\boldsymbol{v} \in \R^n$ as $\boldsymbol{v}^{i}$. Similarly, for a set of indices $\mathbb{I}\subseteq \{1, 2, \dots, n\}$ the corresponding entries are represented as $\boldsymbol{v}^\mathbb{I}$.
Variables, transformed to fulfill bound and box constraints by design, are denoted with over- and underline \bound{x}. 
In optimization problem formulations, comparison operators (e.g., $\leq$) applied to vectors and scalars are interpreted as element-wise operations.
For instance, $\boldsymbol{v} \le \yv$ for $\boldsymbol{v},\yv\in\mathbb{R}^n$ implies $\boldsymbol{v}^i \le \yv^i$ for all $i=1,\dots,n$.
With $\mathcal{C}^k(\R^a, \R^b)$ we denote the set of all $k$-times continuously differentiable, vector-valued functions from $\R^a$ to $\R^b$, with $a,b \in \mathbb{N}_+$.
The derivative of a scalar-valued function $f$ with scalar argument is denoted with $f'$. 
The gradient w.r.t.\ $\xv$ of a scalar-valued function with vectorial argument $f: \xv \mapsto f(\xv)$ is denoted with $\nabla_x f(\xv)$.
Similarly, the Jacobian matrix w.r.t.\ $\xv$ is denoted as $J_x f(\xv)$ and the partial derivative w.r.t.\ the $i$-th entry of $\xv$ is denoted as $\frac{\partial f}{\partial \xv^i}$.

%% file: inc/p_optim_problems.tex
\section{Parametric nonlinear optimization problems}\label{sec:popt}
In this section, we present the foundations for the upcoming approximation of patched solutions of parametric nonlinear optimization problems.
For this, we first summarize optimization theory, mainly adapted from~\cite{nocedal2006numericaloptimization}.
Afterwards, we present the class of optimization problems we handle in this work.

We consider a nonlinear optimization problem of the form
\begin{equation}
	\begin{aligned}
		\min_{\xv\in\R^n} \ &f(\xv, \pv)\\
		\text{s.t.}\ &\gv(\xv, \pv) \le 0,\\
		&\hv(\xv, \pv) = 0,
	\end{aligned}\tag*{OP($\pv$)}\label{eq:pOP}
\end{equation}
where \(\xv\in\setx\) is the vector of decision variables, \(f \in \mathcal{C}^2(\setx\times\setp, \R)\) is the cost function to be minimized, and \(\gv\in \mathcal{C}^2(\setx\times\setp,\setg)\) and \(\hv\in \mathcal{C}^2(\setx\times\setp,\seth)\) are the inequality and equality constraints, respectively.
Each of these functions additionally depends on the parameter vector \(\pv\in\setp\), making it a parametric optimization problem.
This parameter vector contains all entities that can be specified in an exact problem setting, like values of physical constants, the weights in a weighted sum of cost functions, or initial conditions in optimal control.
For each choice of \(\pv\in\setp\), \ref{eq:pOP} states a nonlinear constrained optimization problem.
Let $\xv^\star$ denote a local minimizer of \ref{eq:pOP}.

Among the different variants of \emph{constraint qualifications}, we assume the linear independence constraint qualification (LICQ). 
It requires that, for fixed $\pv$, the gradients of the active inequality constraints, i.e.\ for \( \mathbb{A} \coloneqq \bigl\{ i  \mid {\gv}^i(\xv^\star,\pv)=0\bigr\}\), and the gradients of the equality constraints,
\[ \mathbb{E} = 
\bigl\{ \nabla_{x} \gv^{i}({\xv}^\star,\pv) \mid i\in \mathbb{A} \bigr\} \,\cup\,\bigl\{\nabla_x {\hv}^j(\xv^\star,\pv) \mid j=1,\dots,n_h \bigr\}
\]
are linearly independent.
This means that with the cardinality of $\mathbb{E}$ denoted by $|\mathbb{E}|$,
\begin{align}\label{eq:licq}
	\nexists \, \boldsymbol{\xi}\in\R^{|\mathbb{E}|}\setminus\{0\},\ \text{such that} \ \sum_{i=1}^{|\mathbb{E}|} \boldsymbol{\xi}^i\mathbb{E}^i(\xv^\star, \pv) = 0.
\end{align}
That is, the linearized feasible directions characterize the tangent cone at $\xv^\star$, see \cite{nocedal2006numericaloptimization}.

With Lagrange multipliers \(\mv \in \setg\) and \(\lv\in\seth\), we define the parametric Lagrangian
\[
	L:\setx\times \seth \times \setg \times \setp \to \mathbb{R}
\]
associated with~\ref{eq:pOP} by \(L(\xv, \lv, \mv, \pv) = f(\xv, \pv) + \lv^\top \hv(\xv, \pv) + \mv^\top \gv(\xv, \pv) \).
For a local minimizer $\xv^\star$ of~\ref{eq:pOP}, with associated Lagrange multipliers $\lv^\star, \mv^\star$, the necessary conditions for optimality are given by the KKT conditions
\begin{equation}
	\begin{aligned}
		&\text{\makebox[4.5cm][l]{{Stationarity}:}}
		&&\nabla_x L(\xv^\star, \lv^\star, \mv^\star, \pv) = 0 && \text{\makebox[4.5cm][l]{(\nameStationarity)}}\\
		&\text{\makebox[4.5cm][l]{{Primal feasibility}:}} 
		&&\gv(\xv^\star, \pv) \le 0 && \text{\makebox[4.5cm][l]{(\namePrimalG)}}\\
		& 
		&&\hv(\xv^\star, \pv) = 0 && \text{\makebox[4.5cm][l]{(\namePrimalH)}}\\
		&\text{\makebox[4.5cm][l]{{Dual feasibility}:}} 
		&&{\mv^\star}^i \ge 0 \quad \forall i \in \{1,\dots,n_g\}\\
		&\text{\makebox[4.5cm][l]{{Complementary slackness}:}} 
		&&{\mv^\star}^i {\gv}^i(\xv^\star, \pv) = 0 \quad \forall i \in \{1,\dots,n_g\}. && \text{\makebox[4.5cm][l]{(\nameComplSlack)}}
	\end{aligned}
	\tag*{KKT($\pv$)}\label{eq:kkt}
\end{equation}
Under LICQ, a local minimizer $\xv^\star$ satisfies the KKT conditions.
To ensure that a candidate solution is a local minimizer, second-order sufficient conditions (SOSC) must be satisfied. Specifically, the Hessian of the Lagrangian $\nabla_{xx} L$ evaluated at the local minimizer has to be positive definite on the tangent space defined by the active inequality and equality constraints, that is
\begin{align}\label{eq:sosc}
	\forall \makev{s}\neq 0\in \R^{n_x},	\text{ it holds } \makev{s}^\top \nabla_{xx}L(\xv^\star, \lv^\star, \mv^\star, \pv)\makev{s} > 0, \text{ if }
	\begin{bmatrix}
		J_x \gv^{\mathbb{A}}(\xv^\star, \pv)\\
		J_x \hv(\xv^\star, \pv)
	\end{bmatrix}\makev{s} = 0.
\end{align} 
For further details refer to~\cite{nocedal2006numericaloptimization}.
Decision variables and Lagrange multipliers are often also referred to as primal and dual variables, respectively\MKHcorrection{, and together have dimension $n_y = n_x + n_h + n_g$}.
 
\subsection{Local existence of smooth solution branches}\label{sec:Fiacco}
\MKHcorrection{The Basic Sensitivity Theorem~\cite[Thm.~3.2.2]{fiacco1983sensitivity} is a classical result for parametric optimization problems, and its assumptions are restated in \Cref{assump:fiacco}.
Its first item follows from the twice continuous differentiability of $f$, $\gv$, and $\hv$, see~\ref{eq:pOP}, and is restated to show which regularity \Cref{thm:fiacco} requires.}
Following~\cite{fiacco1983sensitivity}, let $\mathbb{B} \subseteq \R^q$ denote a (possibly bounded) set. For $b\in\mathbb{B}$, a neighborhood of $b$, not necessarily centered at $b$, is denoted $\mathfrak{N}_b\subset \mathbb{B}$. Below, this is applied with $\mathbb{B} = \R^{n_p}$ and $\mathbb{B}=\bar{\mathcal{P}}$.

\begin{assumption}\label{assump:fiacco}
	Without loss of generality, let $\xv^\star$ denote a local minimizer to $\mathrm{OP}(0)$, i.e.\ the optimization problem \ref{eq:pOP} for $\pv=0$. Then, we make the following assumptions.
	\begin{itemize}
		\item \MKHcorrection{On a neighborhood of $(\xv^\star, 0)$ in $\setx\times\setp$, $f$, $\gv$, and $\hv$ are twice continuously differentiable in $\xv$, and $\gv$, $\hv$, $\nabla_{\xv} f$, $\nabla_{\xv} {\gv}^i$, $\nabla_{\xv} {\hv}^j$, $i=1,\dots,n_g,\, j=1,\dots,n_h$, are once continuously differentiable in $\pv$.}
		\item LICQ and SOSC, see \eqref{eq:licq} and \eqref{eq:sosc}, hold at $(\xv^\star,0)$ with associated Lagrange multipliers $\lv^\star$ and $\mv^\star$.
		\item The multipliers $\mv^\star$ satisfy strict complementary slackness, i.e.\ ${\mv^\star}^i>0$ for all $i \in \{i \, | \, {\gv}^i(\xv^\star, 0) = 0\}$.
	\end{itemize}
\end{assumption}\noindent

\begin{theorem}\label{thm:fiacco}
	Under \Cref{assump:fiacco} it holds for $\mathrm{OP}(\pv)$ and its local minimizer $\xv^\star$ that
	\begin{itemize}
		\item $\xv^\star$ is a locally isolated minimizer of $\mathrm{OP}(0)$ and the associated Lagrange multipliers $\lv^\star, \mv^\star$ are unique,
		\item \MKHcorrection{there exist a neighborhood $\mathfrak{N}_{0}\subset \R^{n_p}$ of $0$ and} a unique, once continuously differentiable function
		\[
			{\yv}^\star: \MKHcorrection{\mathfrak{N}_{0} \to \setx\times\seth\times\setg,\ } \pv \mapsto ({\yv}_x^\star(\pv), {\yv}_\lambda^\star(\pv), {\yv}_\mu^\star(\pv)),
		\]
		called the \emph{solution branch}, with $\yv^\star(0) = (\xv^\star, \lv^\star, \mv^\star)$ and such that\MKHcorrection{, for all $\pv \in \mathfrak{N}_{0}$,} the SOSC for a local minimizer of problem $\mathrm{OP}(\pv)$ are satisfied. Hence, $\yv_x^\star(\pv)$ is a locally unique local minimizer of problem $\mathrm{OP}(\pv)$ with associated unique Lagrange multipliers ${\yv}_\lambda^\star(\pv), {\yv}_\mu^\star(\pv)$,
		\item for all $\pv$ in $\mathfrak{N}_{0}$, the set of active inequality constraints remains unchanged,
		strict complementary slackness holds, and the union of gradients of active inequality and equality constraints is linearly independent at $\yv_x^\star(\pv)$.
	\end{itemize}
\end{theorem}\noindent

\subsection{Patched solutions}\label{sec:patched}
Note that Fiacco's setting in \Cref{sec:Fiacco} holds locally on a neighborhood around a nominal solution.
On a global scale, parametric optimization problems are known to exhibit jumps in the decision variables and Lagrange multipliers, for example due to changes in the set of active inequality constraints~\cite{guddat1990jumps}.
For this work, we restrict our application to optimization problems for which the non-smooth points in the patched solution $\pv \mapsto \yv^\star(\pv)$ \MKHcorrection{form a set of measure zero, e.g.\ finitely many points for $n_p = 1$}.
Then, everywhere else,  $\yv^\star(\pv)$ is at least continuously differentiable.

The following conditions on parametric optimization problems are derived in order to guarantee that $\pv \mapsto \yv^\star(\pv)$ is piecewise continuous. 
Then, we extend Fiacco's setting of~\Cref{thm:fiacco} to the case in which multiple nominal solutions are known, so that a patched solution with finitely many jumps can be constructed.

Let $\bar{\mathcal{P}} \subset \R^{n_p}$ denote a closed region of interest for the parameter $\pv$.
Further, let $\pi_{a}(\pv) = a - \pv$ denote a shifting operator.
Note that \Cref{thm:fiacco} considers the nominal solution for parameter value $\pv=0$.
However, this is w.l.o.g., so if
$\xv_i^\star$ is the local minimizer to OP$(\hat{\pv}_i)$, it also is an local minimizer to OP$(\pi_{\hat{\pv}_i} (0))$.
If \Cref{assump:fiacco} holds at $\hat{\pv}_i$, we can apply \Cref{thm:fiacco} to the shifted problem and obtain the existence of a neighborhood $\mathfrak{N}_{\hat{\pv}_i}$.
\begin{assumption}\label{assump:covering}
	Given a parametric optimization problem $\mathrm{OP}(\pv)$ and region of interest $\bar{\mathcal{P}} \subset \R^{n_p}$, we assume that there exist finitely many parameter values $\hat{\pv}_0,\hat{\pv}_1,\dots,\hat{\pv}_M$, with  $M <\infty$, such that \Cref{thm:fiacco} holds for all $\mathrm{OP}(\pi_{\hat{\pv}_i}(0))$, $i=0,1,\dots,M$, and 
	and $\bar{\mathcal{P}} \subset \bigcup_{i=0}^{M} \mathfrak{N}_{\hat{\pv}_i}$, i.e.\ $\bar{\mathcal{P}}$ is covered by finitely many neighborhoods $\mathfrak{N}_{\hat{\pv}_i}$.
\end{assumption}\noindent

The restrictiveness of \Cref{assump:covering} depends on the given problem instance.
As mentioned above, active set changes are one reason that restricts the neighborhoods in which \Cref{thm:fiacco} holds.
Practically, one would choose $\hat{\pv}_i$ in each combination of active constraints and check for local regularity of cost and constraint functions.
The size of the region of interest is an additional design parameter that can be chosen to reduce the minimally required number of parameters.

\begin{definition}\label{def:Patch}
Given a set of solution branches $\yv^\star_i$ on $\mathfrak{N}_{\hat{\pv}_i}$, for $i=0,\dots,M$, satisfying \Cref{assump:covering}, we define a patched solution $\yv^\star:\bar{\mathcal{P}}\to\setx\times\seth\times\setg$ by
	\[ \yv^\star(\pv) = \yv^\star_i(\pv) \text{ with } i= \min \{ j \in \{0,\dots,M\} \, \vert \,  \pv \in \mathfrak{N}_{\hat{p}_j}  \}. \]
\end{definition}
\begin{lemma}
The patched solution is a well-defined function if, for any $\pv\in\bar{\mathcal{P}}$ which lies in the intersection of two neighborhoods, i.e.\ $\pv \in\mathfrak{N}_{\hat{\pv}_i} \cap \mathfrak{N}_{\hat{p}_j}$, $i\neq j, \, 0\leq i,j \leq M$, $\mathrm{OP}(\pv)$ has a unique global minimizer.
\end{lemma}

This follows by construction, since, on the intersection $\mathfrak{N}_{\hat{\pv}_i} \cap \mathfrak{N}_{\hat{p}_j}$, the two solution branches $\yv^\star_i$, $\yv^\star_j$ have to coincide, excluding the possibility of multiple global minimizers.

In case the solution branches on the intersection $\mathfrak{N}_{\hat{\pv}_i} \cap \mathfrak{N}_{\hat{p}_j}$ differ (e.g., because of multiple global or local minimizers), \Cref{def:Patch} guarantees that the patched solution is piecewise continuously differentiable with finitely many jumps.
\KFcorrection{Exploiting the enumeration of the parameter values, exactly one solution branch is selected for each $p\in \bar{\mathcal{P}}$.}
\KFcorrection{For both cases, we can summarize our derivation in the following theorem.}

\begin{theorem}
If a function $\yv^\star$, as defined in \Cref{def:Patch}, exists on $\bar{\mathcal{P}}$, and \Cref{assump:covering} holds, from \Cref{thm:fiacco}, it follows that $\yv^\star$ is piecewise $\mathcal{C}^1$ on $\bar{\mathcal{P}}$ with finitely many jumps in all components, i.e.\ in primal and dual variables.
\end{theorem}

\KFcorrection{Recall that for \Cref{thm:fiacco}, \Cref{assump:fiacco} is crucial.}
\MKHcorrection{For a scalar parameter, the points at which \Cref{assump:fiacco} is violated are called singular points~\cite[Sec.~2.5]{guddat1990jumps}.
	Apart from degenerate cases, these singular points are isolated~\cite[Thm.~2.5.2]{guddat1990jumps}.
	At each of them, the solution branch either continues without being differentiable, e.g.\ at a change of the set of active inequality constraints, or ends, so that the patched solution jumps to another solution branch.
	We denote by $\mathcal{Z}\subset\bar{\mathcal{P}}$ the set of parameter values at which a solution branch reaches a singular point.
	If $\bar{\mathcal{P}}$ and the solution branches are bounded, $\mathcal{Z}$ is finite and thus has Lebesgue measure zero.}

\KFcorrection{This motivates to approximate the patched solution $\yv^\star$  within the set of continuous functions.
	Specifically, we choose neural networks with continuous activation functions, since these networks are continuous by construction.
	Discontinuities of the patched solution will inevitably cause approximation errors.}

\KFcorrection{Turning our perspective from the existence of the patched solution to constructing approaches for approximating it, let us also review \Cref{assump:covering}.
Aiming for data-based learning, the given parameter values $\hat{\pv}_0,\dots,\hat{\pv}_M$ with corresponding nominal solutions shall be chosen such that their neighborhoods in the sense of \Cref{thm:fiacco} cover the region of interest $\bar{\mathcal{P}}$.

Locally about a nominal solution $\yv^\star_i(\hat{\pv}_i)$, the uniqueness of the solution branch $\yv^\star_i$ is guaranteed by \Cref{thm:fiacco}.
As in classical optimization, we interpret the necessary conditions for optimality as a system of constraints and aim for an algorithmic approach for approximating its solution set.

If several local minimizers exist, a chosen nominal solution $\yv^\star_i(\hat{\pv}_i)$ selects a unique solution branch to be approximated.
}


%% file: inc/optinns.tex
\section{Optimality-Informed Neural Networks}\label{sec:optinn}
We introduce Optimality-Informed Neural Networks (OptINNs) as a class of theory-informed neural networks that incorporate optimality conditions into the training process by designing a loss function based on the KKT conditions.
Moreover, some problem constraints are directly transferred into a suitable architectural design of the network.
\Cref{fig:optinn_sketch} provides a schematic overview of the full pipeline.

Given a parametric optimization problem of the form~\ref{eq:pOP} and a finite set $\{\yv^\star_i(\hat{\pv}_i)\}_{i=0}^M$ of nominal solutions for the given parameters $\hat{\pv}_0,\dots,\hat{\pv}_M$ with $0\leq M<\infty$, we aim for approximating the patched solution $\yv^\star$ of the parameterized problem for all $\pv \in \bar{\mathcal{P}}$ as defined in \Cref{sec:patched}.
Our approach is to combine classical nonlinear optimization with \MKHcorrection{neural networks}.
More concretely, we aim to approximate the mapping by a function parametrized in $\tv \in \R^{n_\theta}$ 
\begin{equation}
\begin{aligned}
	\mathcal{N}:\ &\mathbb{R}^{n_p} \times \mathbb{R}^{n_\theta} \to \R^{n_y} = \R^{n_x+n_h+n_g} \\
	&(\pv,\tv)\mapsto \nnt{\pv}={(\N_x^\top(\pv,\tv),\N^\top_\lambda(\pv,\tv),\N_{\mu}^\top(\pv,\tv))}^\top,
\end{aligned}\label{eq:def_NN}
\end{equation}
where $\N(\pv,\tv)$ is decomposed into $\N_x(\pv, \tv)$ for approximating the primal component $\yv_x^\star(\pv)$ of the patched solution, $\N_\lambda(\pv, \tv)$ and $\N_\mu(\pv, \tv)$ for the dual components $\yv^\star_\lambda(\pv)$ and $\yv^\star_\mu(\pv)$, respectively.
We design the loss function as a convex combination of two terms,
\begin{align}\label{eq:combined_loss}
\namedL{}(\tv; \alpha)=\alpha\Lkkt(\tv) + (1-\alpha)\Lmse(\tv)
\end{align}
 with $\alpha \in [0,1]$.
 The mean squared error (MSE) loss, \Lmse, is the standard data-based loss for regression problems.
 It aggregates the pointwise error metric \Dmse{} across all samples, where a single sample $\yv$ is compared with an estimate $\hat{\yv}$:
\begin{align}\label{eq:mse}
    \Lmse(\tv) = \sum_{i=0}^M \Dmse(\N(\hat{\pv}_i, \tv), \MKHcorrection{\yv^\star_i(\hat{\pv}_i)}), \quad \text{with} \quad \Dmse(\hat{\yv}, \yv) = \left\Vert \hat{\yv} - \yv\right\Vert_2^2.
\end{align}
The term $\Lkkt$ encodes optimality information, as described below.
For $\alpha = 0$, we recover the pure regression problem, whereas for $\alpha = 1$, the parameter identification replaces the classical nonlinear optimization, i.e.\ zero-finding of KKT conditions.
\MKHcorrection{If no nominal solutions are available, $\Lmse$ is omitted, and the approach can still be applied as for $\alpha = 1$.}

\subsection{Optimality-informed design of loss functions}\label{sec:loss}
To embed conditions of optimality into the training process, we construct a loss function based on the Karush-Kuhn-Tucker conditions, see~\ref{eq:kkt}.
Recall that the KKT conditions are either inequality or equality constraints, so we can interpret a deviation from equality constraints as a measure for optimality of the candidate solution.
For the inequality constraints $\gv$, we make use of the rectified linear unit \(\relu: \R \to \R, \xv\mapsto\max(0, \xv)\) to only penalize violations of the inequality, i.e.\ penalize $\gv(\xv, \pv) > 0$\MKHcorrection{, so that satisfied inequality constraints yield a penalty of exactly zero}.
Non-negativity of the Lagrange multiplier $\mv$ is guaranteed by architectural design, as discussed in \Cref{sec:architecture}.

\begin{definition}[Unimodal Penalty Function]\label{def:unimodal_penalty}
    A continuous function $\phi: \mathbb{R} \to [0, \infty)$ is called a unimodal penalty function if $\phi(0) = 0$, $\phi$ is differentiable on $\mathbb{R} \setminus \{0\}$, and
    \(
        x \phi'(x) > 0 \quad \text{for all } x \neq 0.
    \)
\end{definition}\noindent
Examples of unimodal penalty functions are the absolute value function or the square function.
\MKHcorrection{Like the ReLU function, the absolute value function is not differentiable at zero, where PyTorch's automatic differentiation uses the subgradient zero.}

With unimodal penalty functions \(P^i\), $i \in \mathcal{T} = \{\mathrm{Stat, FeasG, FeasH, CSl}\}$, for a given point $(\xv, \lv, \mv)$ and associated parameter $\pv$, we can formulate deviations from the KKT conditions~\ref{eq:kkt} as
\begin{equation}\label{eq:kkt_deviations}
	\begin{aligned}
	\Ds(\xv, \lv, \mv, \pv) &= \frac{1}{n_x}\sum_{i=1}^{n_x} \Ps\left(\frac{\partial L}{\partial \xv^i} (\xv, \lv, \mv, \pv)\right),&
	\Dpg(\xv, \pv) &= \frac{1}{n_g} \sum_{i=1}^{n_g} \Ppg\left(\relu(\gv^i(\xv, \pv))\right),\\
	\Dcs(\xv, \mv, \pv) &= \frac{1}{n_g} \sum_{i=1}^{n_g} \Pcs\left(\mv^i \gv^i(\xv, \pv)\right),\label{eq:Lcs}&
	\Dph(\xv, \pv) &= \frac{1}{n_h} \sum_{i=1}^{n_h} \Pph\left(\hv^i(\xv, \pv)\right).
	\end{aligned}
\end{equation}

We sample a set of parameters $\{\pv_i\}_{i=1}^N$ from a uniform distribution $\mathcal{U}$ with lower and upper bounds $\underline{\pv}$ and $\overline{\pv}$.
\MKHcorrection{In each training step, a new set of parameters is uniformly sampled from $\bar{\mathcal{P}}$.}
Using the neural network $\mathcal{N}$ parametrized by $\tv$ (cf.~\eqref{eq:def_NN}), we predict the corresponding primal and dual variables for these samples.
Averaging the deviations from the KKT conditions over this set yields the functions
\begin{equation}\label{eq:kkt_residuals}
    \begin{aligned}
    \Ls(\tv)  &= \frac{1}{N} \sum_{i=1}^{N} \Ds(\nnt[x]{\pv_i}, \nnt[\lambda]{\pv_i}, \nnt[\mu]{\pv_i}, \pv_i),&
    \Lpg(\tv) &= \frac{1}{N} \sum_{i=1}^{N} \Dpg(\nnt[x]{\pv_i}, \pv_i),\\
    \Lcs(\tv) &= \frac{1}{N} \sum_{i=1}^{N} \Dcs(\nnt[x]{\pv_i}, \nnt[\mu]{\pv_i}, \pv_i),&
    \Lph(\tv) &= \frac{1}{N} \sum_{i=1}^{N} \Dph(\nnt[x]{\pv_i}, \pv_i).
\end{aligned} 
\end{equation}
In the sequel, we refer to these functions as KKT-residuals.
They serve as the basis for identifying the parameters $\tv$ of the associated neural network $\mathcal{N}$.
\MKHcorrection{Since the parameters are sampled uniformly, each KKT-residual is a Monte Carlo estimate of the corresponding mean deviation from the KKT conditions over the region of interest, e.g.\
\[
	\Ls(\tv) \approx \frac{1}{\mathrm{vol}(\bar{\mathcal{P}})}\int_{\bar{\mathcal{P}}} \Ds(\nnt[x]{\pv}, \nnt[\lambda]{\pv}, \nnt[\mu]{\pv}, \pv)\,\mathrm{d}\pv.
\]
Hence, the set $\mathcal{Z}$ of measure zero, see \Cref{sec:patched}, does not contribute to these integrals, and the sampled parameters lie in $\mathcal{Z}$ with probability zero.}
The weighted and summed KKT-residuals give the KKT-loss function
\begin{align}\label{eq:lkkt}
	\Lkkt(\tv) = \os\Ls(\tv) + \opg\Lpg(\tv) + \oph\Lph(\tv) + \ocs\Lcs(\tv)
\end{align}
that combines the weighted deviations of the model predictions for all KKT-residuals.
The weights $\omega^i, \, i\in \mathcal{T}$, are used for balancing the different loss terms, and their tuning will be explained in \Cref{sec:weightbalancing}.
Next, we show that only KKT-points, combinations of primal and dual variables satisfying the KKT conditions, produce \MKHcorrection{zero deviations from the KKT conditions}.
\begin{lemma}[\MKHcorrection{Zero deviations for KKT-points}] \label{lem:Dloss_KKT}
	Given an optimization problem of the form \emph{\ref{eq:pOP}} with parameter $\pv \in \R^{n_p}$, primal and dual variables $(\xv^\star, \lv^\star, \mv^\star) \in \mathbb{R}^{n_x} \times  \mathbb{R}^{n_h} \times \mathbb{R}^{n_g}_{\geq 0}$, the \MKHcorrection{deviations \eqref{eq:kkt_deviations} from the KKT conditions} with corresponding unimodal penalty functions $P^i$, $i\in\mathcal{T}$, are zero, i.e.\ \[\Ds(\xv^\star, \lv^\star, \mv^\star, \pv) = \Dpg(\xv^\star, \pv) = \Dph(\xv^\star, \pv) = \Dcs(\xv^\star,\mv^\star, \pv) = 0,\]
	if and only if $(\xv^\star, \lv^\star, \mv^\star)$ and $\pv$ satisfy~\emph{\ref{eq:kkt}}.
\end{lemma}
\begin{proof}
	See Appendix~\ref{app:ProofLemma}.
\end{proof}\noindent
\MKHcorrection{For positive weights, this carries over to the KKT-loss.
\begin{corollary}[Zero KKT-loss for KKT-points]\label{cor:Lkkt_zero}
	Let $\omega^i > 0$, $i\in\mathcal{T}$, and let $\N_\mu(\pv_i, \tv) \geq 0$ for $i = 1,\dots,N$, see \Cref{sec:architecture}.
	Then, $\Lkkt(\tv) = 0$ if and only if $\N(\pv_i, \tv)$ is a KKT-point of $\mathrm{OP}(\pv_i)$ for all sampled parameters $\pv_i$, $i = 1,\dots,N$.
\end{corollary}
\begin{proof}
	Since the unimodal penalty functions are non-negative, all deviations from the KKT conditions and thus all KKT-residuals in~\eqref{eq:kkt_residuals} are non-negative.
	For $\omega^i > 0$, $i\in\mathcal{T}$, it follows that $\Lkkt(\tv) = 0$ if and only if all deviations are zero at all sampled parameters, and the claim follows from \Cref{lem:Dloss_KKT}.
\end{proof}}\noindent

\MKHcorrection{If $\mathrm{OP}(\pv)$ has an empty feasible set, every prediction violates at least one inequality or equality constraint, so that \Dpg or \Dph, and thus the KKT-loss, is positive at $\pv$.
Hence, \Dpg and \Dph indicate such parameter values, similar to the constraint violation of the penalty-based neural networks in \Cref{sec:pmnn}.}

\MKHcorrection{\subsection{Box-and-bound-constraint-satisfying architecture}\label{sec:architecture}}
A comparison of the KKT conditions with the formulated penalties shows that we have not included a penalty for the dual feasibility, cf.\ \ref{eq:kkt} and~\eqref{eq:kkt_residuals}.
We specifically did not include this as a penalty, as it can be fulfilled by architectural design. 
The architecture used in our experiments is based on a standard MLP.
An MLP is the repeated composition of linear-affine transformations, defined by weight matrices \(\makev{W}^{(l)}\in\R^{u_l\times u_{l-1}}\) and bias vectors \(\makev{b}^{(l)}\in\R^{u_l}\), and nonlinear activation functions $\makev{a}^{(l)}: \R^{u_l}\to\R^{u_l}$.
The entries of all $\makev{W}^{(l)}$ and $\makev{b}^{(l)}$ form the parameter vector $\tv$.
 In our setting, the resulting MLP $\mathcal{N}:\ \mathbb{R}^{n_p} \times \mathbb{R}^{n_\theta} \to \R^{n_y}$ with $L$ layers is given as
 $
 \mathcal{N}(\pv,\tv) \coloneq \makev{\nu}^{(L)}(\pv)
 $
 with
\begin{align}\label{eq:mlp}
	\makev{\nu}^{(0)} &= \pv,\\
	\makev{z}^{(l)} &= \makev{W}^{(l)}\,\makev{\nu}^{(l-1)} + \makev{b}^{(l)} \ \forall \ l \in \{1, \dots, L\},\\
	\makev{\nu}^{(l)} &= \makev{a}^{(l)} \left(\makev{z}^{(l)}\right) \ \forall \ l \in \{1, \dots, L\}.
\end{align}

To enforce dual feasibility, the components of \( \makev{z}^{(L)} \) associated with the dual variables of inequality constraints are transformed in the final activation layer \( \makev{a}^{(L)} \) via the Softplus function
\[
\text{Softplus}: \mathbb{R} \to (0,\infty), \quad x \mapsto \ln(1 + e^x),
\]
so that the resulting output satisfies this condition (see~\ref{eq:kkt}).
\MKHcorrection{Unlike the ReLU function, whose derivative is zero for negative arguments, the Softplus function has a positive derivative everywhere, so that predicted Lagrange multipliers close to zero still receive gradients and can increase when the corresponding inequality constraint becomes active.
However, this derivative tends to zero for negative arguments of large magnitude, which may slow down this increase at changes of the set of active inequality constraints.
The same holds for the transformations in the following remark at active bounds.}
\begin{remark}
If some of the inequality constraints in~\ref{eq:pOP} are of the form $x\le \overline{x}, \underline{x}\le x$, or $\underline{x}\le x \le \overline{x}$, i.e.\ lower bounds, upper bounds, or box constraints, these can be fulfilled by design as well.
To enforce feasibility of the predicted decision variables directly through the network architecture, we apply suitable activation functions tailored to the type of constraint.
\MKHcorrection{Let $z$ denote the component of the pre-activation $\makev{z}^{(L)}$ in~\eqref{eq:mlp} associated with such a decision variable $x$, and $\bound{x}$ the corresponding output of the network.}
Lower bounds and upper bounds are fulfilled by $\bound{x} = \underline{x} + \text{{Softplus}}(\MKHcorrection{z})$ and $\bound{x} = \overline{x} - \text{{Softplus}}(\MKHcorrection{z})$, respectively.
For box constraints, the sigmoid function $\sigma: \mathbb{R} \to (0, 1), \ x \mapsto \frac{1}{1 + e^{-x}}$ is shifted and scaled such that \(\bound{x} = \underline{x} + (\overline{x} - \underline{x}) \sigma(\MKHcorrection{z}) \in (\underline{x}, \overline{x})\).
 
Concatenating these transformations for bound- and box-constrained decision variables with identity mappings for unconstrained variables and the Softplus function for dual variables of inequality constraints gives a vector-valued function forming the final nonlinear activation layer $\makev{a}^{(L)}$.
Since bound and box constraints are now always satisfied, their contributions to $\Lpg(\tv)$ are always 0.
\end{remark}
\MKHcorrection{For general inequality and equality constraints, no such closed-form transformation exists.
Projection-based methods such as DC3~\cite{donti2021dc3} and NLPOpt-Net~\cite{roy2026nlpoptnet}, see \Cref{sec:rel_works_nn}, therefore restore feasibility by iterative steps in each evaluation of the network, also at inference.
Since we aim to avoid such iterations, OptINNs enforce constraints by architectural design only where a closed-form transformation exists, and include all remaining constraints in the KKT-loss.}

\input{figs/optinn_sketch}

\subsection{OptINN training}\label{sec:optinn_training}
From the PINN literature, mainly from \cite{wang_experts_2023} for this work, some key insights for training theory-informed neural networks can be drawn.
Here, we summarize the methods used for this work's numerical results.

\subsubsection{Tuning the data loss and optimality trade-off}\label{sec:alpha_scheduling}
As introduced previously, the parameter $\alpha$ tunes the trade-off between $\Lmse$ and $\Lkkt$.
While our main focus is the minimization of $\Lkkt$, the topologies of the $\mathcal{D}^i$, $i\in\mathcal{T}$, can, depending on the optimization problem, be multimodal.
This makes them more involved \MKHcorrection{than \Dmse, which is convex in its arguments, whereas the $\mathcal{D}^i$ are in general not}, see for example~\cite{krishnapriyan2021pinnfailure}.
Thus, we divide training into three phases:
\begin{enumerate}
	\item initialization phase of length $d^\mathrm{init}$: $\alpha$ is kept constant at its lowest value $\underline{\alpha}$
	\item annealing phase of length $d^\mathrm{anneal}$: $\alpha$ is increased following a cosine function
	\item finalization phase of length $d^\mathrm{final}$: $\alpha$ is kept at its highest value $\overline{\alpha}$.
\end{enumerate}
For the epoch $d$, where an epoch refers to one complete pass through the entire training dataset by the learning algorithm, 
\begin{align}
	\alpha(d) = \begin{cases*}
		\underline{\alpha} & if $d \le d^\mathrm{init}$ \\
		\underline{\alpha}+\frac{1}{2}(\overline{\alpha} - \underline{\alpha})(1 - \cos(\pi \frac{d-d^\mathrm{init}}{d^\mathrm{anneal} })) & if $d^\mathrm{init} < d \le d^\mathrm{init} + d^\mathrm{anneal}$ \\
		\overline{\alpha} & if $d^\mathrm{init} + d^\mathrm{anneal} \le d \le d^\mathrm{init} + d^\mathrm{anneal} + d^\mathrm{final}$
	\end{cases*}.
\end{align}
\Cref{fig:cosine_scheduling} shows one cycle for initialization and finalization phases of length 50, $\underline{\alpha} = 0.1$, and $\overline{\alpha} = 0.9$. 
\begin{figure}[htb]
	\centering
	\begin{tikzpicture}
		\begin{axis}[
			xlabel={Epoch $d$},
			ylabel={$\alpha$},
			ytick={0.1, 0.9},
			yticklabels={0.1, 0.9},
			width=0.5\linewidth,
			height=0.3\linewidth,
			]
			\addplot[mark=none, color=udsblue, line width=1pt, domain=0:250, samples=200] 
				{0.1 + 0.4 * (1 - cos(deg(pi * min(max((x - 50) / 150, 0), 1))))};
		\end{axis}
	\end{tikzpicture}
	\caption{Cosine scheduling for $\alpha$ with $\underline{\alpha}=0.1$, $\overline{\alpha}=0.9$, and initialization and finalization phases of lengths $d^\mathrm{init} = d^\mathrm{final}=50$ over 250 epochs}
	\label{fig:cosine_scheduling}
\end{figure}
This progression allows the MSE-loss to dominate early training, while gradually shifting focus to the KKT-loss function.
\MKHcorrection{While the MSE-loss is evaluated on few training data points and may be overfit, the KKT-loss is evaluated on a new set of sampled parameters in each training step. 
Thus, training minimizes the integrals over the region of interest, see \Cref{sec:loss}, rather than the KKT-residuals at a fixed set of parameters.
On a neighborhood on which the deviations from the KKT conditions are zero, fitting a training data point yields the solution branch, see \Cref{sec:patched}.}

\subsubsection{\MKHcorrection{Gradient} balancing}\label{sec:weightbalancing}
Similar to PINNs, a key point in training OptINNs is the imbalance between the gradients of the different losses~\cite{wang_experts_2023}. 
During training, losses with larger gradients are prioritized in minimization, resulting in a potential deterioration of other losses.
To handle this, automated loss balancing has been introduced for scaling the different costs in such a way that the resulting descent direction lies in the cone of common descent for all losses.
In~\cite{wang_experts_2023}, two methods, one based on the so-called neural tangent kernel matrix and one based on the gradients of loss functions, are summarized.
Employing the second approach, the weights $\os, \opg, \oph$, and $\ocs$ \eqref{eq:lkkt} are set according to the magnitude of the gradients with respect to the MLP's parameters $\tv$  \MKHcorrection{as the exponential moving average over}
\begin{align}
	\bar{\omega}^i = \begin{cases*}
		\frac{\sum_{j\in\MKHcorrection{\mathcal{J}}} \Vert \nabla_\theta \mathcal{L}^j(\tv)\Vert}{\Vert\nabla_\theta\mathcal{L}^i(\tv)\Vert}, & $\MKHcorrection{i \in \mathcal{I} \ \text{and} \ }\Vert\nabla_\theta{\namedL{}}^i(\tv)\Vert > \beta$\\
		1, & \text{else}
	\end{cases*}
\end{align}
for some small $\beta$\MKHcorrection{, where $\mathcal{I} \subseteq \mathcal{J}$ are the indices of the weighted losses and of all losses, respectively}.
\MKHcorrection{In the following, $\mathcal{I} = \mathcal{T}$ and $\mathcal{J} = \mathcal{T} \cup \{\mathrm{MSE}\}$.
Since the numerator includes $\Vert\nabla_\theta\Lmse(\tv)\Vert$, each weighted KKT term attains at least the gradient magnitude of the MSE-loss, so that the trade-off between $\Lkkt$ and $\Lmse$ is set by $\alpha$ rather than by their gradient magnitudes.
As a design choice, the MSE-loss is not weighted, so that its influence on the training is controlled by $\alpha$ only.}
\MKHcorrection{In both cases, $\omega^i \geq 1$, since the numerator in the first case contains $\Vert\nabla_\theta\mathcal{L}^i(\tv)\Vert$, so that \Cref{cor:Lkkt_zero} applies.}

Before we numerically validate the OptINN approach in detail, we present the variant based on the quadratic-penalty method, which is known from the literature and serves as a baseline for the OptINN validation.




%% file: figs/optinn_sketch.tex
\begin{figure}[htb]
	\centering
	\ifdefined\tikzexternalize 
		\tikzset{
			external/export next=false,
			>=latex, 
			node/.style={thick,circle,draw=udsblue,minimum size=22,inner sep=0.5,outer sep=0.6},
			node in/.style={node,green!20!black,draw=green!30!black,fill=udsgreen!25},
			node hidden/.style={node,blue!20!black,draw=udsblue!30!black,fill=udsblue!20},
			node out/.style={node,red!20!black,draw=udsred!30!black,fill=udsred!20},
			connect/.style={thick,black}, 
			node 1/.style={node in}, 
			node 2/.style={node hidden},
			node 3/.style={node out},
			box/.style={draw, thick, rounded corners, inner sep=1em, fill=gray!10, align=center},
			box2/.style={draw, thick, rounded corners, inner sep=1em, fill=white, align=center, minimum height=12mm},
			box3/.style={thick, rounded corners, inner sep=1em, fill=none, align=center}
		}
	\else
		\tikzset{
			>=latex, 
			node/.style={thick,circle,draw=udsblue,minimum size=22,inner sep=0.5,outer sep=0.6},
			node in/.style={node,green!20!black,draw=green!30!black,fill=udsgreen!25},
			node hidden/.style={node,blue!20!black,draw=udsblue!30!black,fill=udsblue!20},
			node out/.style={node,red!20!black,draw=udsred!30!black,fill=udsred!20},
			connect/.style={thick,black}, 
			node 1/.style={node in}, 
			node 2/.style={node hidden},
			node 3/.style={node out},
			box/.style={draw, thick, rounded corners, inner sep=1em, fill=gray!10, align=center},
			box2/.style={draw, thick, rounded corners, inner sep=1em, fill=white, align=center, minimum height=12mm},
			box3/.style={thick, rounded corners, inner sep=1em, fill=none, align=center}
		}
	\fi
	
	\def\nstyle{int(\lay<\Nnodlen?min(2,\lay):3)} 
	\newcommand\setAngles[3]{
		\pgfmathanglebetweenpoints{\pgfpointanchor{#2}{center}}{\pgfpointanchor{#1}{center}}
		\pgfmathsetmacro\angmin{\pgfmathresult}
		\pgfmathanglebetweenpoints{\pgfpointanchor{#2}{center}}{\pgfpointanchor{#3}{center}}
		\pgfmathsetmacro\angmax{\pgfmathresult}
		\pgfmathsetmacro\dang{\angmax-\angmin}
		\pgfmathsetmacro\dang{\dang<0?\dang+360:\dang}
	}
	\begin{tikzpicture}[x=2.2cm,y=1.4cm]
		\message{^^JNeural network, shifted}
		\readlist\Nnod{3,5,5,3,3} 
		\readlist\Nstr{n_p,u_{\prev},u_{\prev},n_y,n_y} 
		\readlist\Cstr{\strut p,\nu^{(\prev)},\nu^{(\prev)},\strut y,\strut \bound{y}} 
		\def\yshift{0.5} 

		\message{^^J  Layer}
		\foreachitem \N \in \Nnod{ 
			\def\lay{\Ncnt} 
			\pgfmathsetmacro\prev{int(\Ncnt-1)} 
			\message{\lay,}
			\foreach \i [evaluate={\c=int(\i==\N); \y=\N/2-\i-\c*\yshift;
				\index=(\i<\N?int(\i):"\Nstr[\lay]");
				\x=\lay; \n=(\lay>=\Nstrlen-1?3:\nstyle);}] in {1,...,\N}{ 
				\node[node \n] (N\lay-\i) at (\x,\y) {$\Cstr[\lay]_{\index}$};
				
				\ifnum\lay>1 
				\ifnum\lay=\Nstrlen 
				\draw[connect] (N\prev-\i) -- (N\lay-\i);
				\else 
				\foreach \j in {1,...,\Nnod[\prev]}{ 
					\draw[connect,white,line width=1.2] (N\prev-\j) -- (N\lay-\i);
					\draw[connect] (N\prev-\j) -- (N\lay-\i);
				}
				\fi
				\fi 
			}
			\path (N\lay-\N) --++ (0,1+\yshift) node[midway,scale=1.5] {$\vdots$};
		}
		
		\foreach \i in {1,2,3} {
			\draw[->, thick] (N\Nstrlen-\i) --++ (0.85,0) coordinate (out-\i); 
		}
		
		\draw[->, thick] ($(N\Nstrlen-1)!0.48!(out-1)$) -- ++(0,-1.20cm) 
			arc (270:90:-0.2cm) -- ++(0,-1.7cm)
			arc (270:90:-0.2cm) -- ++(0,-1.9cm);
			\draw[->, thick] ($(N\Nstrlen-2)!0.65!(out-2)$) -- ++(0,-1.9cm) 
				arc (270:90:-0.2cm) -- ++(0,-1.9cm) node[] (MSEnode) {};
			
		\draw[->, thick] ($(N\Nstrlen-3)!0.82!(out-3)$) -- ++(0,-2.1cm);
		
		\begin{pgfonlayer}{background}
			\node[box, fit={(N1-1) (N1-3)}, label={[align=center]above:Parameters}] {};
			\node[box, fit={(N2-1) (N3-5)}, label={[align=center]above:MLP}] (BoxOptINN) {};
			\node[box, fit={(N4-1) (N4-3)}, label={[align=center]above:Trivialization\\ Layer}] (TrivLayer) {};
			\node[box, fit={(N\Nstrlen-1) (N\Nstrlen-3)}, label={[align=center]above: Primal and \\ dual variables}] {}; 
		\end{pgfonlayer}
		
		\node[box, fit={(N4-1) (N4-1)}, minimum height=15mm, below=1.3cm of TrivLayer.south, text depth=2ex] (data) {Data};
		\node[box, anchor=north, minimum height=15mm] at ($(MSEnode |- TrivLayer.south) + (0,-1.3cm)$) (mse) {MSE};
		\draw[->, thick] (data) -- (mse.west);
		\draw[->, thick] (mse) -- ++(0, -1.2cm) node[circle, draw, fill=white, inner sep=0.5mm, anchor=north] (sum) {$\sum$};
		\draw[->, thick] (sum) -| node[pos=0.23,below] {Backpropagation} (BoxOptINN.south);
		
		\def\xmin{6.2}
		\def\xmax{7.2}
		\def\ymin{-2.5}
		\def\ymax{1}
		
		\def\Nx{2} 
		\def\Ny{4} 
		
		\foreach \i in {1,...,\Nx} {
			\foreach \j in {1,...,\Ny} {
				\pgfmathsetmacro\xpos{\xmin + (\xmax - \xmin) * (\i - 1)/(\Nx-1)}
				\pgfmathsetmacro\ypos{\ymin + (\ymax - \ymin) * (\Ny - \j)/(\Ny-1)}
				\node (grid-\j-\i) at (\xpos,\ypos) {};
			}
		}
		\node[box, fit={(grid-1-1) (grid-\Ny-\Nx)}, label={[align=center]above:KKT-loss \eqref{eq:lkkt}}, xscale=1.2, yscale=1.2] (BoxKKT) {};
		\node[box3, fit={(grid-1-1) (grid-1-\Nx)}, label={[align=center]}] {Weighted sum using $\omega^i$};
		\node[box2, fit={(grid-2-1) (grid-2-\Nx)}, label={[align=center]}, text depth=-1ex, text height=1.2ex] {Stationarity};
		\node[box2, fit={(grid-3-1) (grid-3-\Nx)}, label={[align=center]}, text depth=-1ex, text height=1.2ex] {Primal Feasibility};
		\node[box2, fit={(grid-4-1) (grid-4-\Nx)}, label={[align=center]}, text depth=1.5ex] {Complementary Slackness};
		
		\draw[->, thick] (BoxKKT.south) |- (sum.east);
		
	\end{tikzpicture}
	\caption{The architecture and training of a 4-layer OptINN. The problem parameters are passed through the MLP for the given parameter vector $\pv$. The trivialization layer includes all bound and box constraints on the output variables, giving the variables marked with \bound{(\cdot)}. The output $\bound{\yv}\in\R^{n_y}$ corresponds to the concatenation of the primal variables $\xv$ and dual variables $\mv$ and $\lv$}
	\label{fig:optinn_sketch}
\end{figure}

%% file: inc/comparison_pmnn.tex
\section{The KKT-loss function in comparison to constraint-penalization} \label{sec:constraint-penalization}

Although the KKT-loss~\eqref{eq:lkkt} introduces additional complexity, it offers advantages over \MKHcorrection{penalty methods}, as discussed in this section.
We begin by summarizing \MKHcorrection{penalty methods and the derived penalty-based neural networks, including quadratic-penalty-method-based neural networks (QPNNs)}, highlighting the specific weaknesses the latter \MKHcorrection{inherit}.
Subsequently, we compare \MKHcorrection{penalty-based neural networks} (using constraint-penalization) against OptINNs (using the KKT-loss).
Our analysis \MKHcorrection{shows advantages of the OptINN framework in four aspects: (i) the shift of the minimizer of the loss, (ii) the} integration of data-driven terms with the respective optimization losses, \MKHcorrection{(iii)} the suitability of each loss formulation for hyperparameter tuning, and \MKHcorrection{(iv)} the practical restrictions regarding validation and scheduling during training.
\MKHcorrection{Finally, we discuss the limitations of the KKT-loss.}

\MKHcorrection{\subsection{Penalty-based neural networks}\label{sec:pmnn}}
\MKHcorrection{As discussed in \Cref{sec:rel_works_nn}, penalty terms on the constraints are used in ODE-based solver networks, which solve one instance of an optimization problem per run~\cite{rao1992nnmechanical,Lillo1993penaltynn}, and in neural networks that approximate the mapping $\pv \mapsto \xv^\star(\pv)$~\cite{liu2024teaching,donti2021dc3}.
Following the latter, we approximate this mapping with MLPs $\mathcal{N}^\mathrm{P}: \R^{n_p}\times\R^{n_\theta}\to\R^{n_x}$.
Given an optimization problem of the form~\ref{eq:pOP}, penalty methods combine the cost function and constraints into the penalty loss}
\MKHcorrection{%
\begin{align}\label{eq:penalty_loss}
    \namedL{P}_q(\tv) &= \frac{1}{N}\sum_{i=1}^N \namedD{P}_q(\mathcal{N}^\mathrm{P}(\pv_i, \tv), \gamma_g, \gamma_h, \pv_i) \\
	\namedD{P}_q(\xv, \gamma_g, \gamma_h, \pv) &= f(\xv, \pv) + \gamma_g \sum_{i=1}^{n_g} \relu(\gv^i(\xv,\pv))^q + \gamma_h \sum_{j=1}^{n_h} \lvert\hv^j(\xv,\pv)\rvert^q,
\end{align}}
where $\gamma_g$ and $\gamma_h$ are tunable parameters, with higher values resulting in more focus on satisfying the constraints\MKHcorrection{, and $q \in \{1, 2\}$}.
\MKHcorrection{For $q=2$, \eqref{eq:penalty_loss} is the quadratic-penalty loss (QP-loss) $\namedL{QP} = \namedL{P}_2$, $\namedD{QP} = \namedD{P}_2$, and we refer to the MLPs trained on it as QPNNs.
For $q=1$, \eqref{eq:penalty_loss} is an exact penalty function~\cite[Sec.~17.2]{nocedal2006numericaloptimization}.}
Similar to OptINNs, this can be combined with an MSE-loss to incorporate training data by replacing $\Lkkt$ in~\eqref{eq:combined_loss} with \MKHcorrection{$\namedL{P}_q$}.

\MKHcorrection{The penalty-based neural networks in~\cite{rao1992nnmechanical,Lillo1993penaltynn,liu2024teaching,donti2021dc3} use the quadratic penalty, whereas exact penalty functions are a classical alternative in nonlinear programming~\cite[Sec.~17.2]{nocedal2006numericaloptimization}.
We discuss both in this section.
Exact penalty functions share the weaknesses discussed in \Cref{sec:hyperparameter_tuning_metrics,sec:validation}, provide no estimates of the Lagrange multipliers, and avoid the weaknesses discussed in \Cref{sec:penalty_shift,sec:alignment} only if $\gamma_g$ and $\gamma_h$ exceed a bound on the Lagrange multipliers, which is generally unknown for parametric optimization problems.
Therefore, we compare OptINNs against QPNNs in \Cref{sec:numerical_results}.}

\MKHcorrection{\subsection{Shift of the minimizer}\label{sec:penalty_shift}}
One major problem of the QP-loss is that it can only be used to compute the constrained optimum of~\ref{eq:pOP} if the unconstrained minimizer of $f(\xv,\pv)$ is feasible~\cite{nocedal2006numericaloptimization}.
The following example illustrates this limitation.
\begin{example}
	For the equality-constrained optimization problem
	\begin{align}
		\min_{\xv\in\R^2} f(\xv),\ \text{s.t.}\ h(\xv) = 0
	\end{align}
	the associated QP-loss is $f(\xv) + \gamma_h h(\xv)^2$.
	The condition of optimality is $\nabla_\xv f(\xv) + 2\gamma_h h(\xv) \nabla_\xv h(\xv)=0$, where $2\gamma_h h(\xv)$ can be interpreted as an estimate for the Lagrange multiplier $\lambda$.
	Since $h(\xv) = 0$ for any feasible solution, this condition can only be satisfied if $\nabla_\xv f(\xv) = 0$, i.e.\ it can only find the constrained minimizer when it coincides with the unconstrained minimizer.
	Otherwise, the computed solution of the constraint-penalized problem will be infeasible for the original problem.
\end{example}\noindent
Also, if inequality constraints are present in addition, the quadratic-penalty method might yield infeasible solutions.
As a consequence, the minimizer and minimum value of $\mathcal{D}^\mathrm{QP}(\xv, \gamma_g, \gamma_h, \pv)$ vary with the choice of $\gamma_g$ and $\gamma_h$.
For a concrete example, see \Cref{fig:penalty_shift_minima}, where, with increasing $\gamma_g$, the minimizer and minimum of the QP-loss shift towards the constrained optimum ($x=1, f(x)=-1$) without reaching it.
\MKHcorrection{For exact penalty functions, $q=1$, strict local minimizers of the optimization problem are also local minimizers of $\namedD{P}_1$ if $\gamma_g$ and $\gamma_h$ are larger than the absolute values of all Lagrange multipliers~\cite[Thm.~17.3]{nocedal2006numericaloptimization}.
For parametric optimization problems, however, this bound has to hold on the entire region of interest and thus requires an upper bound on the Lagrange multipliers, which is generally unknown.}

In contrast, from the fact that all $P^i, i\in\mathcal{T}$, are unimodal penalty functions, we can deduce that $\Lkkt$, see~\eqref{eq:lkkt}, is lower bounded by zero independent of the values~$\omega^i, i \in \mathcal{T}$.

\input{figs/penalty_method_optimum_shift.tex}

\subsection{Alignment of Loss Minima with Training Data}\label{sec:alignment}
Next, we discuss how  well the loss minima align with training data, i.e.\ how well numerical optimizers' solutions minimize the losses.
\MKHcorrection{Numerical optimizers like IPOPT~\cite{wachter2006ipopt} or WORHP~\cite{Bueskens2013worhp} terminate when the deviations from the stationarity condition, the constraint violation, and the complementary slackness are below prescribed tolerances.}
As a consequence, solutions computed with these solvers \MKHcorrection{exhibit} a KKT-loss \eqref{eq:lkkt} close to 0.
This implies that, for any given data point, the MSE-loss (based on numerical solver data) and the KKT-loss share the same optimum, within the limits of the numerical solver's tolerances.
From \MKHcorrection{\Cref{sec:penalty_shift}}, it is evident that, since the QP-loss predicts infeasible solutions for non-trivial cases, both losses, MSE and QP-loss, will contradict each other.
In such cases, the two losses have different minimizers: for a data point $(\pv, \xv^\star(\pv))$, the MSE-loss is uniquely minimized by $\xv = \xv^\star(\pv)$, whereas $\xv^\star(\pv)$ is not a minimizer of $\mathcal{D}^\mathrm{QP}(\xv, \gamma_g, \gamma_h, \pv)$.
Consequently, no prediction minimizes both losses simultaneously.
\MKHcorrection{For exact penalty functions, this contradiction is avoided only if $\gamma_g$ and $\gamma_h$ exceed the bound on the Lagrange multipliers, see \Cref{sec:penalty_shift}.}

\subsection{Evaluation metrics for hyperparameter tuning}\label{sec:hyperparameter_tuning_metrics}
The objective of hyperparameter tuning is to identify algorithmic parameters that yield optimal model performance according to an evaluation metric $\mathcal{D}^\mathrm{eval}: \R^{n_\theta}\to\R$~\cite{bergstra2011tpe}.
In neural network training, this typically entails tuning learning rates, batch sizes (see~\cite{goodfellow2016deeplearning}), or loss function parameters such as $\underline{\alpha}$, $\overline{\alpha}$, $\gamma_g$, and $\gamma_h$ by optimizing the trained performance.
This process is formally a bilevel optimization problem
\begin{align}
    \min_{\makev{d}_{\mathrm{hyper}}} \quad & \mathcal{D}^\mathrm{eval}(\tv^\star) \\ \label{eq:hyper_opt}
    \text{s.t.} \quad & \tv^\star = \mathcal{A}(\mathcal{L}_{\mathrm{train}}, \makev{d}_{\mathrm{hyper}}),
\end{align}
where $\makev{d}_{\mathrm{hyper}}$ denotes the full set of hyperparameters and $\mathcal{A}$ represents the algorithm that minimizes the training loss $\mathcal{L}_{\mathrm{train}}$.
Note that $\makev{d}_{\mathrm{hyper}}$ may influence both the definition of $\mathcal{L}_{\mathrm{train}}$ and the behavior of $\mathcal{A}$.

We now analyze the suitability of the KKT-loss \eqref{eq:lkkt} versus the \MKHcorrection{penalty loss \eqref{eq:penalty_loss}} as the evaluation metric $\mathcal{D}^\mathrm{eval}$.
The KKT-loss is a robust choice for $\mathcal{D}^\mathrm{eval}$ because\MKHcorrection{, under \Cref{assump:covering},} its global minimum is zero, regardless of the hyperparameter values.
In contrast, the \MKHcorrection{penalty loss} is unsuitable for tuning the penalty weights $\gamma_g$ and $\gamma_h$.
As established previously, the lower bound of the penalty loss is directly dependent on these weights.
Since the penalty terms in~\eqref{eq:penalty_loss} are non-negative, minimizing the \MKHcorrection{penalty loss, for $q=1$ and $q=2$,} with respect to $\makev{d}_{\mathrm{hyper}}$ leads to the trivial solution where $\gamma_g, \gamma_h \to 0$ (or $-\infty$ if $\gamma_g, \gamma_h$ are unbounded).
Consequently, one must rely on alternative strategies: either using a data-driven metric (e.g., MSE on validation data) or heuristically selecting $\gamma_g$ and $\gamma_h$.
The latter means balancing constraint satisfaction with training stability, since large values induce large gradients\MKHcorrection{, as the numerics in~\cite{liu2024teaching} indicate}.
\MKHcorrection{OptINNs also require tuning of hyperparameters, e.g.\ of the schedule bounds $\underline{\alpha}$ and $\overline{\alpha}$, but these can be tuned with the KKT-loss as $\mathcal{D}^\mathrm{eval}$ without requiring additional training data, while the weights $\omega^i$ are determined automatically, see \Cref{sec:weightbalancing}.}

\subsection{Validation, scheduling, and early stopping}\label{sec:validation}
Scheduling describes the active changing of hyperparameters during training, often with regard to previous training results, in order to achieve better performance or accelerate the training~\cite{bishop2023deeplearning}.
One example is the \texttt{ReduceLROnPlateau} learning rate scheduler that reduces the learning rate by a multiplicative factor once the validation loss hits a plateau and stops decreasing.
Again, the lower-boundedness of the KKT-loss allows its application as a validation loss, whereas the \MKHcorrection{penalty loss} is more restricted.
When $\gamma_g$ and $\gamma_h$ are fixed for the whole training, employing it as a validation loss is possible since it is lower-bounded by its minimum, thus disallowing tuning during training.
When $\gamma_g$ and $\gamma_h$ are tuned during training, the \MKHcorrection{penalty loss} does not exhibit the same lower bound for intermediate values of $\gamma_g$ and $\gamma_h$, since \MKHcorrection{its minimum can increase} with rising $\gamma_g$ and $\gamma_h$, see \Cref{fig:penalty_shift_minima}.
This may result in premature reduction of the learning rate, impairing the training.
The same arguments hold when convergence is evaluated on the number of epochs without decreased loss, a technique referred to as early stopping.

\MKHcorrection{%
\subsection{Limitations of the KKT-loss}\label{sec:limitations}
The penalty loss~\eqref{eq:penalty_loss}, for $q=1$ and $q=2$, has two advantages over the KKT-loss \eqref{eq:lkkt}.
First, it does not depend on Lagrange multipliers, thus, it can also hold at points, in which no Lagrange multiplier exists, which satisfy the KKT conditions (cf.~the singular points discussed in \Cref{sec:patched}).
There, the KKT-loss is not zero at the local minimizer, affecting the OptINN training.
Second, it contains the cost function and thus favors local minimizers, whereas the KKT-loss is zero at all KKT-points, e.g.\ also at saddle points.
For OptINNs, the nominal solutions select the solution branch, see \Cref{sec:patched}.
Penalizing violations of the SOSC~\eqref{eq:sosc} in the KKT-loss, e.g.\ via the smallest eigenvalue of the Hessian of the Lagrangian on the tangent space defined by the active inequality and equality constraints, is left for future work.%
}

%% file: figs/penalty_method_optimum_shift.tex
\begin{figure}[htb]
	\centering
	\begin{tikzpicture}
		\begin{axis}[
			width=0.5\linewidth,
			legend pos=north west,
			xlabel={$x$},
			ylabel={$-x + \gamma_g\relu(x-1)^2$},
			legend style={
					/tikz/every even column/.append style={column sep=5mm},
					legend columns=1, 
					legend cell align=left, 
					align=left
				}
			]
			\pgfplotsforeachungrouped \p/\mycolor in {0.5/udsblue, 1/udsred, 2/udsgreen, 5/orange}{
				\pgfmathsetmacro{\xmin}{1/\p+1}
				\pgfmathsetmacro{\ymin}{-\xmin + \p/2*max(0, \xmin-1)^2}
				\edef\plotcommand{
					\noexpand\addplot[
						mark=none,
						color=\mycolor, 
						line width=1pt,
						domain=0:3,
						samples=200
					]
					{-x + \p/2*max(0, x-1)^2} 
					\noexpand; 
					\noexpand\addplot[
						dashed,
						mark=none,
						color=\mycolor, 
						line width=1pt,
						forget plot,
					] coordinates
					{(\xmin, \ymin) (\xmin, 8)}
					\noexpand;
					\noexpand\addplot[
						dotted,
						mark=none,
						color=\mycolor, 
						line width=1pt,
						forget plot,
					] coordinates
					{(0, \ymin) (\xmin, \ymin)}
					\noexpand;
				}
				\plotcommand
				
				\addlegendentryexpanded{$\gamma_g=\p$}
			}
			\addlegendimage{black, line width=1pt}
			\addlegendentry{Constrained minimizer}
			\addplot[
				dashed,
				mark=none,
				color=black, 
				line width=1pt
			] coordinates
			{(1, -1) (1, 8)};
			\addlegendentry{Minimizer}
			\addplot[
				dotted,
				mark=none,
				color=black, 
				line width=1pt
			] coordinates
			{(0, -1) (1, -1)};
		\addlegendentry{Minimum}
		\end{axis}
	\end{tikzpicture}
	\caption{Minima of $-x + \gamma_g\relu(x-1)^2$ under changing $\gamma_g$ for the optimization problem $\min_{x\in\R} -x$ s.t. $x \le 1$. The dashed and dotted lines mark the corresponding minimizer and minimum, respectively. Increasing $\gamma_g$ shifts the minimizer closer to the constrained minimizer $x=1$, but also increases the minimum.
	All suggested minimizers lie in the infeasible region $x>1$}
	\label{fig:penalty_shift_minima}
\end{figure}

%% file: inc/numerical_results.tex
\section{Numerical results}
\label{sec:numerical_results}
\AOcorrection{To compare the performance of OptINNs and QPNNs, we consider four optimization problems of increasing difficulty. We begin with small-scale linear and nonlinear optimization problems ($n_y < 10$) and then move to large-scale linear and nonlinear optimal control problems ($n_y > 100$). For each problem, the training data was generated by formulating the problem in CasADi~\cite{andersson2019casadi} and solving it with IPOPT~\cite{wachter2006ipopt}.
	The parameters of IPOPT describing accuracies were set to $10^{-12}$ to obtain close approximations of the true minimizers.}

\AOcorrection{To ensure comparable expressivity, both OptINNs and QPNNs utilize MLPs with identical depth and uniform width across all hidden layers, where the QPNN implementation follows ~\cite{liu2024teaching, Lillo1993penaltynn}. The sole architectural distinction lies in the final output layer of the OptINN, which contains additional parameters for predicting dual variables. 
Structurally, each linear transformation is followed by layer normalization~\cite{lei2016layernorm} and a ReLU activation. 
The NNs were trained in PyTorch~\cite{paszke2017pytorch} and all parameters are optimized using AdamW~\cite{loshchilov2019adamw} with default settings.
}

The QPNN and OptINN hyperparameters are summarized in \Cref{tab:pmnn_hyperparams,tab:optinn_hyperparams}, respectively.
\AOcorrection{They were tuned with Optuna~\cite{akiba2019optuna} over 100 trials, using as the objective the validation loss averaged over three networks trained with different random seeds.}
Since the QP-loss cannot serve as the hyperparameter tuning metric for QPNNs (see \Cref{sec:hyperparameter_tuning_metrics}), we instead tune their hyperparameters using the mean squared error with respect to training data generated for 256 parameter values sampled from an equidistant grid. This allows a fair comparison between OptINNs and QPNNs with optimized hyperparameters.
\AOcorrection{To reflect settings with only limited solver-generated training data, we fix $\gamma_g$ and $\gamma_h$ in \eqref{eq:penalty_loss} during training, so that validation during training relies solely on the QP-loss.}
\AOcorrection{For OptINNs, both the validation loss and the hyperparameter tuning metric are given by the KKT-loss \eqref{eq:lkkt} with $\omega^i = 1$ and $P^i(x) = |x|$ for $i \in \mathcal{T}$, which requires no additional training data.}
	\AOcorrection{To reduce the dimension of the hyperparameter search space and to ensure that, for $x \approx 0$, constraint satisfaction is weighted at least as strongly as stationarity and complementary slackness, we grouped the penalty functions \eqref{eq:kkt_deviations}.}
\AOcorrection{The penalty function group was tuned only for the first problem, the linear programming problem, and the best-performing group, group 3, was then used for all subsequent problems.}
The optimized hyperparameters for each problem are summarized in Appendix~\ref{app:hparams}.
\begin{table}
    \centering
    \caption{Overview of OptINN hyperparameters.}
    \label{tab:optinn_hyperparams}
    \renewcommand{\arraystretch}{1.2} 
    \begin{tabularx}{\textwidth}{@{}l X@{}}
        \toprule
        \textbf{Category} & \textbf{Parameters / Configurations} \\
        \midrule
        Optimizer & Learning rate, Weight decay \\
        Scheduling & Min/max trade-off ($\underline{\alpha}, \overline{\alpha}$), duration of initialization and annealing phases ($d^\mathrm{init}, d^\mathrm{anneal}$) \\
        Sampling & Number of parameter samples $N$ (per training/validation step) \\
        \addlinespace[1ex]
        Penalty Functions & \textbf{1.} $\Ps(x) = \Pcs(x) = \Ppg(x) = \Pph(x) = |x|$ \\
                          & \textbf{2.} $\Ps(x) = \Pcs(x) = \Ppg(x) = \Pph(x) = x^2$ \\
                          & \textbf{3.} $\Ps(x) = \Pcs(x) = x^2$, $\Ppg(x) = \Pph(x) = |x| + x^2$ \\
                          & \textbf{4.} $\Ps(x) = \Pcs(x) = \Ppg(x) = \Pph(x) = |x| + x^2$ \\
        \bottomrule
    \end{tabularx}
\end{table}

\begin{table}
    \centering
    \caption{Overview of QPNN hyperparameters.}
    \label{tab:pmnn_hyperparams}
    \renewcommand{\arraystretch}{1.2}
    \begin{tabularx}{\textwidth}{@{}l X@{}}
        \toprule
        \textbf{Category} & \textbf{Parameters / Details} \\
        \midrule
        Optimizer & Learning rate, Weight decay set to 0 due to computed value being negligible \\
        Loss Balancing & Trade-off parameter (Data loss vs. Constraint loss) \\
        Constraint Weights & Weights for equality and inequality constraints \\
        Sampling & Number of parameter samples $N$ (per training/validation step) \\
        \bottomrule
    \end{tabularx}
\end{table}

\AOcorrection{After training, we evaluate the predictive performance of the models on $N = 256$ parameter values $\pv$ with corresponding reference solutions $\yv^\star_x(\pv)$. Since QPNNs do not predict dual variables, all metrics depend only on the predicted primal variables. These are the MSE of the primal variables,}
\[
	\namedL{MSE, primal}(\tv) = \frac{1}{N} \sum_{i=1}^N \left\Vert \yv^\star_x(\pv_i) - \N_x(\pv_i, \tv) \right\Vert_2^2,
\]
\MKHcorrection{the mean absolute violations of the inequality and equality constraints
\begin{subequations}
		\begin{align}
			\namedL{IneqViol}(\tv) &= \frac{1}{N n_g}\sum_{i=1}^{N} \left\Vert \relu\left(\gv\left(\N_x(\pv_i, \tv),\pv_i\right)\right) \right\Vert_1, \label{eq:def_L_IneqViol}\\
			\namedL{EqViol}(\tv) &= \frac{1}{N n_h}\sum_{i=1}^{N} \left\Vert \hv\left(\N_x(\pv_i, \tv),\pv_i\right) \right\Vert_1\label{eq:def_L_EqViol},
		\end{align}
\end{subequations}
}%
\MKHcorrection{and the relative cost gap
\begin{equation}
	\delta(\pv_i) = \frac{f\left(\N_x(\pv_i, \tv), \pv_i\right) - f\left(\yv^\star_x(\pv_i), \pv_i\right)}{\max_{j} f\left(\yv^\star_x(\pv_j), \pv_j\right) - \min_{j} f\left(\yv^\star_x(\pv_j), \pv_j\right)}.
	\label{eq:def_delta}
\end{equation}
}%

\MKHcorrection{For $\delta$, we report the mean of $|\delta|$ and the fraction $\delta^-$ of parameter values with $\delta < -10^{-4}$, i.e., the cost of the prediction is below the optimal cost.} \AOcorrection{To assess the robustness of the approaches with respect to randomness, all experiments were repeated with five different random seeds. The seed controls all random components of the training process, including the initialization of the network parameters and the sampling of the parameter values $p$ in each training iteration.}

\subsection{Linear programming}\label{sec:linear_programming}
\AOcorrection{Consider the} linear programming example
\begin{align}
	\min_x &\begin{bmatrix}
		-0.1 & -0.25
	\end{bmatrix}\xv\\
	\st &\begin{bmatrix}
		0.01 & 0.01 \\
		0.04 & 0.12 \\
		0.06 & 0.12 \\
		-0.1 & 0 \\
		0 & -0.1 \\
	\end{bmatrix}\xv \le \begin{bmatrix}
	0.4 \\
	2.4 + \frac{p}{1000}\\
	3.12\\
	0\\
	0	
	\end{bmatrix},
\end{align} 
\AOcorrection{from~\cite{demarchi2023functionapprox} with the closed form solution} 
\[
\left(\xv^\star(p)^\top, \mv^\star(p)^\top\right) =
\begin{cases} 
\left(\begin{pmatrix}
	0 &
	26
\end{pmatrix}, \begin{pmatrix}
	0 & 0 & 2.0833 & 0.25 & 0
\end{pmatrix}\right), & \text{if } 720 \leq p, \\[10pt]
\left(\begin{pmatrix}
	36 - \frac{p}{20} &
	8 + \frac{p}{40}
\end{pmatrix}, \begin{pmatrix}
	0 & 1.25 & 0.833 & 0 & 0
\end{pmatrix}\right), & \text{if } 160 \leq p < 720, \\[10pt]
\left(\begin{pmatrix}
	30 - \frac{p}{80} &
	10 + \frac{p}{80}
\end{pmatrix}, \begin{pmatrix}
	2.5 & 1.875 & 0 & 0 & 0
\end{pmatrix}\right), & \text{if } -800 \leq p < 160, \\[10pt]
\left(\begin{pmatrix}
	60 + \frac{p}{40} &
	0
\end{pmatrix}, \begin{pmatrix}
	0 & 2.5 & 0 & 0 & 0.5
\end{pmatrix}\right), & \text{if } -2400 \leq p < -800, \\[10pt]
\text{empty feasible set}, & \text{if } p < -2400.
\end{cases}
\]
In comparison to~\cite{demarchi2023functionapprox}, the cost and constraints were rescaled to improve numerical stability.
\MKHcorrection{The open intervals between the values in $\mathcal{Z} = \{-2400, -800, 160, 720\}$ and the interval $(720, 2400]$ correspond to the neighborhoods of parameter values as defined in~\Cref{sec:patched}, giving $\bar{\mathcal{P}} = \{p \in \R\ \vert\ -2400 \le p \le 2400\}$.
At each value in $\mathcal{Z}$, three inequality constraints are active, so that LICQ is violated and these are singular points, see \Cref{sec:patched}.}
\AOcorrection{For each of the four neighborhoods}, one value was chosen for training, giving four training data points corresponding to $\hat{p}_i \in \{-1500, -300, 400, 1500\}$ in total.

\Cref{fig:linear_demarchi_nn} shows the results for training an OptINN and a QPNN on this problem using the four data points and an OptINN using no data points\MKHcorrection{ as well as an MLP trained with the MSE-loss only, i.e.\ using the cost function \eqref{eq:combined_loss} with $\alpha = 0$,} for reference.
\input{figs/linear_demarchi_results_optinn_vs_penalty}
\MKHcorrection{Except for the MLP, all models} are able to predict the true values.
\MKHcorrection{The MLP trained with the MSE-loss only overfits the four training data points} \AOcorrection{and is not able to generalize.}
\MKHcorrection{Both} the QPNN and the OptINNs match the reference solution well, even though the OptINNs have in the input and hidden layers the same number of trainable parameters as the QPNN.
\MKHcorrection{With a similar mean squared error in the primal variables, the OptINN additionally predicts the dual variables with a mean squared error of $6.0\cdot10^{-2}$.
The error of the QPNN concentrates at the singular points: $61\,\%$ of its squared error in the primal variables arises within a distance of 50 to the values in $\mathcal{Z}$, which make up $7\,\%$ of the evaluated parameter values, compared to $30\,\%$ for the OptINN.}
\AOcorrection{The predictions of the OptINNs violate the inequality constraints less than those of the QPNNs, as indicated by the lower values of the mean inequality violation $\namedL{IneqViol}$ \eqref{eq:def_L_IneqViol} in \Cref{tab:numerical_summary}.}.
\MKHcorrection{
The difference is most pronounced in $\delta^-$: for $70\,\%$ of the parameter values, the cost of the QPNN prediction is below the optimal cost, i.e.\ the prediction lies outside of the feasible region, as discussed in from \Cref{sec:penalty_shift}.}
\AOcorrection{For the OptINNs, the predicted dual variables change continuously across the points where the active constraint switches, whereas the exact dual variables are discontinuous there.
	We attribute this to the Softplus activation, which yields strictly positive dual variables. Complementary slackness is therefore only approximately satisfied, and the dual variables of inactive constraints remain non-zero near the switching points.}

\Cref{fig:linear_demarchi_convergence} shows the training curves for both the OptINN trained with four data points and the OptINN trained without data points.
\MKHcorrection{From \Cref{tab:numerical_summary}, it can be concluded that with an OptINN with MSE loss yields better results and reduces the spread of the primal error over the five trainings.
}

\input{figs/linear_demarchi_convergence.tex}


\subsection{Nonlinear problem with nonconvex inequality constraints}\label{sec:quadratic_nonl}
\AOcorrection{As a second example, we consider the Euclidean projection of a parameter point onto a nonconvex set. Let $f \colon \R \to \R$ be the polynomial $f(x) = 2x^3 + 3x^2 + 2x + 1$, and let the parameter region of interest be $\bar{\mathcal{P}} = \{\pv \in \R^2\ \vert\ \Vert\pv\Vert_\infty \le 1\}$. For a given $\pv \in \bar{\mathcal{P}}$, we seek the point $\makev{z} = \begin{bmatrix} z_1 & z_2 \end{bmatrix}^\top \in \R^2$ closest to $\pv$ that satisfies four nonconvex constraints, which leads to the optimization problem
	\begin{equation}\label{eq:nonconvex_example}
		\begin{aligned}
			\min_{\makev{z} \in \R^2} \quad & \Vert \makev{z} - \pv \Vert_2^2 \\
			\st \quad & z_2 \le f(z_1), \\
			& -z_2 \le f(z_1), \\
			& z_2 \le f(-z_1), \\
			& -z_2 \le f(-z_1).
		\end{aligned}
\end{equation}}
\AOcorrection{\Cref{fig:quadratic_nonl_problem_results} shows the feasible set, which is the region enclosed by the boundaries of the four constraints.
If $\pv$ lies in the feasible set, the unique minimizer is $\makev{z}^\star = \pv$.
Otherwise, $\makev{z}^\star$ lies on the boundary of the feasible set. If only one constraint is active at $\makev{z}^\star$, the vector $\pv - \makev{z}^\star$ is orthogonal to the tangent of that constraint's boundary curve at $\makev{z}^\star$.

For both the QPNN and the OptINN, \Cref{fig:quadratic_nonl_problem_results} shows the training inputs, the test inputs, and the corresponding test predictions.
Again, both methods learn the overall structure of the patched solution: for $\pv$ inside the feasible set, the predictions nearly coincide with $\pv$, whereas for $\pv$ outside, they are mapped onto the boundary of the feasible set.
Here, the OptINN outperforms the QPNN in all metrics of \Cref{tab:numerical_summary}, most notably in the violation of the inequality constraints, which is lower by a factor of 4.5.}


\input{figs/quadratic_nonl_solution_diff_kkt}

\subsection{Linear optimal control problem -- Double integrator}\label{sec:rocketcar}
\AOcorrection{Next, we present the results for the optimal control of the double integrator $\ddot{z}(t) = u(t)$, with state $\xv = {[z, v]}^\top \in \R^2$, consisting of position $z$ and velocity $v = \dot{z}$, and input $u \in \R$, the acceleration.
	The task is to move the system from rest at the origin to rest at the target position $p \in \R$ within the fixed time horizon $t_f = \SI{12.8}{\second}$, while minimizing the control effort $\int_0^{t_f} u^2(t) \,\mathrm{d}t$ and respecting the input bounds $-1 \le u(t) \le 1$.
	The target position $p$ is the parameter of the problem.
	
	We solve this task with full discretization, a direct optimal control method. The horizon is divided into $N = 32$ intervals of length $\Delta t = t_f / N = \SI{0.4}{\second}$, and the input is assumed to be piecewise constant on each interval (zero-order hold). Under this assumption, integrating the dynamics over one interval gives the exact discrete-time system $\xv[k+1] = \makev{A}_d \xv[k] + \makev{b}_d u[k]$ with
	\begin{equation}
		\makev{A}_d = \begin{bmatrix} 1 & \Delta t \\ 0 & 1 \end{bmatrix}, \qquad
		\makev{b}_d = \begin{bmatrix} \frac{\Delta t^2}{2} \\ \Delta t \end{bmatrix},
	\end{equation}
	where $\xv[k]$ denotes the state at time $t = k \Delta t$ and $u[k]$ the input on the interval $[k\Delta t, (k+1)\Delta t)$. Approximating the integral in the cost by the corresponding sum yields the finite-dimensional optimization problem
	\begin{equation}\label{eq:double_integrator}
		\begin{aligned}
			\min_{\xv[0], \ldots, \xv[N],\, u[0], \ldots, u[N-1]} \quad & \Delta t \sum_{k=0}^{N-1} u[k]^2 \\
			\st \quad & \xv[k+1] = \makev{A}_d \xv[k] + \makev{b}_d u[k], && k = 0, \ldots, N-1, \\
			& \xv[0] = \begin{bmatrix} 0 \\ 0 \end{bmatrix}, \quad \xv[N] = \begin{bmatrix} p \\ 0 \end{bmatrix}, \\
			& -1 \le u[k] \le 1, && k = 0, \ldots, N-1.
		\end{aligned}
\end{equation}}

where \(\cdot[k]\) indicates the discretized quantity at the $k$-th point in time $t_k = k\Delta t$\MKHcorrection{, with $k\in\{0,\dots,N\}$ for the state and $k\in\{0,\dots,N-1\}$ for the input}.
The terminal position $p \in \bar{\mathcal{P}} = \{p \in \R\ \vert\ 0 \le p \le 40\}$ is the parameter in this problem, giving $n_x=98$, $n_p=1$, $n_g=64$, and $n_h=68$. 
For more information regarding optimal control refer to~\cite{gerdts2012optimalcontrol}.

\AOcorrection{For this experiment, we generated training data for the MSE loss at the three parameter values $\hat{p}_i \in \{0, 20, 40\}$.
In contrast to the previous examples, these samples do not satisfy \Cref{assump:covering}.
The magnitude of the optimal input grows with $p$, and for $p \ge 28.16$, its largest entries reach the input bounds $\pm 1$. 
As $p$ increases further, more and more entries of $u$ saturate, so the set of active inequality constraints changes repeatedly.
The region $p \ge 28.16$ therefore consists of many neighborhoods with different active sets, which cannot all be covered by the single sample $\hat{p} = 40$. Since, except for $p=40$, these cases are not part of the training data, the QPNN and OptINN have to learn the activity of the constraints from the QP-loss and the KKT-loss, respectively.
}

\AOcorrection{\Cref{fig:heatmap_rocketcar_diff} shows the accuracy of the OptINN and the QPNN.
For each decision variable, i.e., position $z$, velocity $v$, and acceleration $u$, a heatmap shows the deviation of the prediction from the reference solution computed with IPOPT, as a function of the parameter $p$ and the time step $k$.
Both methods achieve low errors also for parameter values far from the training samples, indicating that they generalize rather than merely memorize the training data.
As expected from the violation of \Cref{assump:covering}, the errors increase for $p \gtrsim 28$, where the input bounds become active. In this region, the OptINN achieves lower errors than the QPNN.}
\input{figs/rocket_car_heatmaps_sparse_data}
\Cref{fig:rocket_car_1d} shows the results for \MKHcorrection{$p=39.6$, the value with the highest deviation between OptINN prediction and reference solution.}
Both methods generally follow the reference solution well. However, even at its worst-performing input value, the OptINN produces lower errors than the QPNN, especially regarding acceleration.
\input{figs/rocket_car_1d_results}

\AOcorrection{\Cref{fig:rocketcar_difference} shows, for $p \in [0, 40]$, the mean, minimum, and maximum over the five training runs of the deviation between the predicted and the reference primal variables.}
\AOcorrection{
We find that the OptINN achieves lower errors than the QPNN \MKHcorrection{for all parameter values.
The error of the OptINN remains small up to $p \approx 28$, increases once the box constraints on $u$ become active, and decreases again near the training point $p = 40$.
The error of the QPNN, in contrast, increases steadily with $p$ and is largest at $p = 40$, even though this is a training point.}
Owing to the trivialization layer, the OptINN satisfies the box constraints on $u$ by construction for all parameter values.
\MKHcorrection{Overall, the OptINN outperforms the QPNN in all metrics of \Cref{tab:numerical_summary}.
The largest difference occurs in the dynamics: the mean violation of the equality constraints is 18 times lower for the OptINN than for the QPNN.
The larger violation of the QPNN also explains why its predicted cost lies below the optimal cost for $94\,\%$ of the parameter values ($\delta^-$ in \Cref{tab:numerical_summary}). Such a cost can only be attained by predictions that violate the constraints, so it indicates an error rather than an advantage.}}
\input{figs/rocketcar_diff_primal_variables.tex}

\newpage
\subsection{Nonlinear optimal control problem -- Pendulum swing-up}\label{sec:pendulum}
\AOcorrection{Finally, we consider the classical problem of swinging up a pendulum by an applied torque.
The pendulum consists of a point mass $m$ attached to a massless rod of length $l$ and is governed by the nonlinear ODE
\begin{equation}
	m l^2 \ddot{\varphi}(t) = \tau(t) - m g l \sin\bigl(\varphi(t)\bigr),
\end{equation}
where $\varphi$ is the angle measured from the downward position, $\tau$ the applied torque, and $g$ the gravitational acceleration.
With the state $\xv = {[\varphi, \dot{\varphi}]}^\top \in \R^2$ and the input $\tau \in \R$, the dynamics take the form $\dot{\xv}(t) = \makev{f}\bigl(\xv(t), \tau(t)\bigr)$ with
\begin{equation}
	\makev{f}(\xv, \tau) = \begin{bmatrix} \xv^2 \\ -\frac{g}{l} \sin(\xv^1) \end{bmatrix} + \begin{bmatrix} 0 \\ \frac{1}{m l^2} \end{bmatrix} \tau.
\end{equation}
The task is to move the pendulum from rest in the downward position to rest in the upright position within the time horizon $p$, while minimizing the control effort and respecting the torque bounds $-2 \le \tau(t) \le 2$. In contrast to the double integrator, the parameter of the problem is now the final time $p$, which leads to the optimal control problem
\begin{equation}\label{eq:pendulum}
	\begin{aligned}
		\min_{\tau(\cdot)} \quad & \int_0^{p} \tau(t)^2 \,\mathrm{d}t \\
		\st \quad & \dot{\xv}(t) = \makev{f}\bigl(\xv(t), \tau(t)\bigr), && t \in [0, p], \\
		& \xv(0) = \begin{bmatrix} 0 \\ 0 \end{bmatrix}, \quad \xv(p) = \begin{bmatrix} \pi \\ 0 \end{bmatrix}, \\
		& -2 \le \tau(t) \le 2, && t \in [0, p].
	\end{aligned}
\end{equation}}
\MKHcorrection{As in the previous example, the input is held constant over each time step.
While this yields the exact state transition for the double integrator, we approximate the state transition of the nonlinear pendulum by} the classical Runge-Kutta method of order 4 ($\rk$) (cf.~\cite{gerdts2012optimalcontrol}).
For the resulting optimization problem 
\begin{align}
	\min_{\xv,\tau} \ \Delta t(p) \sum_{k=0}^{N-1}& {\tau[k]}^2\\
	\st \ \xv[k+1] &= \rk(\xv[k], \tau[k], \Delta t(p)) \ \forall k\in\{0, \dots, N-1\}\\
	 \xv[0] &= \begin{bmatrix}
		0 \\ 0
	\end{bmatrix}\\
	 \xv[N] &= \begin{bmatrix}
	\pi \\ 0
	\end{bmatrix}\\
	-2 &\le \tau[k] \le 2 \ \forall k\in\{0, \dots, N-1\},
\end{align}
the time step \(\Delta t(p) = \frac{p}{N}\) depends on the parameter $p\in\mathcal{\bar{P}} = \{p\in\R\ \vert\ 6\le p \le 15\}$, the time to reach the upright position.
\AOcorrection{In contrast to the double integrator, we do not impose the boundary conditions as equality constraints.
Instead, we fix $\xv[0]$ and $\xv[N]$ to their prescribed values and remove them from the decision variables, which reduces the number of equality constraints by four.
With $N = 100$ time steps, the resulting problem has $n_x = 298$ decision variables, $n_p = 1$ parameter, $n_g = 200$ inequality constraints, and $n_h = 200$ equality constraints.
This formulation corresponds to an $\varepsilon$-constraint scalarization (see, e.g.,~\cite{Marler2004}) of the multi-objective optimal control problem of finding optimal trade-offs between minimum control effort and minimum transition time, with the transition time prescribed by the parameter $p$.}

Unlike the previous linear-quadratic (and thus convex) problem, the swing-up of the pendulum has multiple local minimizers.
\MKHcorrection{Its local minimizers differ in the number of swings, and which of them has the lowest cost changes with $p$.
Since IPOPT converges to a local minimizer that depends on its initial guess, we started IPOPT for each of the 256 validation parameter values from 22 initial guesses, among them swing-ups with 0 to 20 half swings, and from the solutions at neighboring parameter values, and kept the solution with the lowest cost.
This does not guarantee a global minimizer, but warm-starting IPOPT from the predictions of all trained networks yields no lower cost.
\AOcorrection{For an OptINN, this expensive search is needed only at the few training parameter values, cf.\ \Cref{sec:patched}, whereas a purely data-driven network requires it at every training point.}
The trained network should then predict the solutions from the same branch, so that deploying the optimal control requires a single evaluation of the network instead of expensive warm-started solves of IPOPT.}

\AOcorrection{Both the QPNN and the OptINN rely on optimality conditions that only characterize local minimizers. Since the pendulum problem has several local minimizers with comparable cost, a prediction may approximate a different local minimizer than the reference solution.
It would then show a large deviation in the primal variables, even though it is an equally good solution.
We therefore do not evaluate the deviation in the primal variables.
Instead, we assess how close the predicted cost is to the cost of the reference solution, using the relative error, and how strongly the predictions violate the equality constraints, using the mean absolute error.}
\MKHcorrection{Due to the trivialization layer, the OptINN satisfies the inequality constraints per construction, whereas the QPNN violates them.
	
\Cref{fig:pendulum_example} shows the violation of equality constraints (a), the cost of the prediction (b), and the relative cost gap $\delta$ defined in \eqref{eq:def_delta} (c) for a given transition time $p$.}
\MKHcorrection{For almost all values of $p$, the OptINN produces a lower mean absolute equality constraint violation $\namedL{EqViol}$ defined in \eqref{eq:def_L_IneqViol} than the QPNN, i.e.\ it better adheres to the dynamics and reaches closer to the goal.}
\MKHcorrection{At the training points, the cost of the QPNN prediction closely matches the optimal cost. Between them, however, it deviates by up to $+21\,\%$ and $-1.2\,\%$, measured relative to the range of the optimal cost as in the definition of $\delta$.
The cost of the OptINN prediction stays close to the optimal cost over the whole parameter range, with deviations of at most $+3.4\,\%$ and $-1.8\,\%$. 
\AOcorrection{The largest deviations occur near the eight kinks of the optimal cost, which are spaced roughly $1.2$ apart between $p = 6.7$ and $p = 14.9$.}
Predicted costs below the optimal cost occur more than twice as often for the QPNN as for the OptINN, see $\delta^-$ in \Cref{tab:numerical_summary}.}

\input{figs/pendulum_cost_comparison}
\input{tables/numeric_metrics.tex}

\subsection{Summary}
\MKHcorrection{In our numerical experiments on small-scale optimization problems, OptINNs predicted the primal variables with an accuracy comparable to or better than that of QPNNs while better satisfying the constraints, with the added advantage of providing estimates for dual variables.
For the larger-scale problems, the OptINN clearly outperforms the QPNN in all evaluated metrics.}

Inspecting the tables in Appendix~\ref{app:hparams}, we observe that the tuned hyperparameters for QPNNs are significantly more extreme than those for OptINNs.
The disparity is most pronounced in the trade-off parameter $\alpha$, which balances the data-based loss against the penalty term. 
While the bounds $\underline{\alpha}$ and $\overline{\alpha}$ for the OptINN span a broad range of moderate magnitude, the QPNN necessitates an increasingly dominant weighting of the data-based loss.

%% file: figs/linear_demarchi_results_optinn_vs_penalty.tex

\def\optinnFile{data/results/lp/optinn/lp_optinn_predictions.csv}
\def\penaltyFile{data/results/lp/qpnn/lp_qpnn_predictions.csv}
\def\nodataFile{data/results/lp_nodata/optinn/lp_nodata_optinn_predictions.csv}
\def\mseonlyFile{data/results/lp_mseonly/optinn/lp_mseonly_optinn_predictions.csv}
\definecolor{mseonlycolor}{RGB}{148,103,189}
\begin{figure}[htb]
	\centering
	\begin{tikzpicture}
		\begin{groupplot}[
			group style={
				group size=3 by 3,
				horizontal sep=1.2cm,
				vertical sep=1cm
			},
			ylabel shift=-4pt,
			scale only axis,
			xmin=-2400,
			xmax=2400,
			axis background/.style={fill=white},
			height=0.18\linewidth,
			width=0.25\linewidth,
			enlarge x limits=0.05,
            grid=major,
			]
			\nextgroupplot[
			ylabel={${x}^{1}$},
			]
			\addplot [color=udsred, line width=1.0pt, forget plot]
			table[x=input, y=x0, col sep=comma] {data/results/lp/linear_problem_data_resolution.csv};
			
			\addplot [dashed, color=udsgreen, line width=1.0pt, forget plot]
			table[x=p, y=x1_mean, col sep=comma] {\optinnFile};
			
			\addplot [dashed, name path=upper_optinn, color=udsgreen, line width=1.0pt, forget plot, draw=none]
			table[x=p, y=x1_max, col sep=comma] {\optinnFile};
			
			\addplot [dashed, name path=lower_optinn, color=udsgreen, line width=1.0pt, forget plot, draw=none]
			table[x=p, y=x1_min, col sep=comma] {\optinnFile};
			
			\addplot [fill=udsgreen, draw=none, opacity=0.2]
				fill between[
					of=upper_optinn and lower_optinn
			];
			
			\addplot [dashdotted, color=udsblue, line width=1.0pt, forget plot]
			table[x=p, y=x1_mean, col sep=comma] {\penaltyFile};
			
			\addplot [dashdotted, name path=upper_penalty, color=udsblue, line width=1.0pt, forget plot, draw=none]
			table[x=p, y=x1_max, col sep=comma] {\penaltyFile};
			
			\addplot [dashdotted, name path=lower_penalty, color=udsblue, line width=1.0pt, forget plot, draw=none]
			table[x=p, y=x1_min, col sep=comma] {\penaltyFile};
			
			\addplot [fill=udsblue, draw=none, opacity=0.2]
			fill between[
			of=upper_penalty and lower_penalty
			];

			\addplot [dotted, color=black, line width=1.0pt, forget plot]
			table[x=p, y=x1_mean, col sep=comma] {\nodataFile};

			\addplot [dotted, name path=upper_no_data, color=black, line width=1.0pt, forget plot, draw=none]
			table[x=p, y=x1_max, col sep=comma] {\nodataFile};

			\addplot [dotted, name path=lower_no_data, color=black, line width=1.0pt, forget plot, draw=none]
			table[x=p, y=x1_min, col sep=comma] {\nodataFile};

			\addplot [fill=black, draw=none, forget plot, opacity=0.2]
			fill between[
			of=upper_no_data and lower_no_data
			];

			\addplot [dash dot dot, color=mseonlycolor, line width=1.0pt, forget plot]
			table[x=p, y=x1_mean, col sep=comma] {\mseonlyFile};

			\addplot [dash dot dot, name path=upper_mse_only, color=mseonlycolor, line width=1.0pt, forget plot, draw=none]
			table[x=p, y=x1_max, col sep=comma] {\mseonlyFile};

			\addplot [dash dot dot, name path=lower_mse_only, color=mseonlycolor, line width=1.0pt, forget plot, draw=none]
			table[x=p, y=x1_min, col sep=comma] {\mseonlyFile};

			\addplot [fill=mseonlycolor, draw=none, forget plot, opacity=0.2]
			fill between[
			of=upper_mse_only and lower_mse_only
			];
			
			\addplot [color=orange, line width=1.0pt, only marks, mark=x, mark options={solid, orange}]
			table[x=input, y=x0, col sep=comma] {data/results/lp/linear_problem_data_optinn_trained_on.csv};
			
			\nextgroupplot[
			ymin=-0.330680524518895,
			ymax=30,
			ylabel={${x}^{2}$},
			legend to name=commonlegendLinNN, 
			legend style={
				/tikz/every even column/.append style={column sep=5mm},
				legend columns=3, 
				at={(0.03,0.97)}, 
				anchor=north west, 
				legend cell align=left, 
				align=left, 
				draw=white!15!black
			}
			]
			\addplot [color=udsred, line width=1.0pt]
			table[x=input, y=x1, col sep=comma] {data/results/lp/linear_problem_data_resolution.csv};
			\addlegendentry{Numerical optimization}
			
			\addplot [dashed, color=udsgreen, line width=1.0pt]
			table[x=p, y=x2_mean, col sep=comma] {\optinnFile};
			\addlegendentry{OptINN}
			
			\addplot [dashed, name path=upper_optinn, color=udsgreen, line width=1.0pt, forget plot, draw=none]
			table[x=p, y=x2_max, col sep=comma] {\optinnFile};
			
			\addplot [dashed, name path=lower_optinn, color=udsgreen, line width=1.0pt, forget plot, draw=none]
			table[x=p, y=x2_min, col sep=comma] {\optinnFile};
			
			\addplot [fill=udsgreen, draw=none, forget plot, opacity=0.2]
			fill between[
			of=upper_optinn and lower_optinn
			];
			
			\addplot [dashdotted, color=udsblue, line width=1.0pt]
			table[x=p, y=x2_mean, col sep=comma] {\penaltyFile};
			\addlegendentry{QPNN}
			
			\addplot [dashdotted, name path=upper_penalty, color=udsblue, line width=1.0pt, forget plot, draw=none]
			table[x=p, y=x2_max, col sep=comma] {\penaltyFile};
			
			\addplot [dashdotted, name path=lower_penalty, color=udsblue, line width=1.0pt, forget plot, draw=none]
			table[x=p, y=x2_min, col sep=comma] {\penaltyFile};
			
			\addplot [fill=udsblue, draw=none, forget plot, , opacity=0.2]
			fill between[
			of=upper_penalty and lower_penalty
			];

			\addplot [dotted, color=black, line width=1.0pt]
			table[x=p, y=x2_mean, col sep=comma] {\nodataFile};
			\addlegendentry{OptINN (no data)}

			\addplot [dotted, name path=upper_no_data, color=black, line width=1.0pt, forget plot, draw=none]
			table[x=p, y=x2_max, col sep=comma] {\nodataFile};

			\addplot [dotted, name path=lower_no_data, color=black, line width=1.0pt, forget plot, draw=none]
			table[x=p, y=x2_min, col sep=comma] {\nodataFile};

			\addplot [fill=black, draw=none, forget plot, opacity=0.2]
			fill between[
			of=upper_no_data and lower_no_data
			];

			\addplot [dash dot dot, color=mseonlycolor, line width=1.0pt]
			table[x=p, y=x2_mean, col sep=comma] {\mseonlyFile};
			\addlegendentry{MLP (MSE only)}

			\addplot [dash dot dot, name path=upper_mse_only, color=mseonlycolor, line width=1.0pt, forget plot, draw=none]
			table[x=p, y=x2_max, col sep=comma] {\mseonlyFile};

			\addplot [dash dot dot, name path=lower_mse_only, color=mseonlycolor, line width=1.0pt, forget plot, draw=none]
			table[x=p, y=x2_min, col sep=comma] {\mseonlyFile};

			\addplot [fill=mseonlycolor, draw=none, forget plot, opacity=0.2]
			fill between[
			of=upper_mse_only and lower_mse_only
			];
			
			\addplot [color=orange, line width=1.0pt, only marks, mark=x, mark options={solid, orange}]
			table[x=input, y=x1, col sep=comma] {data/results/lp/linear_problem_data_optinn_trained_on.csv};
			\addlegendentry{Training data}
			
			\nextgroupplot[
			ymin=-0.4,
			ymax=3.0,
			ylabel={$\mu{}^{1}$},
			]
			\addplot [color=udsred, line width=1.0pt, forget plot]
			table[x=input, y=l0, col sep=comma] {data/results/lp/linear_problem_data_resolution.csv};
			
			\addplot [dashed, color=udsgreen, line width=1.0pt, forget plot]
			table[x=p, y=mu1_mean, col sep=comma] {\optinnFile};
			
			\addplot [dashed, name path=upper_optinn, color=udsgreen, line width=1.0pt, forget plot, draw=none]
			table[x=p, y=mu1_max, col sep=comma] {\optinnFile};
			
			\addplot [dashed, name path=lower_optinn, color=udsgreen, line width=1.0pt, forget plot, draw=none]
			table[x=p, y=mu1_min, col sep=comma] {\optinnFile};
			
			\addplot [fill=udsgreen, draw=none, opacity=0.2]
			fill between[
			of=upper_optinn and lower_optinn
			];

			\addplot [dotted, color=black, line width=1.0pt, forget plot]
			table[x=p, y=mu1_mean, col sep=comma] {\nodataFile};

			\addplot [dotted, name path=upper_no_data, color=black, line width=1.0pt, forget plot, draw=none]
			table[x=p, y=mu1_max, col sep=comma] {\nodataFile};

			\addplot [dotted, name path=lower_no_data, color=black, line width=1.0pt, forget plot, draw=none]
			table[x=p, y=mu1_min, col sep=comma] {\nodataFile};

			\addplot [fill=black, draw=none, forget plot, opacity=0.2]
			fill between[
			of=upper_no_data and lower_no_data
			];

			\addplot [dash dot dot, color=mseonlycolor, line width=1.0pt, forget plot]
			table[x=p, y=mu1_mean, col sep=comma] {\mseonlyFile};

			\addplot [dash dot dot, name path=upper_mse_only, color=mseonlycolor, line width=1.0pt, forget plot, draw=none]
			table[x=p, y=mu1_max, col sep=comma] {\mseonlyFile};

			\addplot [dash dot dot, name path=lower_mse_only, color=mseonlycolor, line width=1.0pt, forget plot, draw=none]
			table[x=p, y=mu1_min, col sep=comma] {\mseonlyFile};

			\addplot [fill=mseonlycolor, draw=none, forget plot, opacity=0.2]
			fill between[
			of=upper_mse_only and lower_mse_only
			];
			
			\addplot [color=orange, line width=1.0pt, only marks, mark=x, mark options={solid, orange}]
			table[x=input, y=l0, col sep=comma] {data/results/lp/linear_problem_data_optinn_trained_on.csv};
			
			\nextgroupplot[
			ymin=-0.123496284707616,
			ymax=3,
			ylabel={$\mu{}^{2}$},
			xlabel style={font=\color{white!15!black}},
			xlabel={$p$},
			]
			\addplot [color=udsred, line width=1.0pt, forget plot]
			table[x=input, y=l1, col sep=comma] {data/results/lp/linear_problem_data_resolution.csv};
			
			\addplot [dashed, color=udsgreen, line width=1.0pt, forget plot]
			table[x=p, y=mu2_mean, col sep=comma] {\optinnFile};
			
			\addplot [dashed, name path=upper_optinn, color=udsgreen, line width=1.0pt, forget plot, draw=none]
			table[x=p, y=mu2_max, col sep=comma] {\optinnFile};
			
			\addplot [dashed, name path=lower_optinn, color=udsgreen, line width=1.0pt, forget plot, draw=none]
			table[x=p, y=mu2_min, col sep=comma] {\optinnFile};
			
			\addplot [fill=udsgreen, draw=none, opacity=0.2]
			fill between[
			of=upper_optinn and lower_optinn
			];

			\addplot [dotted, color=black, line width=1.0pt, forget plot]
			table[x=p, y=mu2_mean, col sep=comma] {\nodataFile};

			\addplot [dotted, name path=upper_no_data, color=black, line width=1.0pt, forget plot, draw=none]
			table[x=p, y=mu2_max, col sep=comma] {\nodataFile};

			\addplot [dotted, name path=lower_no_data, color=black, line width=1.0pt, forget plot, draw=none]
			table[x=p, y=mu2_min, col sep=comma] {\nodataFile};

			\addplot [fill=black, draw=none, forget plot, opacity=0.2]
			fill between[
			of=upper_no_data and lower_no_data
			];

			\addplot [dash dot dot, color=mseonlycolor, line width=1.0pt, forget plot]
			table[x=p, y=mu2_mean, col sep=comma] {\mseonlyFile};

			\addplot [dash dot dot, name path=upper_mse_only, color=mseonlycolor, line width=1.0pt, forget plot, draw=none]
			table[x=p, y=mu2_max, col sep=comma] {\mseonlyFile};

			\addplot [dash dot dot, name path=lower_mse_only, color=mseonlycolor, line width=1.0pt, forget plot, draw=none]
			table[x=p, y=mu2_min, col sep=comma] {\mseonlyFile};

			\addplot [fill=mseonlycolor, draw=none, forget plot, opacity=0.2]
			fill between[
			of=upper_mse_only and lower_mse_only
			];
			
			\addplot [color=orange, line width=1.0pt, only marks, mark=x, mark options={solid, orange}]
			table[x=input, y=l1, col sep=comma] {data/results/lp/linear_problem_data_optinn_trained_on.csv};
			
			\nextgroupplot[
			ymin=-0.0390108262490451,
			ymax=2.20250808028335,
			ylabel={$\mu{}^{3}$},
			xlabel style={font=\color{white!15!black}},
			xlabel={$p$},
			]
			\addplot [color=udsred, line width=1.0pt, forget plot]
			table[x=input, y=l2, col sep=comma] {data/results/lp/linear_problem_data_resolution.csv};
			
			\addplot [dashed, color=udsgreen, line width=1.0pt, forget plot]
			table[x=p, y=mu3_mean, col sep=comma] {\optinnFile};
			
			\addplot [dashed, name path=upper_optinn, color=udsgreen, line width=1.0pt, forget plot, draw=none]
			table[x=p, y=mu3_max, col sep=comma] {\optinnFile};
			
			\addplot [dashed, name path=lower_optinn, color=udsgreen, line width=1.0pt, forget plot, draw=none]
			table[x=p, y=mu3_min, col sep=comma] {\optinnFile};
			
			\addplot [fill=udsgreen, draw=none, opacity=0.2]
			fill between[
			of=upper_optinn and lower_optinn
			];

			\addplot [dotted, color=black, line width=1.0pt, forget plot]
			table[x=p, y=mu3_mean, col sep=comma] {\nodataFile};

			\addplot [dotted, name path=upper_no_data, color=black, line width=1.0pt, forget plot, draw=none]
			table[x=p, y=mu3_max, col sep=comma] {\nodataFile};

			\addplot [dotted, name path=lower_no_data, color=black, line width=1.0pt, forget plot, draw=none]
			table[x=p, y=mu3_min, col sep=comma] {\nodataFile};

			\addplot [fill=black, draw=none, forget plot, opacity=0.2]
			fill between[
			of=upper_no_data and lower_no_data
			];

			\addplot [dash dot dot, color=mseonlycolor, line width=1.0pt, forget plot]
			table[x=p, y=mu3_mean, col sep=comma] {\mseonlyFile};

			\addplot [dash dot dot, name path=upper_mse_only, color=mseonlycolor, line width=1.0pt, forget plot, draw=none]
			table[x=p, y=mu3_max, col sep=comma] {\mseonlyFile};

			\addplot [dash dot dot, name path=lower_mse_only, color=mseonlycolor, line width=1.0pt, forget plot, draw=none]
			table[x=p, y=mu3_min, col sep=comma] {\mseonlyFile};

			\addplot [fill=mseonlycolor, draw=none, forget plot, opacity=0.2]
			fill between[
			of=upper_mse_only and lower_mse_only
			];
			
			\addplot [color=orange, line width=1.0pt, only marks, mark=x, mark options={solid, orange}]
			table[x=input, y=l2, col sep=comma] {data/results/lp/linear_problem_data_optinn_trained_on.csv};
			
			\nextgroupplot[
			ymin=-0.04,
			ymax=0.3,
			ylabel={$\mu{}^{4}$},
			xlabel style={font=\color{white!15!black}},
			xlabel={$p$},
			]
			\addplot [color=udsred, line width=1.0pt, forget plot]
			table[x=input, y=l3, col sep=comma] {data/results/lp/linear_problem_data_resolution.csv};
			
			\addplot [dashed, color=udsgreen, line width=1.0pt, forget plot]
			table[x=p, y=mu4_mean, col sep=comma] {\optinnFile};
			
			\addplot [dashed, name path=upper_optinn, color=udsgreen, line width=1.0pt, forget plot, draw=none]
			table[x=p, y=mu4_max, col sep=comma] {\optinnFile};
			
			\addplot [dashed, name path=lower_optinn, color=udsgreen, line width=1.0pt, forget plot, draw=none]
			table[x=p, y=mu4_min, col sep=comma] {\optinnFile};
			
			\addplot [fill=udsgreen, draw=none, opacity=0.2]
			fill between[
			of=upper_optinn and lower_optinn
			];

			\addplot [dotted, color=black, line width=1.0pt, forget plot]
			table[x=p, y=mu4_mean, col sep=comma] {\nodataFile};

			\addplot [dotted, name path=upper_no_data, color=black, line width=1.0pt, forget plot, draw=none]
			table[x=p, y=mu4_max, col sep=comma] {\nodataFile};

			\addplot [dotted, name path=lower_no_data, color=black, line width=1.0pt, forget plot, draw=none]
			table[x=p, y=mu4_min, col sep=comma] {\nodataFile};

			\addplot [fill=black, draw=none, forget plot, opacity=0.2]
			fill between[
			of=upper_no_data and lower_no_data
			];

			\addplot [dash dot dot, color=mseonlycolor, line width=1.0pt, forget plot]
			table[x=p, y=mu4_mean, col sep=comma] {\mseonlyFile};

			\addplot [dash dot dot, name path=upper_mse_only, color=mseonlycolor, line width=1.0pt, forget plot, draw=none]
			table[x=p, y=mu4_max, col sep=comma] {\mseonlyFile};

			\addplot [dash dot dot, name path=lower_mse_only, color=mseonlycolor, line width=1.0pt, forget plot, draw=none]
			table[x=p, y=mu4_min, col sep=comma] {\mseonlyFile};

			\addplot [fill=mseonlycolor, draw=none, forget plot, opacity=0.2]
			fill between[
			of=upper_mse_only and lower_mse_only
			];
			
			\addplot [color=orange, line width=1.0pt, only marks, mark=x, mark options={solid, orange}]
			table[x=input, y=l3 , col sep=comma] {data/results/lp/linear_problem_data_optinn_trained_on.csv};

			\nextgroupplot[group/empty plot]

			\nextgroupplot[
			ymin=-0.04,
			ymax=0.6,
			ylabel={$\mu{}^{5}$},
			xlabel style={font=\color{white!15!black}},
			xlabel={$p$},
			]
			\addplot [color=udsred, line width=1.0pt, forget plot]
			table[x=input, y=l4, col sep=comma] {data/results/lp/linear_problem_data_resolution.csv};
			
			\addplot [dashed, color=udsgreen, line width=1.0pt, forget plot]
			table[x=p, y=mu5_mean, col sep=comma] {\optinnFile};
			
			\addplot [dashed, name path=upper_optinn, color=udsgreen, line width=1.0pt, forget plot, draw=none]
			table[x=p, y=mu5_max, col sep=comma] {\optinnFile};
			
			\addplot [dashed, name path=lower_optinn, color=udsgreen, line width=1.0pt, forget plot, draw=none]
			table[x=p, y=mu5_min, col sep=comma] {\optinnFile};
			
			\addplot [fill=udsgreen, draw=none, opacity=0.2]
			fill between[
			of=upper_optinn and lower_optinn
			];

			\addplot [dotted, color=black, line width=1.0pt, forget plot]
			table[x=p, y=mu5_mean, col sep=comma] {\nodataFile};

			\addplot [dotted, name path=upper_no_data, color=black, line width=1.0pt, forget plot, draw=none]
			table[x=p, y=mu5_max, col sep=comma] {\nodataFile};

			\addplot [dotted, name path=lower_no_data, color=black, line width=1.0pt, forget plot, draw=none]
			table[x=p, y=mu5_min, col sep=comma] {\nodataFile};

			\addplot [fill=black, draw=none, forget plot, opacity=0.2]
			fill between[
			of=upper_no_data and lower_no_data
			];

			\addplot [dash dot dot, color=mseonlycolor, line width=1.0pt, forget plot]
			table[x=p, y=mu5_mean, col sep=comma] {\mseonlyFile};

			\addplot [dash dot dot, name path=upper_mse_only, color=mseonlycolor, line width=1.0pt, forget plot, draw=none]
			table[x=p, y=mu5_max, col sep=comma] {\mseonlyFile};

			\addplot [dash dot dot, name path=lower_mse_only, color=mseonlycolor, line width=1.0pt, forget plot, draw=none]
			table[x=p, y=mu5_min, col sep=comma] {\mseonlyFile};

			\addplot [fill=mseonlycolor, draw=none, forget plot, opacity=0.2]
			fill between[
			of=upper_mse_only and lower_mse_only
			];
			
			\addplot [color=orange, line width=1.0pt, only marks, mark=x, mark options={solid, orange}]
			table[x=input, y=l4 , col sep=comma] {data/results/lp/linear_problem_data_optinn_trained_on.csv};
			
		\end{groupplot}
		
		\node[anchor=north] at (current bounding box.south) {\pgfplotslegendfromname{commonlegendLinNN}};
	\end{tikzpicture}%
	\caption{\AOcorrection{Predictions of a QPNN, OptINNs trained with four and with no data points, and an MLP trained with the MSE-loss only. Lines show the mean and shaded areas the minimum and maximum over five training runs. All methods except the MLP track the primal variables well. The OptINNs also give reasonable dual estimates, with a larger spread at the jumps.}}
	\label{fig:linear_demarchi_nn}
\end{figure}

%% file: figs/linear_demarchi_convergence.tex


%
\def\filedata{data/results/lp/optinn/lp_optinn_val_curves_10.csv}
\def\filenodata{data/results/lp_nodata/optinn/lp_nodata_optinn_val_curves_10.csv}
\begin{figure}[!hbt]
	\centering
	\begin{tikzpicture}
		\begin{groupplot}[
			group style={
				group size=2 by 2,
				horizontal sep=1.2cm,
				vertical sep=1cm
			},
			ylabel shift=-4pt,
			scale only axis,
            height=0.18\linewidth,
			width=0.4\linewidth,
            ymode=log,
            xmax=200000,
			enlarge x limits=0.05,
            grid=major,
            unbounded coords=jump,
			]
			\nextgroupplot[
                title={Validation Data Loss},
                ylabel={4 data points},
				ymin=5e-3, ymax=2e2
			]
			\addplot [color=udsgreen, line width=1.0pt, forget plot]
			table[x=epoch, y=seed3_data, col sep=comma] {\filedata};

			\addplot [color=udsred, line width=1.0pt, forget plot]
			table[x=epoch, y=seed4_data, col sep=comma] {\filedata};

			\addplot [color=udsblue, line width=1.0pt, forget plot]
			table[x=epoch, y=seed5_data, col sep=comma] {\filedata};

			\addplot [color=udsyellow, line width=1.0pt, forget plot]
			table[x=epoch, y=seed6_data, col sep=comma] {\filedata};

			\addplot [color=orange, line width=1.0pt, forget plot]
			table[x=epoch, y=seed7_data, col sep=comma] {\filedata};

			\nextgroupplot[
                title={Validation KKT Penalty},
                ymin=3e-4, ymax=2
            ]
            \addplot [color=udsgreen, line width=1.0pt, forget plot]
			table[x=epoch, y=seed3_kkt, col sep=comma] {\filedata};

			\addplot [color=udsred, line width=1.0pt, forget plot]
			table[x=epoch, y=seed4_kkt, col sep=comma] {\filedata};

			\addplot [color=udsblue, line width=1.0pt, forget plot]
			table[x=epoch, y=seed5_kkt, col sep=comma] {\filedata};

			\addplot [color=udsyellow, line width=1.0pt, forget plot]
			table[x=epoch, y=seed6_kkt, col sep=comma] {\filedata};

			\addplot [color=orange, line width=1.0pt, forget plot]
			table[x=epoch, y=seed7_kkt, col sep=comma] {\filedata};

            \nextgroupplot[
                xlabel={Epoch},
                ylabel={no data points},
				ymin=5e-3, ymax=2e2
			]
			\addplot [color=udsgreen, line width=1.0pt, forget plot]
			table[x=epoch, y=seed3_data, col sep=comma] {\filenodata};

			\addplot [color=udsred, line width=1.0pt, forget plot]
			table[x=epoch, y=seed4_data, col sep=comma] {\filenodata};

			\addplot [color=udsblue, line width=1.0pt, forget plot]
			table[x=epoch, y=seed5_data, col sep=comma] {\filenodata};

			\addplot [color=udsyellow, line width=1.0pt, forget plot]
			table[x=epoch, y=seed6_data, col sep=comma] {\filenodata};

			\addplot [color=orange, line width=1.0pt, forget plot]
			table[x=epoch, y=seed7_data, col sep=comma] {\filenodata};

			\nextgroupplot[
                xlabel={Epoch},
                ymin=3e-4, ymax=2
            ]
            \addplot [color=udsgreen, line width=1.0pt, forget plot]
			table[x=epoch, y=seed3_kkt, col sep=comma] {\filenodata};

			\addplot [color=udsred, line width=1.0pt, forget plot]
			table[x=epoch, y=seed4_kkt, col sep=comma] {\filenodata};

			\addplot [color=udsblue, line width=1.0pt, forget plot]
			table[x=epoch, y=seed5_kkt, col sep=comma] {\filenodata};

			\addplot [color=udsyellow, line width=1.0pt, forget plot]
			table[x=epoch, y=seed6_kkt, col sep=comma] {\filenodata};

			\addplot [color=orange, line width=1.0pt, forget plot]
			table[x=epoch, y=seed7_kkt, col sep=comma] {\filenodata};
		\end{groupplot}
	\end{tikzpicture}%
	\caption{\AOcorrection{Training curves of OptINNs trained with 4 data points (top) and without data (bottom): data-based validation loss on 256 points, used only for evaluation (left), and KKT-loss (right). Colors denote random seeds. \MKHcorrection{One run with data and three without} stop early after 20000 epochs without improvement. No overfitting is observed.}
	}
	\label{fig:linear_demarchi_convergence}
\end{figure}

%% file: figs/quadratic_nonl_solution_diff_kkt.tex

\def\quadrSeed{3}
\def\quadrOptinnFile{data/results/nonconvex/optinn/nonconvex_optinn_projections.csv}
\def\quadrPenaltyFile{data/results/nonconvex/qpnn/nonconvex_qpnn_projections.csv}
\begin{figure}[htb]
	\centering
	\begin{tikzpicture}[
			declare function={
				f(\x) = 2*(\x)^3 + 3*(\x)^2 + 2*(\x) + 1;
			}
		]
		\definecolor{crimson2143940}{RGB}{214,39,40}
		\definecolor{darkgray176}{RGB}{176,176,176}
		\definecolor{darkorange25512714}{RGB}{255,127,14}
		\definecolor{darkturquoise23190207}{RGB}{23,190,207}
		\definecolor{forestgreen4416044}{RGB}{44,160,44}
		\definecolor{goldenrod18818934}{RGB}{188,189,34}
		\definecolor{gray127}{RGB}{127,127,127}
		\definecolor{mediumpurple148103189}{RGB}{148,103,189}
		\definecolor{orchid227119194}{RGB}{227,119,194}
		\definecolor{sienna1408675}{RGB}{140,86,75}
		\definecolor{steelblue31119180}{RGB}{31,119,180}
		
		\begin{groupplot}[
			group style={
				group size=2 by 1, 
				horizontal sep=3cm, 
			},
			tick align=outside,
			tick pos=left,
			x grid style={darkgray176},
			y grid style={darkgray176},
			xmin=-1.1, xmax=1.1,
			ymin=-1.1, ymax=1.1,
			xtick style={color=black},
			ytick style={color=black},
			axis equal image, 
			width=0.45\textwidth, 
			legend style={
				/tikz/every even column/.append style={column sep=0.3cm},
				legend columns=-1, 
				at={(0.00,0.97)}, 
				anchor=north west, 
				legend cell align=left, 
				align=left, 
				draw=white!15!black
			},
			xlabel={$\makev{z}^1, \pv^1$},
			ylabel={$\makev{z}^2, \pv^2$},
			]
			
			\nextgroupplot[
				title={QPNN Results},
				legend to name=commonlegendQuad,
			]
			\addplot[thick, udsred, domain=-1:0, name path=A]
			{%
				f(x)
			};
			\addlegendentry{Constraints}
			\addplot[thick, udsred, domain=-1:0, forget plot, name path=B]
			{%
				-f(x)
			};
			\addplot[thick, udsred, domain=0:1, forget plot, name path=C]
			{%
				f(-x)
			};
			\addplot[thick, udsred, domain=0:1, forget plot, name path=D]
			{%
				-f(-x)
			};
			\addplot [udsred!20!white, forget plot] fill between [of=A and B];
			\addplot [udsred!20!white, forget plot] fill between [of=C and D];
			
			\addplot [thick, goldenrod18818934, mark=diamond*, mark size=3, mark options={solid}, only marks]
			table[x=p0, y=p1, col sep=comma] {data/results/nonconvex/quadr_nonl_data_optinn_trained_on.csv};
			\addlegendentry{Training data $\pv_\mathrm{train}\in\mathcal{\bar{P}}$};
			
			\addplot [thick, steelblue31119180, mark=*, mark size=1.5, mark options={solid}, only marks]
			table[x=p1, y=p2, col sep=comma] {\quadrPenaltyFile};
			\addlegendentry{Test data $\pv_\mathrm{test}\in\mathcal{\bar{P}}$};

			\addplot [thick, darkorange25512714, mark=x, mark size=3, mark options={solid}, only marks]
			table[x=z1_seed\quadrSeed, y=z2_seed\quadrSeed, col sep=comma] {\quadrPenaltyFile};
			\addlegendentry{Prediction $\nnt[x]{\pv_\mathrm{test}}$};

			\pgfplotstableread[col sep=comma]{\quadrPenaltyFile}{\datatableB}
			\pgfplotstablegetrowsof{\datatableB}
			\pgfmathtruncatemacro{\lastrow}{\pgfplotsretval - 1}
			\foreach \i in {0,...,\lastrow}{
				\pgfplotstablegetelem{\i}{p1}\of\datatableB
				\let\pzero\pgfplotsretval
				\pgfplotstablegetelem{\i}{p2}\of\datatableB
				\let\pone\pgfplotsretval
				\pgfplotstablegetelem{\i}{z1_seed\quadrSeed}\of\datatableB
			\let\xzero\pgfplotsretval
			\pgfplotstablegetelem{\i}{z2_seed\quadrSeed}\of\datatableB
			\let\xone\pgfplotsretval
			\addplot [black] coordinates {(\pzero,\pone) (\xzero,\xone)};
		}
		
		\nextgroupplot[
		title={OptINN Results}
		]
		\addplot[thick, udsred, domain=-1:0, forget plot, name path=A]
			{%
				f(x)
			};
			\addplot[thick, udsred, domain=-1:0, forget plot, name path=B]
			{%
				-f(x)
			};
			\addplot[thick, udsred, domain=0:1, forget plot, name path=C]
			{%
				f(-x)
			};
			\addplot[thick, udsred, domain=0:1, forget plot, name path=D]
			{%
				-f(-x)
			};
		\addplot [udsred!20!white, forget plot] fill between [of=A and B];
		\addplot [udsred!20!white, forget plot] fill between [of=C and D];
		
		\addplot [thick, goldenrod18818934, mark=diamond*, mark size=3, mark options={solid}, only marks]
		table[x=p0, y=p1, col sep=comma] {data/results/nonconvex/quadr_nonl_data_optinn_trained_on.csv};
		
		\addplot [thick, steelblue31119180, mark=*, mark size=1.5, mark options={solid}, only marks]
		table[x=p1, y=p2, col sep=comma] {\quadrOptinnFile};

		\addplot [thick, darkorange25512714, mark=x, mark size=3, mark options={solid}, only marks]
		table[x=z1_seed\quadrSeed, y=z2_seed\quadrSeed, col sep=comma] {\quadrOptinnFile};

		\pgfplotstableread[col sep=comma]{\quadrOptinnFile}{\datatableB}
		\pgfplotstablegetrowsof{\datatableB}
		\pgfmathtruncatemacro{\lastrow}{\pgfplotsretval - 1}
		\foreach \i in {0,...,\lastrow}{
			\pgfplotstablegetelem{\i}{p1}\of\datatableB
			\let\pzero\pgfplotsretval
			\pgfplotstablegetelem{\i}{p2}\of\datatableB
			\let\pone\pgfplotsretval
			\pgfplotstablegetelem{\i}{z1_seed\quadrSeed}\of\datatableB
			\let\xzero\pgfplotsretval
			\pgfplotstablegetelem{\i}{z2_seed\quadrSeed}\of\datatableB
			\let\xone\pgfplotsretval
			\addplot [black] coordinates {(\pzero,\pone) (\xzero,\xone)};
		}
		
	\end{groupplot}
	
	\node[anchor=north] at (current bounding box.south) {\pgfplotslegendfromname{commonlegendQuad}};
\end{tikzpicture}

\caption{Comparison of OptINN and QPNN results for the quadratic problem with nonconvex constraints. \MKHcorrection{Both networks} learned that for $\pv$ outside of the feasible set large changes are necessary and for $\pv$ inside, the results are much closer to an identity mapping. Moreover, the vector between $\pv$ and $\makev{z}$ is perpendicular to the corresponding constraint's tangent as predicted to be optimal.}\label{fig:quadratic_nonl_problem_results}
\end{figure}

%% file: figs/rocket_car_heatmaps_sparse_data.tex
\def\dirOptinn{data/results/rocket_car/optinn}
\def\dirQpnn{data/results/rocket_car/qpnn}
\input{\dirOptinn/rocket_car_optinn_x_error_colormap.tex}
\input{\dirOptinn/rocket_car_optinn_v_error_colormap.tex}
\input{\dirOptinn/rocket_car_optinn_u_error_colormap.tex}
\input{\dirQpnn/rocket_car_qpnn_x_error_colormap.tex}
\input{\dirQpnn/rocket_car_qpnn_v_error_colormap.tex}
\input{\dirQpnn/rocket_car_qpnn_u_error_colormap.tex}
\tikzset{
	rocketcar training data/.style={thick, dashed, mark=none, color=black},
	rocketcar active bound/.style={thick, dotted, mark=none, color=black},
}
\begin{figure}[ht]
	\centering
	\begin{tikzpicture}
		\begin{groupplot}[
			group style={
				group size=3 by 2,
				horizontal sep=3.0cm,
				vertical sep=1.0cm
			},
			width=0.24\linewidth,
			height=0.28\linewidth,
			colorbar,
			xmin=-0.5,
			xmax=32.5,
			ymin=-0.2,
			ymax=40.2,
			xlabel={Time step $k$},
			axis on top,
			legend style={
				/tikz/every even column/.append style={column sep=0.2cm},
				legend columns=-1,
			},
			]
			\nextgroupplot[title={$z$},
			ylabel style={align=center},
			ylabel={OptINN\\ Terminal position $p$},
			rocket_car-optinn-x,
			xlabel={},
			legend to name=commonlegendRocketcarHeatmap]
			\addplot[forget plot] graphics [xmin=-0.5, xmax=32.5, ymin=-0.2, ymax=40.2] {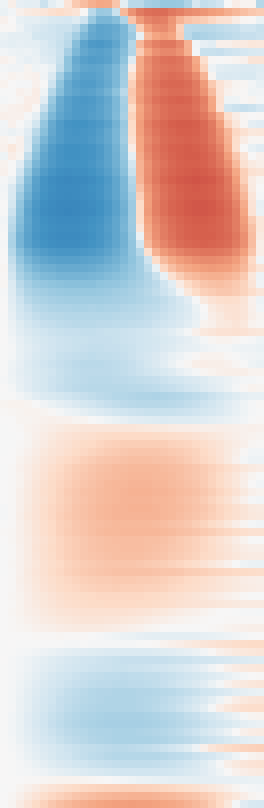};
			\addplot[rocketcar training data] coordinates {(-0.5, 0) (32.5, 0)};
			\addlegendentry{Training data}
			\addplot[rocketcar training data, forget plot] coordinates {(-0.5, 20) (32.5, 20)};
			\addplot[rocketcar training data, forget plot] coordinates {(-0.5, 40) (32.5, 40)};
			\addplot[rocketcar active bound] coordinates {(-0.5, 28.16) (32.5, 28.16)};
			\addlegendentry{$p=28.16$}

			\nextgroupplot[title={$v$},
			rocket_car-optinn-v,
			xlabel={}]
			\addplot[forget plot] graphics [xmin=-0.5, xmax=32.5, ymin=-0.2, ymax=40.2] {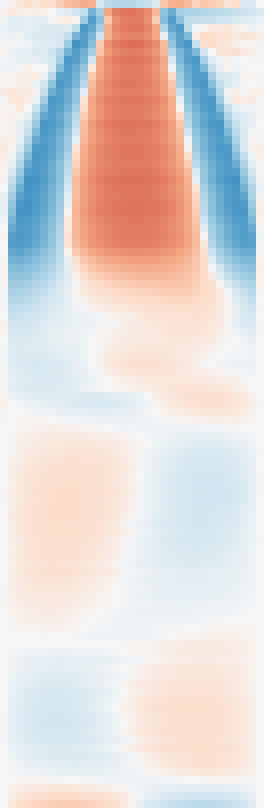};
			\pgfplotsinvokeforeach{0, 20, 40}{
				\addplot[rocketcar training data] coordinates {(-0.5, #1) (32.5, #1)};
			}
			\addplot[rocketcar active bound] coordinates {(-0.5, 28.16) (32.5, 28.16)};

			\nextgroupplot[title={$u$},
			xmax=31.5,
			rocket_car-optinn-u,
			xlabel={}]
			\addplot[forget plot] graphics [xmin=-0.5, xmax=31.5, ymin=-0.2, ymax=40.2] {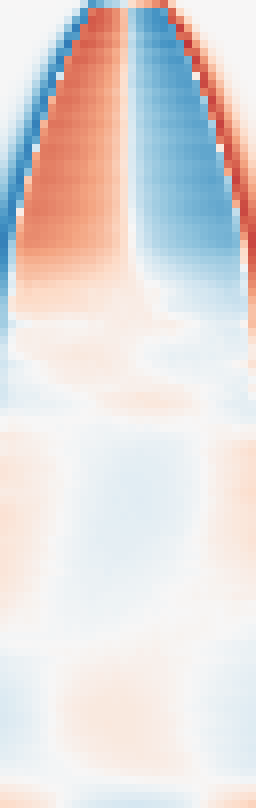};
			\pgfplotsinvokeforeach{0, 20, 40}{
				\addplot[rocketcar training data] coordinates {(-0.5, #1) (31.5, #1)};
			}
			\addplot[rocketcar active bound] coordinates {(-0.5, 28.16) (31.5, 28.16)};

			\nextgroupplot[
			ylabel style={align=center},
			ylabel={QPNN\\ Terminal position $p$},
			rocket_car-qpnn-x]
			\addplot[forget plot] graphics [xmin=-0.5, xmax=32.5, ymin=-0.2, ymax=40.2] {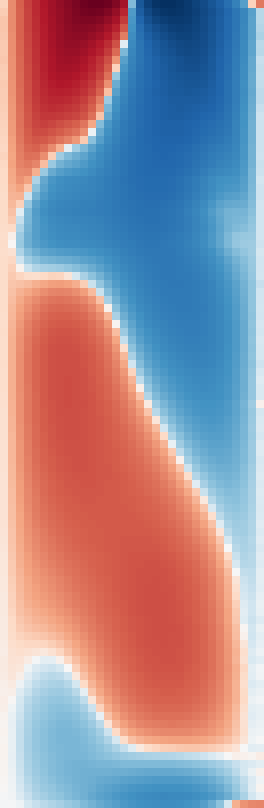};
			\pgfplotsinvokeforeach{0, 20, 40}{
				\addplot[rocketcar training data] coordinates {(-0.5, #1) (32.5, #1)};
			}
			\addplot[rocketcar active bound] coordinates {(-0.5, 28.16) (32.5, 28.16)};

			\nextgroupplot[
			rocket_car-qpnn-v]
			\addplot[forget plot] graphics [xmin=-0.5, xmax=32.5, ymin=-0.2, ymax=40.2] {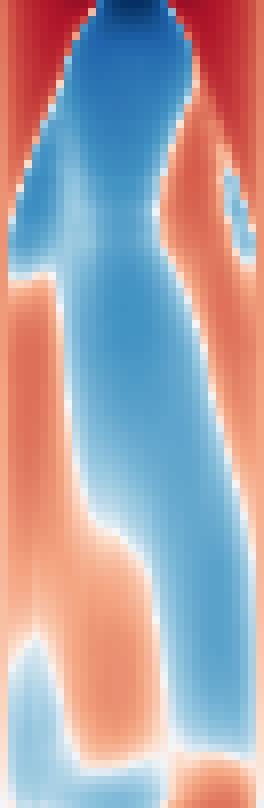};
			\pgfplotsinvokeforeach{0, 20, 40}{
				\addplot[rocketcar training data] coordinates {(-0.5, #1) (32.5, #1)};
			}
			\addplot[rocketcar active bound] coordinates {(-0.5, 28.16) (32.5, 28.16)};

			\nextgroupplot[
			xmax=31.5,
			rocket_car-qpnn-u]
			\addplot[forget plot] graphics [xmin=-0.5, xmax=31.5, ymin=-0.2, ymax=40.2] {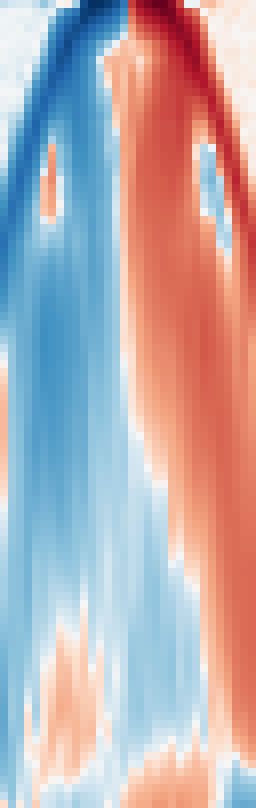};
			\pgfplotsinvokeforeach{0, 20, 40}{
				\addplot[rocketcar training data] coordinates {(-0.5, #1) (31.5, #1)};
			}
			\addplot[rocketcar active bound] coordinates {(-0.5, 28.16) (31.5, 28.16)};
		\end{groupplot}
		\node[anchor=north] at (current bounding box.south) {\pgfplotslegendfromname{commonlegendRocketcarHeatmap}};
	\end{tikzpicture}
	\caption{\AOcorrection{Signed deviation of the OptINN and QPNN predictions from the reference solution for the double integrator at 101 equidistant values of $p${, shown for the training run with the median primal error}. Dashed lines mark training points. {The color scale is logarithmic in magnitude and linear near zero.} The largest deviations occur for $p \ge 28.16$, where the input bounds become active.}}
	\label{fig:heatmap_rocketcar_diff}
\end{figure}

%% file: figs/rocket_car_1d_results.tex
\def\dirOptinn{data/results/rocket_car/optinn}
\def\dirQpnn{data/results/rocket_car/qpnn}
\def\pWorst{39.6}
\begin{figure}[htb]
	\centering
	\definecolor{darkgray176}{RGB}{176,176,176}
	\begin{tikzpicture}
		\begin{groupplot}[
			group style={
				group size=3 by 1,
				vertical sep=2cm
			},
			yticklabel shift={1mm},
			width=0.35\linewidth,
			xlabel={Time step $k$},
			tick align=outside,
			tick pos=left,
			x grid style={darkgray176},
			xtick style={color=black},
			y grid style={darkgray176},
			ytick style={color=black},
			yticklabel style={anchor=center},
			xmajorgrids=true,
			ymajorgrids=true,
			legend style={
				/tikz/every even column/.append style={column sep=1cm},
				legend columns=-1, 
				at={(0.03,0.97)},
				anchor=north west,
				legend cell align=left,
				align=left,
				draw=white!15!black
			}, xmin=0, xmax=32,
			enlarge x limits=0.05,
			table/x=k, table/col sep=comma,
			x filter/.expression={abs(\thisrow{p}-\pWorst) < 0.01 ? x : nan},
			unbounded coords=discard,
			]
			\nextgroupplot[
			ylabel={$z$},
			legend to name=commonlegendRocketcar,
			yticklabel shift=6pt
			]
			\addplot[mark=none, color=udsred, line width=1pt] table[y=truth] {\dirOptinn/rocket_car_optinn_x.csv};
			\addlegendentry{Reference solution}

			\addplot[dashed, mark=none, color=udsgreen, line width=1pt] table[y=seed3] {\dirOptinn/rocket_car_optinn_x.csv};
			\addlegendentry{OptINN prediction}

			\addplot[dotted, mark=none, color=udsorange, line width=1pt] table[y=seed7] {\dirQpnn/rocket_car_qpnn_x.csv};
			\addlegendentry{QPNN prediction}

			\nextgroupplot[
			ylabel={$v$},
			]
			\addplot[mark=none, color=udsred, line width=1pt] table[y=truth] {\dirOptinn/rocket_car_optinn_v.csv};

			\addplot[dashed, mark=none, color=udsgreen, line width=1pt] table[y=seed3] {\dirOptinn/rocket_car_optinn_v.csv};

			\addplot[dotted, mark=none, color=udsorange, line width=1pt] table[y=seed7] {\dirQpnn/rocket_car_qpnn_v.csv};

			\nextgroupplot[
			ylabel={$u$},
			]

			\addplot[mark=none, color=udsred, line width=1pt] table[y=truth] {\dirOptinn/rocket_car_optinn_u.csv};

			\addplot[dashed, mark=none, color=udsgreen, line width=1pt] table[y=seed3] {\dirOptinn/rocket_car_optinn_u.csv};

			\addplot[dotted, mark=none, color=udsorange, line width=1pt] table[y=seed7] {\dirQpnn/rocket_car_qpnn_u.csv};
		\end{groupplot}
		\node[anchor=north] at (current bounding box.south) {\pgfplotslegendfromname{commonlegendRocketcar}};
	\end{tikzpicture}
	\caption{\AOcorrection{Reference solution and predictions for the double integrator at $p = 39.6$, the test point with the largest OptINN deviation. The OptINN deviates less than the QPNN, especially in $u$.}}
	\label{fig:rocket_car_1d}
\end{figure}

%% file: figs/rocketcar_diff_primal_variables.tex


\def\fileOptinn{data/results/rocket_car/optinn/rocket_car_optinn_deviation_summary.csv}
\def\filePenalty{data/results/rocket_car/qpnn/rocket_car_qpnn_deviation_summary.csv}
\begin{figure}[htb]
	\centering
	\begin{tikzpicture}
		\begin{groupplot}[
			group style={
				group size=3 by 1,
				horizontal sep=1.2cm,
			},
			ylabel shift=-4pt,
			xlabel={Terminal position $p$},
			scale only axis,
			xmin=0,
			xmax=40,
			axis background/.style={fill=white},
			height=0.18\linewidth,
			width=0.25\linewidth,
            enlarge x limits=0.05,
            grid=major
			]
			\nextgroupplot[
			title={$z$},
			ylabel={Signed deviation}
			]
			\addplot [color=udsgreen, line width=1.0pt, forget plot]
			table[x=p, y=mean_z, col sep=comma] {\fileOptinn};

			\addplot [color=udsgreen, line width=1.0pt, forget plot, draw=none, name path=lower_optinn]
			table[x=p, y=min_z, col sep=comma] {\fileOptinn};

			\addplot [color=udsgreen, line width=1.0pt, forget plot, draw=none, name path=upper_optinn]
			table[x=p, y=max_z, col sep=comma] {\fileOptinn};

            \addplot [fill=udsgreen, draw=none, opacity=0.2, forget plot, postaction={pattern=north east lines, pattern color=udsgreen}]
				fill between[
					of=upper_optinn and lower_optinn
			];

			\addplot [dashed, color=udsorange, line width=1.0pt, forget plot]
			table[x=p, y=mean_z, col sep=comma] {\filePenalty};

			\addplot [color=udsorange, line width=1.0pt, forget plot, draw=none, name path=lower_penalty]
			table[x=p, y=min_z, col sep=comma] {\filePenalty};

			\addplot [color=udsorange, line width=1.0pt, forget plot, draw=none, name path=upper_penalty]
			table[x=p, y=max_z, col sep=comma] {\filePenalty};

            \addplot [fill=udsorange, draw=none, opacity=0.2, forget plot, postaction={pattern=north west lines, pattern color=udsorange}]
				fill between[
					of=upper_penalty and lower_penalty
			];
			
			\nextgroupplot[
			title={$v$},
			legend to name=commonRCDiff, 
			legend style={
				/tikz/every even column/.append style={column sep=5mm},
				legend columns=-1, 
				at={(0.03,0.97)}, 
				anchor=north west, 
				legend cell align=left, 
				align=left, 
				draw=white!15!black
			}
			]
            \addplot [color=udsgreen, line width=1.0pt]
			table[x=p, y=mean_v, col sep=comma] {\fileOptinn};
			\addlegendentry{OptINN - mean}

			\addplot [color=udsgreen, line width=1.0pt, forget plot, draw=none, name path=lower_optinn]
			table[x=p, y=min_v, col sep=comma] {\fileOptinn};

			\addplot [color=udsgreen, line width=1.0pt, forget plot, draw=none, name path=upper_optinn]
			table[x=p, y=max_v, col sep=comma] {\fileOptinn};

			\addplot [fill=udsgreen, draw=none, opacity=0.2, postaction={pattern=north east lines, pattern color=udsgreen}]
										fill between[
											of=upper_optinn and lower_optinn
									];
			\addlegendentry{OptINN - min/max}

			\addplot [dashed, color=udsorange, line width=1.0pt]
			table[x=p, y=mean_v, col sep=comma] {\filePenalty};
			\addlegendentry{QPNN - mean}

			\addplot [color=udsorange, line width=1.0pt, forget plot, draw=none, name path=lower_penalty]
			table[x=p, y=min_v, col sep=comma] {\filePenalty};

			\addplot [color=udsorange, line width=1.0pt, forget plot, draw=none, name path=upper_penalty]
			table[x=p, y=max_v, col sep=comma] {\filePenalty};

            \addplot [fill=udsorange, draw=none, opacity=0.2, postaction={pattern=north west lines, pattern color=udsorange}]
				fill between[
					of=upper_penalty and lower_penalty
			];
			\addlegendentry{QPNN - min/max}
			
			\nextgroupplot[
			title={$u$},
			]
            \addplot [color=udsgreen, line width=1.0pt, forget plot]
			table[x=p, y=mean_u, col sep=comma] {\fileOptinn};

			\addplot [color=udsgreen, line width=1.0pt, forget plot, draw=none, name path=lower_optinn]
			table[x=p, y=min_u, col sep=comma] {\fileOptinn};

			\addplot [color=udsgreen, line width=1.0pt, forget plot, draw=none, name path=upper_optinn]
			table[x=p, y=max_u, col sep=comma] {\fileOptinn};

            \addplot [fill=udsgreen, draw=none, opacity=0.2, forget plot, postaction={pattern=north east lines, pattern color=udsgreen}]
				fill between[
					of=upper_optinn and lower_optinn
			];

			\addplot [dashed, color=udsorange, line width=1.0pt, forget plot]
			table[x=p, y=mean_u, col sep=comma] {\filePenalty};

			\addplot [color=udsorange, line width=1.0pt, forget plot, draw=none, name path=lower_penalty]
			table[x=p, y=min_u, col sep=comma] {\filePenalty};

			\addplot [color=udsorange, line width=1.0pt, forget plot, draw=none, name path=upper_penalty]
			table[x=p, y=max_u, col sep=comma] {\filePenalty};

            \addplot [fill=udsorange, draw=none, opacity=0.2, postaction={pattern=north west lines, pattern color=udsorange}]
				fill between[
					of=upper_penalty and lower_penalty
			];			
		\end{groupplot}
		\node[anchor=north] at (current bounding box.south) {\pgfplotslegendfromname{commonRCDiff}};
	\end{tikzpicture}%
	\caption{Mean, minimum, and maximum signed deviation from the reference solution for both methods over the parameter value.
	The OptINN produces smaller deviations for most parameters and a significantly lower error over all parameters}
	\label{fig:rocketcar_difference}
\end{figure}

%% file: figs/pendulum_cost_comparison.tex
\def\optinnfile{data/results/pendulum/optinn/pendulum_optinn_curves.csv}
\def\penaltyfile{data/results/pendulum/qpnn/pendulum_qpnn_curves.csv}

\begin{figure}[htb]
	\centering
	\definecolor{darkgray176}{RGB}{176,176,176}
	\begin{tikzpicture}
		\begin{groupplot}[
			group style={
				group size=3 by 1,
				horizontal sep=1.3cm
			},
			width=0.34\linewidth,
			xlabel={Transition time $p$},
			xtick={6,9,12,15},
			tick align=outside,
			tick pos=left,
			x grid style={darkgray176},
			xtick style={color=black},
			y grid style={darkgray176},
			ytick style={color=black},
			xmajorgrids=true,
			ymajorgrids=true,
			title style={font=\small},
			extra x ticks={6,6.6,7.2,7.8,8.4,9,9.6,10.2,10.8,11.4,12,12.6,13.2,13.8,14.4,15},
			extra x tick labels={},
			extra x tick style={grid=major, major grid style={thick, dotted, black}, major tick length=0pt},
			legend style={
				/tikz/every even column/.append style={column sep=1cm},
				legend columns=3,
				legend cell align=left,
				align=left,
				draw=white!15!black
			}
			]
			\nextgroupplot[
			title={(a) Equality constraint violation},
			ymode=log
			]
			\addplot[dashed, mark=none, color=udsgreen, line width=1pt] table[x=p, y=eq_mean, col sep=comma] {\optinnfile};

			\addplot[draw=none, name path=optinn_min] table[x=p, y=eq_min, col sep=comma, forget plot] {\optinnfile};
			\addplot[draw=none, name path=optinn_max] table[x=p, y=eq_max, col sep=comma, forget plot] {\optinnfile};
			\addplot[fill=udsgreen, fill opacity=0.2] fill between[of=optinn_min and optinn_max];

			\addplot[dotted, mark=none, color=udsorange, line width=1pt] table[x=p, y=eq_mean, col sep=comma] {\penaltyfile};

			\addplot[draw=none, name path=pmnn_min] table[x=p, y=eq_min, col sep=comma, forget plot] {\penaltyfile};
			\addplot[draw=none, name path=pmnn_max] table[x=p, y=eq_max, col sep=comma, forget plot] {\penaltyfile};
			\addplot[fill=udsorange, fill opacity=0.2] fill between[of=pmnn_min and pmnn_max];

			\nextgroupplot[
			title={(b) Cost of prediction},
			legend to name=commonlegendPendulum,
			]
			\addplot[mark=none, color=udsred, line width=1pt] table[x=p, y=cost_truth, col sep=comma] {\optinnfile};
			\addlegendentry{Num. Optim.}

			\addplot[dashed, mark=none, color=udsgreen, line width=1pt] table[x=p, y=cost_mean, col sep=comma] {\optinnfile};
			\addlegendentry{OptINN mean}

			\addplot[draw=none, name path=optinn_min] table[x=p, y=cost_min, col sep=comma, forget plot] {\optinnfile};
			\addplot[draw=none, name path=optinn_max] table[x=p, y=cost_max, col sep=comma, forget plot] {\optinnfile};
			\addplot[fill=udsgreen, fill opacity=0.2] fill between[of=optinn_min and optinn_max];
			\addlegendentry{OptINN min/max}

			\addlegendimage{thick, dotted, black}
			\addlegendentry{Training data}

			\addplot[dotted, mark=none, color=udsorange, line width=1pt] table[x=p, y=cost_mean, col sep=comma] {\penaltyfile};
			\addlegendentry{QPNN mean}

			\addplot[draw=none, name path=pmnn_min] table[x=p, y=cost_min, col sep=comma, forget plot] {\penaltyfile};
			\addplot[draw=none, name path=pmnn_max] table[x=p, y=cost_max, col sep=comma, forget plot] {\penaltyfile};
			\addplot[fill=udsorange, fill opacity=0.2] fill between[of=pmnn_min and pmnn_max];
			\addlegendentry{QPNN min/max}

			\nextgroupplot[
			title={(c) Relative cost gap $\delta$},
			]
			\addplot[mark=none, color=udsred, line width=1pt, domain=6:15, samples=2] {0};

			\addplot[dashed, mark=none, color=udsgreen, line width=1pt] table[x=p, y=relcost_mean, col sep=comma] {\optinnfile};

			\addplot[draw=none, name path=optinn_min] table[x=p, y=relcost_min, col sep=comma, forget plot] {\optinnfile};
			\addplot[draw=none, name path=optinn_max] table[x=p, y=relcost_max, col sep=comma, forget plot] {\optinnfile};
			\addplot[fill=udsgreen, fill opacity=0.2] fill between[of=optinn_min and optinn_max];

			\addplot[dotted, mark=none, color=udsorange, line width=1pt] table[x=p, y=relcost_mean, col sep=comma] {\penaltyfile};

			\addplot[draw=none, name path=pmnn_min] table[x=p, y=relcost_min, col sep=comma, forget plot] {\penaltyfile};
			\addplot[draw=none, name path=pmnn_max] table[x=p, y=relcost_max, col sep=comma, forget plot] {\penaltyfile};
			\addplot[fill=udsorange, fill opacity=0.2] fill between[of=pmnn_min and pmnn_max];
		\end{groupplot}
		\node[anchor=north] at (current bounding box.south) {\pgfplotslegendfromname{commonlegendPendulum}};
	\end{tikzpicture}
	\caption{
\AOcorrection{Results for the pendulum over the transition time $p$: (a) mean absolute violation of the equality constraints defined  in \eqref{eq:def_L_EqViol}, (b) cost of the prediction compared with the cost of the reference solution, and (c) relative cost gap $\delta$ \eqref{eq:def_delta}. Lines show the mean and shaded areas the minimum and maximum over the five training runs. Dotted vertical lines mark the training points}
}
	\label{fig:pendulum_example}
\end{figure}

%% file: tables/numeric_metrics.tex
\begin{table}[ht]
    \centering
    \caption{Summary of numerical evaluation metrics. The value shown is the mean over all parameter values and all five validation trainings.
    The relative cost gap $\delta$ is the difference between the cost of the prediction and the optimal cost, divided by the range of the optimal cost over all parameter values.
    The column $\delta^-$ gives the fraction of parameter values with $\delta < -10^{-4}$, for which the cost of the prediction is below the optimal cost}
    \label{tab:numerical_summary}
    \renewcommand{\arraystretch}{1.2} 
    \small
    \setlength{\tabcolsep}{4pt}
    \begin{tabularx}{\textwidth}{@{} l l *{4}{>{\centering\arraybackslash}X} >{\centering\arraybackslash\hsize=0.8\hsize}X >{\centering\arraybackslash\hsize=1.2\hsize}X @{}}
        \toprule
        \textbf{Problem} & \textbf{Method} & \textbf{Primal MSE } & \textbf{Eq. Violation } & \textbf{Ineq. Violation } & \textbf{Mean $|\delta|$ } & \textbf{$\delta^-$ } \\
        \midrule

        \multirow{4}{*}{Lin. Prog.}
             & QPNN & $\boldsymbol{7.17 \times 10^{-3}}$ & -- & $1.13 \times 10^{-4}$ & $2.27 \times 10^{-4}$ & $70.2\,\%$ \\
             & OptINN & $9.55 \times 10^{-3}$ & -- & $\boldsymbol{4.04 \times 10^{-5}}$ & $\boldsymbol{1.78 \times 10^{-4}}$ & $\boldsymbol{8.2\,\%}$ \\
             & OptINN (no data) & $2.72 \times 10^{-2}$ & -- & $6.98 \times 10^{-5}$ & $2.92 \times 10^{-4}$ & $14.6\,\%$ \\
             & MLP (MSE only) & $6.18 \times 10^{1}$ & -- & $2.50 \times 10^{-2}$ & $6.52 \times 10^{-2}$ & $49.1\,\%$ \\
        \addlinespace
        \multirow{2}{*}{Nonconvex}
             & QPNN & $2.19 \times 10^{-5}$ & -- & $1.68 \times 10^{-4}$ & $9.41 \times 10^{-4}$ & $42.2\,\%$ \\
             & OptINN & $\boldsymbol{1.31 \times 10^{-5}}$ & -- & $\boldsymbol{3.73 \times 10^{-5}}$ & $\boldsymbol{3.07 \times 10^{-4}}$ & $\boldsymbol{15.3\,\%}$ \\
        \addlinespace
        \multirow{2}{*}{Rocket car}
             & QPNN & $7.60 \times 10^{-2}$ & $1.42 \times 10^{-3}$ & $7.36 \times 10^{-8}$ & $1.27 \times 10^{-2}$ & $93.6\,\%$ \\
             & OptINN & $\boldsymbol{2.30 \times 10^{-3}}$ & $\boldsymbol{7.91 \times 10^{-5}}$ & \textbf{0} & $\boldsymbol{2.07 \times 10^{-4}}$ & $\boldsymbol{13.8\,\%}$ \\

        \toprule
        &  & \textbf{Cost MSE } &  &  &  &  \\
        \midrule

        \multirow{2}{*}{Pendulum}
             & QPNN & $2.55 \times 10^{-1}$ & $1.77 \times 10^{-4}$ & $1.71 \times 10^{-3}$ & $2.36 \times 10^{-2}$ & $36.4\,\%$ \\
             & OptINN & $\boldsymbol{6.42 \times 10^{-3}}$ & $\boldsymbol{7.67 \times 10^{-5}}$ & \textbf{0} & $\boldsymbol{4.42 \times 10^{-3}}$ & $\boldsymbol{16.25\,\%}$ \\

        \bottomrule
    \end{tabularx}
\end{table}

%% file: inc/summary.tex
\section{Conclusion and outlook}\label{sec:sum}
In this work, we proposed Optimality-Informed Neural Networks (OptINNs) to learn the parameter-to-optimal-solution mapping of parametric optimization problems.
Unlike previous works that rely primarily on the quadratic-penalty method, we introduced a problem-specific architecture and a loss function rooted in the Karush-Kuhn-Tucker (KKT) conditions.
This formulation yields two key advantages: first, it enables the simultaneous prediction of optimal decision variables and the corresponding Lagrange multipliers.
Second, unlike data-based validation metrics, the KKT-residual can be computed from randomly sampled parameter values, serving as a self-supervised metric.
Since the KKT-loss is lower-bounded by zero, it is directly applicable for both hyperparameter tuning and model validation without requiring additional ground-truth data.
\MKHcorrection{Moreover, after training, the KKT-residuals can be evaluated at any parameter value to check whether the prediction satisfies the KKT conditions.}
By adopting the perspective of theory-informed \MKHcorrection{machine learning}, we leveraged established techniques such as adaptive weight balancing to train these networks to high accuracy even in low-data regimes.

Our numerical experiments across four benchmark problems, linear programming, nonconvex constraints, and both linear and nonlinear optimal control, demonstrate the efficacy of the proposed method.
For smaller problems, OptINNs achieve a performance comparable to that of penalty-based networks, with the added benefit of predicting dual variables.
However, for the larger optimal control tasks, OptINNs demonstrate superior performance, achieving lower constraint violations and higher fidelity to numerical solvers.
We attribute this superiority to the explicit KKT guidance, which avoids the known tendency of quadratic-penalty methods to converge towards infeasible local minimizers.

The methodological framework presented here is model-agnostic and can be generalized to other parametric models.
Currently, the reliance on Multi-Layer Perceptrons restricts OptINNs to optimization problems with fixed dimensionality.
Future research should therefore investigate dynamic architectures, such as Recurrent Neural Networks or Transformers, to mitigate this restriction and accommodate variable-sized problem instances.

Furthermore, several directions remain for improving training efficiency and applicability.
Using penalty-based methods for pretraining OptINNs could significantly accelerate early convergence.
\MKHcorrection{Surrogate gradients~\cite{neftci2019surrogategradient}, which replace the derivative of an activation function in the backward pass, may improve the gradients of the output layer at changes of the set of active inequality constraints.}
Regarding applications, the rapid inference capabilities of OptINNs make them highly suitable for real-time optimal control; their integration into model predictive control workflows warrants further study.
Additionally, exploiting parametric sensitivity analysis offers a promising avenue for data efficiency.
Computing the derivatives of the solution branches with respect to problem parameters would enable Sobolev training, extracting richer gradient information from the available training data.
Finally, to strengthen the theoretical foundations, future work should investigate the alignment of this problem setting with the field of tame optimization to establish convergence properties.

%% file: inc/appendix.tex
\section{Hyperparameters for numerical examples}
\label{app:hparams}
\noindent
Here, we summarize the hyperparameters used for training Optimality-informed (OptINN) and quadratic-penalty-method-based neural networks (QPNN) in the numerical examples. The depth describes the number of hidden layers, excluding input, output and trivialization layers. The width is the number of neurons in each hidden layer. The learning rate (LR) is the initial learning rate used in the AdamW optimizer~\cite{loshchilov2019adamw} with default parameters, using a learning rate scheduler (ReduceLROnPlateau) that reduces the learning rate by a factor of 0.8 if the validation loss does not improve for 2000 epochs. Additionally, the training was terminated early if the validation loss did not improve for 20000 epochs.
For OptINNs the validation loss is the KKT-loss $\Lkkt$ evaluated on \#Validation samples of $p$ with $P^i(x) = |x|, i\in\mathcal{T}$, for QPNNs it is $\namedL{QP}$ for fixed $\gamma_g$ and $\gamma_h$.
\#Data is the number of training data points and the batch size.
\#Training is the number of $p$ samples ($N$ in the equations) used for the training loss.

\subsection*{Linear programming}
\input{tables/linear_hparam}
\vspace{-6mm}
\input{tables/linear_hparam_no_data}
\vspace{-6mm}
\input{tables/linear_hparam_penalty}
\subsection*{Nonconvex inequality constraints}
\input{tables/quadr_hparam_optinn}
\vspace{-6mm}
\input{tables/quadr_hparam_penalty}
\subsection*{Linear optimal control problem -- Rocket car}
\input{tables/rocketcar_hparam_optinn}
\vspace{-6mm}
\input{tables/rocketcar_hparam_penalty}
\subsection*{Nonlinear optimal control problem -- Pendulum swing-up}
\input{tables/pendulum_hparam_optinn}
\vspace{-6mm}
\input{tables/pendulum_hparam_penalty}


\section{Proof of~\Cref{lem:Dloss_KKT}}	\label{app:ProofLemma}
	If the conditions~\ref{eq:kkt} are satisfied, it follows that for all $i\in\{1, \dots, n_x\}$, $j\in\{1, \dots, n_h\}$, and $k\in\{1, \dots, n_g\}$
	\[
	\frac{\partial L}{\partial \xv^i} (\xv^\star, \lv^\star, \mv^\star, \pv) = \hv^j(\xv^\star, \pv) = \relu(\gv^k(\xv^\star, \pv)) = {\mv^\star}^k \gv^k(\xv^\star, \pv) = 0.
	\]
	From \Cref{def:unimodal_penalty}, since $P^i(0) = 0$, $i\in\mathcal{T}$, it follows that 
	\[\Ds(\xv^\star, \lv^\star, \mv^\star, \pv) = \Dpg(\xv^\star, \pv) = \Dph(\xv^\star, \pv) = \Dcs(\xv^\star,\mv^\star, \pv) = 0,\]
	i.e.\ KKT-points result in \MKHcorrection{zero deviations from the KKT conditions}.
	
	\noindent Moreover, assume there exists $(\tilde{\xv}, \tilde{\lv}, \tilde{\mv}) \in \mathbb{R}^{n_x} \times \mathbb{R}^{n_h} \times \mathbb{R}^{n_g}_{\geq 0}$ and $\pv \in \R^{n_p}$ not satisfying~\ref{eq:kkt}, with \[\Ds(\tilde{\xv}, \tilde{\lv}, \tilde{\mv}, \pv) = \Dpg(\tilde{\xv}, \pv) = \Dph(\tilde{\xv}, \pv) = \Dcs(\tilde{\xv},\tilde{\mv}, \pv) = 0.\]
	Since $(\tilde{\xv}, \tilde{\lv}, \tilde{\mv})$ do not satisfy~\ref{eq:kkt}, for at least one $i\in\{1, \dots, n_x\}$, $j\in\{1, \dots, n_h\}$, or $k\in\{1, \dots, n_g\}$ any of the quantities\ 
	\[
	\frac{\partial L}{\partial \xv^i} (\tilde{\xv}, \tilde{\lv}, \tilde{\mv}, \pv),\ \hv^j(\tilde{\xv}, \pv),\ \relu(\gv^k(\tilde{\xv}, \pv)),\ \tilde{\mv}^k \gv^k(\tilde{\xv}, \pv),
	\] 
	is not equal to 0.
	This contradicts the assumption that the $P^i$ are unimodal penalty functions, since $P^i(x) = 0$ only at $x=0$ and $P^i(x) > 0$ everywhere else.

%% file: tables/linear_hparam.tex
\def\file{data/results/lp/linear_problem_hparams_optinn.csv}

\pgfplotstableset{
	col sep=comma,
	columns/Data/.style={column name={\#Data}},
	columns/Training/.style={column name={\#Training}},
	columns/Validation/.style={column name={\#Validation}},
	columns/L_MIN/.style={column name={$\underline{\alpha}$}},
	columns/L_MAX/.style={column name={$\overline{\alpha}$}},	
	columns/LB_NORM/.style={column name={$\beta$}},
	columns/init. phase/.style={column name={$d^\mathrm{init}$}, fixed, precision=0, 1000 sep={\,}},
	columns/anneal. phase/.style={column name={$d^\mathrm{anneal}$}, fixed, precision=0, 1000 sep={\,}},
	columns/final. phase/.style={column name={$d^\mathrm{final}$}, fixed, precision=0, 1000 sep={\,}},
}

\begin{table}[H]
    \caption{Hyperparameters for training the OptINN on the linear programming problem.}
    \label{tab:linear_optinn_hparam}
    \centering
    
    \pgfplotstabletypeset[
        col sep=comma,
        begin table={\begin{tabularx}{\linewidth}},
        end table={\end{tabularx}},
        columns={Epochs, Data, Training, Validation, Penalty type},
        every head row/.style={before row=\hline,after row=\hline},
        every last row/.style={after row=\hline},
        every column/.style={column type=|Y},
        columns/Penalty type/.append style={column type=|Y|},
    ]{\file}
    
    \pgfplotstabletypeset[
        col sep=comma,
        begin table={\begin{tabularx}{\linewidth}},
        end table={\end{tabularx}},
        columns={Depth, Width, LR, weight decay, LB_NORM},
        every head row/.style={before row=\hline,after row=\hline},
        every last row/.style={after row=\hline},
        every column/.style={column type=|Y},
        columns/LB_NORM/.append style={column type=|Y|},
    ]{\file}
    
    \pgfplotstabletypeset[
        col sep=comma,
        begin table={\begin{tabularx}{\linewidth}},
        end table={\end{tabularx}},
        columns={L_MIN, L_MAX, init. phase, anneal. phase, final. phase},
        every head row/.style={before row=\hline,after row=\hline},
        every last row/.style={after row=\hline},
        every column/.style={column type=|Y},
        columns/final. phase/.append style={column type=|Y|},
    ]{\file}
\end{table}

%% file: tables/linear_hparam_no_data.tex
\def\file{data/results/lp/linear_problem_hparams_optinn_no_data.csv}

\pgfplotstableset{
	col sep=comma,
	columns/Data/.style={column name={\#Data}},
	columns/Training/.style={column name={\#Training}},
	columns/Validation/.style={column name={\#Validation}},
	columns/L_MIN/.style={column name={$\underline{\alpha}$}},
	columns/L_MAX/.style={column name={$\overline{\alpha}$}},	
	columns/LB_NORM/.style={column name={$\beta$}},
	columns/init. phase/.style={column name={$d^\mathrm{init}$}, fixed, precision=0, 1000 sep={\,}},
	columns/anneal. phase/.style={column name={$d^\mathrm{anneal}$}, fixed, precision=0, 1000 sep={\,}},
	columns/final. phase/.style={column name={$d^\mathrm{final}$}, fixed, precision=0, 1000 sep={\,}},
}

\begin{table}[H]
    \caption{Hyperparameters for training the OptINN without training data on the linear programming problem.}
    \label{tab:linear_optinn_nodata_hparam}
    \centering
    
    \pgfplotstabletypeset[
        col sep=comma,
        begin table={\begin{tabularx}{\linewidth}},
        end table={\end{tabularx}},
        columns={Epochs, Data, Training, Validation, Penalty type},
        every head row/.style={before row=\hline,after row=\hline},
        every last row/.style={after row=\hline},
        every column/.style={column type=|Y},
        columns/Penalty type/.append style={column type=|Y|},
    ]{\file}
    
    \pgfplotstabletypeset[
        col sep=comma,
        begin table={\begin{tabularx}{\linewidth}},
        end table={\end{tabularx}},
        columns={Depth, Width, LR, weight decay, LB_NORM},
        every head row/.style={before row=\hline,after row=\hline},
        every last row/.style={after row=\hline},
        every column/.style={column type=|Y},
        columns/LB_NORM/.append style={column type=|Y|},
    ]{\file}
    
    \pgfplotstabletypeset[
        col sep=comma,
        begin table={\begin{tabularx}{\linewidth}},
        end table={\end{tabularx}},
        columns={L_MIN, L_MAX, init. phase, anneal. phase, final. phase},
        every head row/.style={before row=\hline,after row=\hline},
        every last row/.style={after row=\hline},
        every column/.style={column type=|Y},
        columns/final. phase/.append style={column type=|Y|},
    ]{\file}
\end{table}

%% file: tables/linear_hparam_penalty.tex
\def\file{data/results/lp/linear_problem_hparams_penalty.csv}

\pgfplotstableset{
    col sep=comma,
    columns/Data/.style={column name={\#Data}},
    columns/Training/.style={column name={\#Training}},
    columns/Validation/.style={column name={\#Validation}},
    columns/alpha/.style={column name={$\alpha$}},
}

\begin{table}[H]
    \caption{Hyperparameters for training the QPNN on the linear programming problem.}
    \label{tab:linear_penalty_hparam}
    \centering
    
    \pgfplotstabletypeset[
        col sep=comma,
        begin table={\begin{tabularx}{\linewidth}},
        end table={\end{tabularx}},
        columns={Epochs, Data, Width, LR, wg},
        every head row/.style={before row=\hline,after row=\hline},
        every last row/.style={after row=\hline},
        every column/.style={column type=|Y},
        columns/wg/.style={column type=|Y|, column name={$\gamma_g$}},
    ]{\file}
    
    \pgfplotstabletypeset[
        col sep=comma,
        begin table={\begin{tabularx}{\linewidth}},
        end table={\end{tabularx}},
        columns={Training, Validation, Depth, alpha, wh},
        every head row/.style={before row=\hline,after row=\hline},
        every last row/.style={after row=\hline},
        every column/.style={column type=|Y},
        columns/wh/.style={column type=|Y|, column name={$\gamma_h$}},
    ]{\file}
\end{table}

%% file: tables/quadr_hparam_optinn.tex
\def\file{data/results/nonconvex/quadr_nonl_hparams_optinn.csv}

\pgfplotstableset{
	col sep=comma,
	columns/Data/.style={column name={\#Data}},
	columns/Training/.style={column name={\#Training}},
	columns/Validation/.style={column name={\#Validation}},
	columns/L_MIN/.style={column name={$\underline{\alpha}$}},
	columns/L_MAX/.style={column name={$\overline{\alpha}$}},	
	columns/LB_NORM/.style={column name={$\beta$}},
	columns/init. phase/.style={column name={$d^\mathrm{init}$}, fixed, precision=0, 1000 sep={\,}},
	columns/anneal. phase/.style={column name={$d^\mathrm{anneal}$}, fixed, precision=0, 1000 sep={\,}},
	columns/final. phase/.style={column name={$d^\mathrm{final}$}, fixed, precision=0, 1000 sep={\,}},
}

\begin{table}[H]
    \caption{Hyperparameters for training the OptINN on the problem with nonconvex inequality constraints.}
    \label{tab:quadr_optinn_hparam}
    \centering
    
    \pgfplotstabletypeset[
        col sep=comma,
        begin table={\begin{tabularx}{\linewidth}},
        end table={\end{tabularx}},
        columns={Epochs, Data, Training, Validation, Penalty type},
        every head row/.style={before row=\hline,after row=\hline},
        every last row/.style={after row=\hline},
        every column/.style={column type=|Y},
        columns/Penalty type/.append style={column type=|Y|},
    ]{\file}
    
    \pgfplotstabletypeset[
        col sep=comma,
        begin table={\begin{tabularx}{\linewidth}},
        end table={\end{tabularx}},
        columns={Depth, Width, LR, weight decay, LB_NORM},
        every head row/.style={before row=\hline,after row=\hline},
        every last row/.style={after row=\hline},
        every column/.style={column type=|Y},
        columns/LB_NORM/.append style={column type=|Y|},
    ]{\file}
    
    \pgfplotstabletypeset[
        col sep=comma,
        begin table={\begin{tabularx}{\linewidth}},
        end table={\end{tabularx}},
        columns={L_MIN, L_MAX, init. phase, anneal. phase, final. phase},
        every head row/.style={before row=\hline,after row=\hline},
        every last row/.style={after row=\hline},
        every column/.style={column type=|Y},
        columns/final. phase/.append style={column type=|Y|},
    ]{\file}
\end{table}

%% file: tables/quadr_hparam_penalty.tex
\pgfplotstableset{
	col sep=comma,
	columns/Data/.style={column name={\#Data}},
	columns/Training/.style={column name={\#Training}},
	columns/Validation/.style={column name={\#Validation}},
	columns/alpha/.style={column name={$\alpha$}},
}

\def\cellwidth{1.88cm}
\def\file{data/results/nonconvex/quadr_nonl_hparams_penalty.csv}

\begin{table}[H]
	\caption{Hyperparameters for training the QPNN on the problem with nonconvex inequality constraints.}
	\centering
    \pgfplotstabletypeset[
        col sep=comma,
        begin table={\begin{tabularx}{\linewidth}},
        end table={\end{tabularx}},
        columns={Epochs, Data, Width, LR, wg},
        every head row/.style={before row=\hline,after row=\hline},
        every last row/.style={after row=\hline},
        every column/.style={column type=|Y},
        columns/wg/.style={column type=|Y|, column name={$\gamma_g$}},
    ]{\file}
    
    \pgfplotstabletypeset[
        col sep=comma,
        begin table={\begin{tabularx}{\linewidth}},
        end table={\end{tabularx}},
        columns={Training, Validation, Depth, alpha, wh},
        every head row/.style={before row=\hline,after row=\hline},
        every last row/.style={after row=\hline},
        every column/.style={column type=|Y},
        columns/wh/.style={column type=|Y|, column name={$\gamma_h$}},
    ]{\file}
	\label{tab:quadr_penalty_hparam}
\end{table}

%% file: tables/rocketcar_hparam_optinn.tex
\def\file{data/results/rocket_car/rocketcar_hparams_optinn.csv}

\pgfplotstableset{
	col sep=comma,
	columns/Data/.style={column name={\#Data}},
	columns/Training/.style={column name={\#Training}},
	columns/Validation/.style={column name={\#Validation}},
	columns/L_MIN/.style={column name={$\underline{\alpha}$}},
	columns/L_MAX/.style={column name={$\overline{\alpha}$}},	
	columns/LB_NORM/.style={column name={$\beta$}},
	columns/init. phase/.style={column name={$d^\mathrm{init}$}, fixed, precision=0, 1000 sep={\,}},
	columns/anneal. phase/.style={column name={$d^\mathrm{anneal}$}, fixed, precision=0, 1000 sep={\,}},
	columns/final. phase/.style={column name={$d^\mathrm{final}$}, fixed, precision=0, 1000 sep={\,}},
}

\begin{table}[H]
    \caption{Hyperparameters for training the OptINN on the optimal control for the rocket car.}
    \label{tab:rocketcar_optinn_hparam}
    \centering
    
    \pgfplotstabletypeset[
        col sep=comma,
        begin table={\begin{tabularx}{\linewidth}},
        end table={\end{tabularx}},
        columns={Epochs, Data, Training, Validation, Penalty type},
        every head row/.style={before row=\hline,after row=\hline},
        every last row/.style={after row=\hline},
        every column/.style={column type=|Y},
        columns/Penalty type/.append style={column type=|Y|},
    ]{\file}
    
    \pgfplotstabletypeset[
        col sep=comma,
        begin table={\begin{tabularx}{\linewidth}},
        end table={\end{tabularx}},
        columns={Depth, Width, LR, weight decay, LB_NORM},
        every head row/.style={before row=\hline,after row=\hline},
        every last row/.style={after row=\hline},
        every column/.style={column type=|Y},
        columns/LB_NORM/.append style={column type=|Y|},
    ]{\file}
    
    \pgfplotstabletypeset[
        col sep=comma,
        begin table={\begin{tabularx}{\linewidth}},
        end table={\end{tabularx}},
        columns={L_MIN, L_MAX, init. phase, anneal. phase, final. phase},
        every head row/.style={before row=\hline,after row=\hline},
        every last row/.style={after row=\hline},
        every column/.style={column type=|Y},
        columns/final. phase/.append style={column type=|Y|},
    ]{\file}
\end{table}

%% file: tables/rocketcar_hparam_penalty.tex
\pgfplotstableset{
	col sep=comma,
	columns/Data/.style={column name={\#Data}},
	columns/Training/.style={column name={\#Training}},
	columns/Validation/.style={column name={\#Validation}},
	columns/alpha/.style={column name={$\alpha$}},
}

\def\cellwidth{1.88cm}
\def\file{data/results/rocket_car/rocketcar_hparams_penalty.csv}

\begin{table}[H]
	\caption{Hyperparameters for training the QPNN on the optimal control for the rocket car.}
	\centering
    \pgfplotstabletypeset[
        col sep=comma,
        begin table={\begin{tabularx}{\linewidth}},
        end table={\end{tabularx}},
        columns={Epochs, Data, Width, LR, wg},
        every head row/.style={before row=\hline,after row=\hline},
        every last row/.style={after row=\hline},
        every column/.style={column type=|Y},
        columns/wg/.style={column type=|Y|, column name={$\gamma_g$}},
    ]{\file}
    
    \pgfplotstabletypeset[
        col sep=comma,
        begin table={\begin{tabularx}{\linewidth}},
        end table={\end{tabularx}},
        columns={Training, Validation, Depth, alpha, wh},
        every head row/.style={before row=\hline,after row=\hline},
        every last row/.style={after row=\hline},
        every column/.style={column type=|Y},
        columns/wh/.style={column type=|Y|, column name={$\gamma_h$}},
    ]{\file}
	\label{tab:rocketcar_penalty_hparam}
\end{table}

%% file: tables/pendulum_hparam_optinn.tex
\def\file{data/results/pendulum/hparams_optinn.csv}

\pgfplotstableset{
	col sep=comma,
	columns/Data/.style={column name={\#Data}},
	columns/Training/.style={column name={\#Training}},
	columns/Validation/.style={column name={\#Validation}},
	columns/L_MIN/.style={column name={$\underline{\alpha}$}},
	columns/L_MAX/.style={column name={$\overline{\alpha}$}},	
	columns/LB_NORM/.style={column name={$\beta$}},
	columns/init. phase/.style={column name={$d^\mathrm{init}$}, fixed, precision=0, 1000 sep={\,}},
	columns/anneal. phase/.style={column name={$d^\mathrm{anneal}$}, fixed, precision=0, 1000 sep={\,}},
	columns/final. phase/.style={column name={$d^\mathrm{final}$}, fixed, precision=0, 1000 sep={\,}},
}

\begin{table}[H]
    \caption{Hyperparameters for training the OptINN on the optimal control for the pendulum.}
    \label{tab:pendulum_hparam}
    \centering
    
    \pgfplotstabletypeset[
        col sep=comma,
        begin table={\begin{tabularx}{\linewidth}},
        end table={\end{tabularx}},
        columns={Epochs, Data, Training, Validation, Penalty type},
        every head row/.style={before row=\hline,after row=\hline},
        every last row/.style={after row=\hline},
        every column/.style={column type=|Y},
        columns/Penalty type/.append style={column type=|Y|},
    ]{\file}
    
    \pgfplotstabletypeset[
        col sep=comma,
        begin table={\begin{tabularx}{\linewidth}},
        end table={\end{tabularx}},
        columns={Depth, Width, LR, weight decay, LB_NORM},
        every head row/.style={before row=\hline,after row=\hline},
        every last row/.style={after row=\hline},
        every column/.style={column type=|Y},
        columns/LB_NORM/.append style={column type=|Y|},
    ]{\file}
    
    \pgfplotstabletypeset[
        col sep=comma,
        begin table={\begin{tabularx}{\linewidth}},
        end table={\end{tabularx}},
        columns={L_MIN, L_MAX, init. phase, anneal. phase, final. phase},
        every head row/.style={before row=\hline,after row=\hline},
        every last row/.style={after row=\hline},
        every column/.style={column type=|Y},
        columns/final. phase/.append style={column type=|Y|},
    ]{\file}
\end{table}

%% file: tables/pendulum_hparam_penalty.tex
\pgfplotstableset{
	col sep=comma,
	columns/Data/.style={column name={\#Data}},
	columns/Training/.style={column name={\#Training}},
	columns/Validation/.style={column name={\#Validation}},
	columns/alpha/.style={column name={$\alpha$}},
}

\def\cellwidth{1.88cm}
\def\file{data/results/pendulum/hparams_penalty.csv}

\begin{table}[H]
	\caption{Hyperparameters for training the QPNN on the optimal control for the pendulum.}
	\centering
    \pgfplotstabletypeset[
        col sep=comma,
        begin table={\begin{tabularx}{\linewidth}},
        end table={\end{tabularx}},
        columns={Epochs, Data, Width, LR, wg},
        every head row/.style={before row=\hline,after row=\hline},
        every last row/.style={after row=\hline},
        every column/.style={column type=|Y},
        columns/wg/.style={column type=|Y|, column name={$\gamma_g$}},
    ]{\file}
    
    \pgfplotstabletypeset[
        col sep=comma,
        begin table={\begin{tabularx}{\linewidth}},
        end table={\end{tabularx}},
        columns={Training, Validation, Depth, alpha, wh},
        every head row/.style={before row=\hline,after row=\hline},
        every last row/.style={after row=\hline},
        every column/.style={column type=|Y},
        columns/wh/.style={column type=|Y|, column name={$\gamma_h$}},
    ]{\file}
	\label{tab:pendulum_penalty_hparam}
\end{table}